\documentclass[a4paper,12pt]{article}
\usepackage{amssymb,amsmath,amsthm}
\usepackage{indentfirst}
\usepackage{hyperref,hypbmsec}
\usepackage[pdftex]{graphicx}

\numberwithin{equation}{section}

\newtheorem{Th}{\hskip\parindent Theorem}[section]
\newtheorem{Le}{\hskip\parindent Lemma}[section]
\newtheorem{Sl}{\hskip\parindent Corollary}[section]
\newtheorem{Zam}{\hskip\parindent Remark}[section]
\newtheorem{Hyp}{\hskip\parindent Conjecture}[section]

\newcommand{\A}{\mathcal{A}}
\newcommand{\E}{\mathfrak{C}}
\newcommand{\R}{\mathfrak{R}}
\newcommand{\s}{\mathbb{S}}
\newcommand{\D}{\mathfrak{D}}
\newcommand{\N}{\mathbb{N}}
\newcommand{\M}{\mathfrak{M}}
\newcommand{\MM}{\mathfrak{M}^{*}}
\newcommand{\NN}{\mathfrak{N}}
\newcommand{\NNN}{\mathfrak{N}^{*}}
\newcommand{\1}{\mathbf{1}}
\newcommand{\jj}{\hat{j}}
\newcommand{\h}{\hat{h}}
\newcommand{\Om}{\widetilde{\Omega}}
\newcommand{\q}{\mathbf{q}}
\newcommand{\KK}{\overline{K}}
\newcommand{\Ll}{\lambda}
\newcommand{\LL}{\Lambda}
\newcommand{\K}{|K|}
\newcommand{\QN}{N_{1/2}}
\newcommand{\KN}{N^{*}_{1/2}}
\newcommand{\Z}{\mathbb{Z}}
\newcommand{\rr}{\mathbb{R}}
\newcounter{propet}

\renewcommand{\le}{\leqslant}\renewcommand{\ge}{\geqslant}
\renewcommand{\proofname}{Proof}
\renewcommand{\abstractname}{{\bf Abstract}}

\begin{document}
\author{ D.\,A.\,Frolenkov\footnote{
Research is supported by the Dynasty Foundation, by the Russian Foundation for Basic Research
(grants no. 11-01-00759-a and no. 12-01-31165).} \quad
I.\,D.\,Kan,\footnote{Research is supported by RFFI  (grant № 12-01-00681-а)} }

\title{
\begin{flushright}
\small{Dedicated to the memory of\\professor N.M. Korobov.}
\end{flushright}
A reinforcement of the Bourgain-Kontorovich's theorem by elementary methods II.
}
\date{}
\maketitle

\begin{abstract}
Zaremba's conjecture (1971) states that every positive integer number $d$ can be represented as a denominator (continuant) of a finite continued fraction  $\frac{b}{d}=[d_1,d_2,\ldots,d_{k}],$ with all partial quotients $d_1,d_2,\ldots,d_{k}$ being bounded by an absolute constant $A.$ Recently (in 2011) several new theorems concerning this conjecture were proved by Bourgain and Kontorovich. The easiest of them states that the set of numbers satisfying Zaremba's conjecture with $A=50$ has positive proportion in $\N.$ In this paper,using only elementary methods, the same theorem is proved with $A=5.$\\
\noindent
\textbf{Bibliography:} 15 titles. \par

\end{abstract}
\renewcommand{\abstractname}{{\bf Abstract}}
\renewcommand{\proofname}{{\bf Proof}}
\renewcommand{\refname}{{\bf References}}
\renewcommand{\contentsname}{{\bf Contents}}
\setcounter{Zam}0
\tableofcontents
\part{Introduction}
\section{Historical background}
Let $\R_{A}$ be the set of rational numbers whose continued fraction expansion has all partial quotients being bounded by  $A:$
$$
\R_{A}=\left\{\frac{b}{d}=[d_1,d_2,\ldots,d_{k}]\Bigl| 1\le d_{j}\le A \,\mbox{for}\, j=1,\ldots,k\right\},
$$
where
\begin{equation}\label{cont.fraction}
[d_1,\ldots,d_k]=\cfrac{1}{d_1+{\atop\ddots\,\displaystyle{+\cfrac{1}{d_k}}}}
\end{equation}
is a finite continued fraction, $d_1,\ldots,d_{k}$ are partial quotients. Let $\D_{A}$ be the set of denominators of numbers in $\R_{A}$:
$$
\D_{A}=\left\{d\Bigl| \exists b:\, (b,d)=1, \frac{b}{d}\in\R_{A}\right\}.
$$
\begin{Hyp}(\textbf{Zaremba's conjecture} ~\cite[p. 76]{Zaremba},\,1971 ).
For sufficiently large $A$  one has
$$\D_{A}=\N.$$
\end{Hyp}
That is, every $d\ge1$ can be represented as a denominator of a finite continued fraction $\frac{b}{d}$ whose partial quotients are bounded by $A.$ In fact Zaremba conjectured that $A=5$ is already large enough. A bit earlier, in the 1950-th, studying problems concerning  numerical integration Bahvalov, Chensov and N.M. Korobov also made the same assumption. But they did not publish it anywhere. Korobov ~\cite{Korobov2} proved that for a prime $p$ there exists $a,$ such that the  greatest partial quotient of $\frac{a}{p}$ is smaller than $\log p.$ A detailed survey on results concerning Zaremba's conjecture can be found in ~\cite{BK},~\cite{NG}.\par
Bourgain and Kontorovich suggested that the problem should be generalized in the following way. Let $\A \in \N$ be any finite alphabet ($|\A|\ge2$) and let $\R_{\A}$ and $\E_{\A}$  be the set of finite and infinite continued fractions
whose partial quotients belong to $\A:$
$$\R_{\A}=\left\{[d_1,\ldots,d_k]: d_j\in\A,\,j=1,\ldots,k\right\},$$
$$\E_{\A}=\left\{[d_1,\ldots,d_j,\ldots]: d_j\in\A,\,j=1,\ldots\right\}.$$
And let
$$\D_{\A}(N)=\left\{d\Bigl| d\le N,\, \exists b: (b,d)=1,\, \frac{b}{d}\in\R_{\A}\right\}$$
be the set of denominators bounded by $N.$ Let $\delta_{\A}$ be the Hausdorff dimension of $\E_{\A}.$
Then the Bourgain-Kontorovich's theorem ~\cite[p. 13, Theorem 1.25]{BK} ia as follows
\begin{Th}\label{BKtheorem}
For any alphabet $\A$  with
\begin{equation}\label{BKcondition}
\delta_{\A}>1-\frac{5}{312}=0,983914\ldots,
\end{equation}
the following inequality (positive proportion)
\begin{equation}\label{BKresult}
\#\D_{\A}(N)\gg N.
\end{equation}
holds.
\end{Th}
For some alphabets the condition \eqref{BKcondition} can be verified by two means. For an alphabet $\A=\left\{1,2,\ldots,A-1,A\right\}$ Hensley ~\cite{Hen1} proved that
\begin{equation}\label{hen-hausdorff-dimension}
\delta_{\A}=\delta_{A}=1-\frac{6}{\pi^2}\frac{1}{A}-\frac{72}{\pi^4}\frac{\log A}{A^2}+O\left(\frac{1}{A^2}\right).
\end{equation}
Moreover Jenkinson ~\cite{Jenkinson} obtained approximate values for some $\delta_{\A}.$ In view of these results the alphabet
\begin{equation}\label{specalpha}
\left\{1,2,\ldots,A-1,A\right\}
\end{equation}
with $A=50$ is assumed to satisfy \eqref{BKcondition}. Several results improving \eqref{BKresult} were also proved in ~\cite{BK}.
However, we do not consider them in our work.
\section{Statement of the main result}
The main result of the paper is the following theorem.
\begin{Th}\label{uslov}
For any alphabet $\A$ with
\begin{equation}\label{KFcondition1}
\delta_{\A}>1-\frac{1}{6}=0,8333\ldots,
\end{equation}
the following inequality (positive proportion) holds
\begin{equation}\label{KFresult}
\#\D_{\A}(N)\gg N.
\end{equation}
\end{Th}
\begin{Zam}
It is proved ~\cite{Jenkinson} that $\delta_{5}=0,8368\ldots,$ From this follows that the alphabet $\left\{1,2,\ldots,5\right\}$
satisfies the condition of Theorem \ref{uslov}.
\end{Zam}
\begin{Zam}
The proof of the Theorem \ref{uslov} is  based on the method of Bourgain-Kontorovich~\cite{BK}. In this article a new technique  for estimating exponential sums taking over a set of continuants is devised. We improve on this method by refining the main set $\Omega_N.$ In ~\cite{BK} this set is named as ensemble.
\end{Zam}

\section{Acknowledgment}
We thank prof. N.G. Moshchevitin for numerous discussions of our results. It was he who took our notice of ~\cite{BK}. We are also grateful to prof. I.D. Shkredov and prof. I.S. Rezvyakova for their questions and comments during our talks.
\section{Notation}
Throughout $\epsilon_0=\epsilon_0(\A)\in(0,\frac{1}{2500}).$ For two functions $f(x), g(x)$ the Vinogradov notation $f(x)\ll g(x)$  means that there exists a constant $C,$ depending on $A,\epsilon_0,$ such that $|f(x)|\le Cg(x).$  The notation $f(x)=O(g(x))$ means the same. The notation $f(x)=O_1(g(x))$ means that $|f(x)|\le g(x).$ The notation $f(x)\asymp g(x)$ means that
$f(x)\ll g(x)\ll f(x).$ Also a traditional notation $e(x)=\exp(2\pi ix)$ is used. The cardinality of a finite set $S$ is denoted either $|S|$ or $\#S.$ $[\alpha]$ and $\|\alpha\|$ denote the integral part of $\alpha$ and the distance from $\alpha$ to the nearest integer:
$[\alpha]=\max\left\{z\in\mathbb{Z}|\,z\le \alpha \right\},$
$\|\alpha\|=\min\left\{|z-\alpha|\Bigl|\,z\in\mathbb{Z}\right\}. $
The following sum $\sum_{d|q}1$ is denoted as $\tau(q).$
\part{Preparation for estimating exponential sums}
\section{Continuants and matrices}
In this section we recall the simplest techniques concerning continuants. As a rule, all of them can be found in any research dealing with continued fractions. To begin with we define several operations on finite sequences. Let
\begin{equation}\label{def of seq}
D=\{d_1,d_2,\ldots,d_k\},\, B=\{b_1,b_2,\ldots,b_n\}
\end{equation}
then $D,B$ denotes the following sequence
$$D,B=\{d_1,d_2,\ldots,d_k,b_1,b_2,\ldots,b_n\}.$$
For every $D$ from \eqref{def of seq} let define $D_{-},D^{-},\overleftarrow{D}$as follows
\begin{equation*}
D_{-}=\{d_2,d_3,\ldots,d_k\},\,
D^{-}=\{d_1,d_2,\ldots,d_{k-1}\},\,
\overleftarrow{D}=\{d_k,d_{k-1},\ldots,d_2,d_1\}.
\end{equation*}
We denote by $[D]$ the continued fraction \eqref{cont.fraction}, that is $[D]=[d_1,\ldots,d_k].$ And by $\langle D\rangle$ we denote its denominator $\langle D\rangle=\langle d_1,\ldots,d_k\rangle.$ This denominator is called the continuant of the sequence $D.$ The continuant of the sequence can also be defined as follows
\begin{equation*}
\langle\,\rangle=1,\,\langle d_1\rangle=d_1,
\end{equation*}
\begin{equation*}
\langle d_1,\ldots,d_k\rangle=\langle d_1,\ldots,d_{k-1}\rangle d_k+\langle d_1,\ldots,d_{k-2}\rangle,\,\mbox{for}\, k\ge2.
\end{equation*}
It is well known ~\cite{knut} that
\begin{equation}\label{continuant properties}
\langle D\rangle=\langle\overleftarrow{D}\rangle,\,
[D]=\frac{\langle D_{-}\rangle}{\langle D\rangle},\,
[\overleftarrow{D}]=\frac{\langle D^{-}\rangle}{\langle D\rangle},
\end{equation}
\begin{equation}\label{continuant identity}
\langle D,B\rangle=\langle D\rangle\langle B\rangle(1+[\overleftarrow{D}][B]).
\end{equation}
It follows from this that
\begin{equation}\label{continuant inequality}
\langle D\rangle\langle B\rangle\le
\langle D,B\rangle\le
2\langle D\rangle\langle B\rangle,
\end{equation}
and that the elements of the matrix
\begin{equation}\label{matrix}
\gamma=
\begin{pmatrix}
a & b \\
c & d
\end{pmatrix}=
\begin{pmatrix}
0 & 1 \\
1 & d_1
\end{pmatrix}
\begin{pmatrix}
0 & 1 \\
1 & d_2
\end{pmatrix}\ldots
\begin{pmatrix}
0 & 1 \\
1 & d_k
\end{pmatrix}
\end{equation}
can be expressed by continuants
\begin{gather}\label{matrix=continuant}
a=\langle d_2,d_3,\ldots,d_{k-1}\rangle,\,
b=\langle d_2,d_3,\ldots,d_{k}\rangle,\,
c=\langle d_1,d_2,\ldots,d_{k-1}\rangle,\,
d=\langle d_1,d_2,\ldots,d_{k}\rangle.
\end{gather}
For the matrix $\gamma$ from \eqref{matrix} we use the following norm
\begin{equation*}
\|\gamma\|=\max\{|a|,|b|,|c|,|d|\},
\end{equation*}
It follows from \eqref{matrix=continuant} that
\begin{equation}\label{matrix norm}
\|\gamma\|=d=\langle d_1,d_2,\ldots,d_{k}\rangle.
\end{equation}
For $\gamma$ from \eqref{matrix} we have $|\det \gamma|=1.$
Let $\Gamma_{\A}\subseteq SL\left(2,\mathbb{Z}\right)$ be a semigroup generated by
$
\begin{pmatrix}
0 & 1 \\
1 & a_i
\end{pmatrix}
$
where $a_i\in\A.$ It follows from \eqref{matrix norm} that to prove a positive density of continuants in $\N$
\begin{equation}\label{BKresult1,1}
\#\D_{\A}(N)\gg N,
\end{equation}
it is enough to obtain the same property of the set
\begin{equation}\label{set matrix norm}
\|\Gamma_{\A}\|=\left\{\|\gamma\|\Bigl| \gamma\in\Gamma_{\A}\right\}.
\end{equation}
In fact, for proving inequality  \eqref{BKresult1,1} only a part of the semigroup $\Gamma_{\A},$ a so called ensemble $\Omega_N,$ will be used. A preparation for constructing $\Omega_N$ will start in the next section.
\section{Hensley's method}\label{method of hensley}
Before estimating the amount of continuants not exceeding $N$ it might be well to assess the amount of continued fractions with denominator being bounded by $N.$ Though these problems are similar, in the second case every continuant should be counted in view of its multiplicity. Let $F_{\A}(x)=\D_{\A}([x])$ and
\begin{gather}\label{maxA}
A=\max\A.
\end{gather}
Generalizing Hensley's method~\cite{Hen2} one can prove that the following theorem holds.
\begin{Th}\label{Kan-Hen}
Let $\delta_{\A}>\frac{1}{2},$ then for any $x\ge4A^2$ one has
\begin{gather}\label{F_A less}
\frac{1}{32A^4}x^{2\delta_{\A}}\le F_{\A}(x)-F_{\A}\left(\frac{x}{4A^2}\right)\le F_{\A}(x)\le8 x^{2\delta_{\A}}.
\end{gather}
\end{Th}
Hensley~\cite{Hen2} proved this theorem for the alphabet $\A$ of the form \eqref{specalpha}.
\section{The basic ideas of the Bourgain-Kontorovich's method}\label{section idea BK}
In this section a notion of constructions necessary for proving Theorem  \ref{uslov} will be given.\par
In view of exponential sums, for studying the density of the set \eqref{set matrix norm} it is natural to estimate the absolute value of the sum
\begin{gather}\label{def Sn 1}
S_N(\theta):=\sum_{\gamma\in\Omega_{N}}e(\theta\|\gamma\|),
\end{gather}
where $\Omega_{N}\subseteq\Gamma_{\A}\cap\{\|\gamma\|\le N\}$ is a proper set of matrices (ensemble), $\theta\in[0,1],$ and the norm $\|\gamma\|$ is defined in \eqref{matrix norm}. As usual, the Fourier coefficient of the function $S_N(\theta)$ is defined by
\begin{gather*}
\widehat{S}_N(n)=\int_0^1 S_N(\theta)e(-n\theta)d\theta=\sum_{\gamma\in \Omega_{N}}\1_{\{\|\gamma\|=n\}}.
\end{gather*}
Note that if ${\widehat{S}_N(n)>0}$ then $n\in\D_{\A^2}(N).$ Since
\begin{gather*}
S_N(0)=\sum_{\gamma\in\Omega_{N}}1=\sum_{n=1}^N \sum_{\gamma\in \Omega_{N}}\1_{\{\|\gamma\|=n\}}=
\sum_{n=1}^N \widehat{S}_N(n)=\sum_{n=1}^N \widehat{S}_N(n)\1_{\{\widehat{S}_N(n)\neq0\}},
\end{gather*}
then applying Cauchy-Schwarz inequality one has
\begin{gather}\label{8-4}
(S_N(0))^2\le \sum_{n=1}^N \1_{\{\widehat{S}_N(n)\neq0\}}\sum_{m=1}^N \left(\widehat{S}_N(m)\right)^2.
\end{gather}
The first factor of the right hand side of the inequality \eqref{8-4} is bounded from above by $\#\D_{\A^2}(N).$ Applying Parseval for the second factor one has
\begin{gather*}
\sum_{n=1}^N \left(\widehat{S}_N(n)\right)^2=\int_0^1\left|S_N(\theta)\right|^2d\theta.
\end{gather*}
Consequently
\begin{gather}\label{8-5}
(S_N(0))^2\le\#\D_{\A}(N)\int_0^1\left|S_N(\theta)\right|^2d\theta.
\end{gather}
Thus a lower bound on the magnitude of the set $\D_{\A^2}(N)$ follows from \eqref{8-5}
\begin{gather}\label{8-6}
\#\D_{\A}(N)\ge\frac{(S_N(0))^2}{\int\limits_0^1\left|S_N(\theta)\right|^2d\theta}.
\end{gather}
Thus, the estimate
\begin{gather*}
\#\D_{\A}(N)\gg N
\end{gather*}
will be proved, if we are able to assess exactly as possible the integral from \eqref{8-6}
\begin{gather}\label{8-7}
\int_0^1\left|S_N(\theta)\right|^2d\theta\ll\frac{(S_N(0))^2}{N}=\frac{|\Omega_{N}|^2}{N}.
\end{gather}
It follows from the Dirichlet's theorem that for any $\theta\in[0,1]$ there exist $a,q\in\N\cup\{0\}$ and $\beta\in\rr$
such that
\begin{gather*}
\theta=\frac{a}{q}+\beta,\;(a,q)=1,\; 0\le a\le q\le N^{1/2},\;\beta=\frac{K}{N},\; |K|\le\frac{N^{1/2}}{q},
\end{gather*}
with $a=0$ and $a=q$ being possible if only $q=1.$
Following the article ~\cite{BK}, to obtain the estimate \eqref{8-7} we represent the integral as the sum of integrals over different domains in variables $(q,K).$  Each of them will be estimated in a special way depending on the domain.\par
It remains to define ensemble $\Omega_{N}.$  To begin with we determine a concept "pre-ensemble". The subset $\Xi$ of matrices $\gamma\in\Gamma_{\A}$ is referred to as $N$~--\textbf{pre-ensemble}, if the following conditions hold
\begin{enumerate}
  \item for any matrix $\gamma\in\Xi$ its norm is of the order of $N:$
\begin{gather}\label{8-18a}
\|\gamma\|\asymp N;
\end{gather}
  \item for any $\epsilon>0$ the set $\Xi$ contains $\epsilon$~--full amount of elements, that is
\begin{gather}\label{8-18b}
\#\Xi\gg_{\epsilon} N^{2\delta_{\A}-\epsilon}.
\end{gather}
\end{enumerate}
By the \textbf{product} of two pre-ensembles $\Xi^{(1)}\Xi^{(2)}$ we mean the set of all possible  products of matrices $\gamma_1\gamma_2$ such that $\gamma_1\in\Xi^{(1)},\,\gamma_2\in\Xi^{(2)}.$ The product of pre-ensembles has an \textbf{unique} \textbf{expansion} if it follows from the relations
\begin{gather}\label{8-19}
\gamma_1\gamma_2=\gamma_1'\gamma_2', \quad \gamma_1,\gamma_1'\in\Xi^{(1)},\,\gamma_2,\gamma_2'\in\Xi^{(2)}
\end{gather}
that
\begin{gather}\label{8-20}
\gamma_1=\gamma_1',\quad \gamma_2=\gamma_2'.
\end{gather}
Let $\epsilon_0$ be a fixed number such that $0<\epsilon_0\le\frac{1}{2}.$ Then $N$~-- pre-ensemble $\Omega$ is called the right (left) $(\epsilon_0,N)$~-- \textbf{ensemble} if for any $M,$ such that
\begin{gather}\label{8-21}
1\ll M\le N^{\frac{1}{2}},
\end{gather}
there exist positive numbers $N_1$ and $N_2$ such that
\begin{gather}\label{8-22}
N_1N_2\asymp N,\quad N_2^{1-\epsilon_0}\ll M\ll N_2,
\end{gather}
and there exist $N_1$~--pre-ensemble $\Xi^{(1)}$ and $N_2$~--pre-ensemble $\Xi^{(2)}$ such that the pre-ensemble $\Omega$ is equal to the product $\Xi^{(1)}\Xi^{(2)}$ ($\Xi^{(2)}\Xi^{(1)}$ respectively) having an unique expansion. Such terminology allows us to say that in the article ~\cite{BK} the $(\frac{1}{2},N)$~--ensemble has been constructed while we will construct $(\epsilon_0,N)$~--ensemble, being simultaneously the right and the left (bilateral) ensemble, for $\epsilon_0\in\left(0,\frac{1}{2500}\right).$
\begin{Zam}\label{zam813}
Observe that there is no use to require an upper bound in \eqref{8-18b}. According to the Theorem \ref{Kan-Hen} it follows from \eqref{8-18a} that
\begin{gather*}
\#\Xi\ll N^{2\delta_{\A}}.
\end{gather*}
\end{Zam}
\section{Pre-ensemble $\Xi(M).$}\label{predansambl}
Let $\delta:=\delta_{\A}>\frac{1}{2},$ $\Gamma:=\Gamma_{\A}$ and as usual $A=\max\A.$
Let also $M$ be a fixed parameter satisfying the inequality
\begin{gather}\label{9-1}
M\ge2^9A^3\log^3M.
\end{gather}
In this section we construct a pre-ensemble $\Xi(M)\subset\Gamma$ It is the key element which will be used to construct the ensemble $\Omega_N.$ To generate $\Xi(M)$ we use an algorithm. The number $M$ is an input parameter. During the algorithm we generate the following numbers
\begin{gather*}
L=L(M)\asymp M,\,p=p(M)\asymp\log\log M,\, k=k(M)\asymp\log M,
\end{gather*}
being responsible for the properties of the elements of $\Xi(M).$ We now proceed to the description of the algorithm consisting of four steps.
\begin{description}
  \item[Step 1]
First consider the set $S_1\subset\Gamma$ of matrices $\gamma\in\Gamma,$ such that
\begin{gather}\label{9-2}
\frac{M}{64A^2}\le\|\gamma\|\le M.
\end{gather}
According to the Theorem \ref{Kan-Hen}\, $\#S_1\asymp M^{2\delta}.$
  \item[Step 2]
Let $F_0=0,\,F_1=1,\,F_{n+1}=F_{n}+F_{n-1}$ for $n\ge1$ be Fibonacci numbers.  Define an integer number $p=p(M)$ by the relation
\begin{gather}\label{9-4}
F_{p-1}\le\log^{\frac{1}{2}}M\le F_p.
\end{gather}
Note that then
\begin{gather}\label{9-5}
F_{p}\le2F_{p-1}\le2\log^{\frac{1}{2}}M.
\end{gather}
Let consider the set $S_2\subset S_1$ of matrices $\gamma\in S_1$ of the form \eqref{matrix} for witch
\begin{gather}\label{9-6}
d_1=d_2=\ldots=d_p=1,\quad
d_k=d_{k-1}=\ldots=d_{k-p+1}=1.
\end{gather}
To put it another way, the first $p$ and the last $p$ elements of the sequence $D=D(\gamma)=\{d_1,d_2,\ldots,d_k\}$ are equal to one. At the moment we have to interrupt for a while the description of the algorithm in order to prove the following lemma.
\begin{Le}\label{lemma9-1}
One has the estimate
\begin{gather}\label{9-8}
\#S_2\ge\frac{M^{2\delta}}{2^{13}A^4\log^2M}.
\end{gather}
\begin{proof}
The sequence $D$ can be represented in the form
Последовательность $D$  представим как
\begin{gather*}
D=\{\underbrace{1,1,\ldots,1}_p,b_1,b_2,\ldots,b_n,\underbrace{1,1,\ldots,1}_p\}
\end{gather*}
then the required inequality \eqref{9-2} can be represented in the form
\begin{gather}\label{9-9}
\frac{M}{64A^2}\le\langle\underbrace{1,1,\ldots,1}_p,b_1,b_2,\ldots,b_n,\underbrace{1,1,\ldots,1}_p\rangle\le M.
\end{gather}
Let prove that the inequality \eqref{9-9} follows from the inequality
\begin{gather}\label{9-10}
\frac{M}{64A^2\log M}\le\langle b_1,b_2,\ldots,b_n\rangle\le \frac{M}{16\log M}.
\end{gather}
Indeed, let the inequality \eqref{9-10} be true. Then on the one side it follows from inequalities \eqref{continuant inequality} and \eqref{9-5} that
\begin{gather*}
\langle\underbrace{1,1,\ldots,1}_p,b_1,b_2,\ldots,b_n,\underbrace{1,1,\ldots,1}_p\rangle\le4F_p^2
\langle b_1,b_2,\ldots,b_n\rangle\le16\langle b_1,b_2,\ldots,b_n\rangle\log M\le M,
\end{gather*}
and on the other side it follows in a similar way from the inequality \eqref{9-4} that
\begin{gather*}
\langle\underbrace{1,1,\ldots,1}_p,b_1,b_2,\ldots,b_n,\underbrace{1,1,\ldots,1}_p\rangle\ge F_p^2
\langle b_1,b_2,\ldots,b_n\rangle\ge\langle b_1,b_2,\ldots,b_n\rangle\log M\ge\frac{M}{64A^2}.
\end{gather*}
Thus the implication $\eqref{9-10}\Rightarrow\eqref{9-9}$ is proved. It remains to obtain a lower bound for the amount of sequences $B=\langle b_1,b_2,\ldots,b_n\rangle$ satisfying the inequality \eqref{9-10}. We set $x=\frac{M}{16\log M}$ and note that the condition $x\ge4A^2$ in Theorem \ref{Kan-Hen} follows from the inequality \eqref{9-1}. Consequently, considering \eqref{F_A less} one has
\begin{gather*}
\#S_2\ge F_{\A}(x)-F_{\A}\left(\frac{x}{4A^2}\right)\ge\frac{1}{32A^4}\left(\frac{M}{16\log M}\right)^{2\delta}\ge
\frac{M^{2\delta}}{2^{13}A^4\log^2M},
\end{gather*}
since $\delta<1.$ This completes the proof of the lemma.
\end{proof}
\end{Le}
Let return to the description of the algorithm.
  \item[Step 3]
For any $L$ in the interval
\begin{gather}\label{9-11}
\left[\frac{M}{64A^2}, M\right]
\end{gather}
consider the set $S_3(L)\subset S_2$ of matrices $\gamma\in S_2,$ for which the following inequality holds
\begin{gather}\label{9-12}
\max\left\{\frac{M}{64A^2}, L(1-\log^{-1}L)\right\}\le\|\gamma\|\le L.
\end{gather}
Here, we also have to interrupt the description of the algorithm in order to prove the following lemma.
\begin{Le}\label{lemma9-2}
There is a number $L$ in the interval \eqref{9-11} such that
\begin{gather}\label{9-13}
|S_3(L)|\ge\frac{L^{2\delta}}{2^{16}A^5\log^3L}.
\end{gather}
\begin{proof}
Let $t$ be the minimal positive integer number satisfying the inequality
\begin{gather}\label{9-14}
(1-\log^{-1}M)^t\le\frac{1}{64A^2}.
\end{gather}
Note that $t\le8A\log M.$ For $j=1,2,\ldots,t$ consider sets $s(j)$ each of them consists of matrices  $\gamma\in S_2,$ such that
\begin{gather*}
M(1-\log^{-1}M)^j\le\|\gamma\|\le M(1-\log^{-1}M)^{j-1}
\end{gather*}
Since $S_2\subset\bigcup_{1\le j\le t}s(j),$ by the pigeonhole principle there is a set among $s(j)$ containing at least $\frac{|S_2|}{t}$ matrices. Let $$L=M(1-\log^{-1}M)^{j_0-1},$$ then $L$ belongs to the segment \eqref{9-11} and $s(j_0)\subset S_3(L).$ Hence $|S_3(L)|\ge\frac{|S_2|}{t}.$ Using the bound \eqref{9-8} and the restriction on $t$ one has
\begin{gather}\label{9-17}
|S_3(L)|\ge\frac{M^{2\delta}}{(2^{13}A^4\log^2M)8A\log M}=\frac{M^{2\delta}}{2^{16}A^5\log^3M}.
\end{gather}
Because the function $f(x)=x^{2\delta}\log^{-3}x$ increases and since $M\ge L,$ then replacing in \eqref{9-17} the parameter $M$ by $L$ one has the inequality \eqref{9-13}. This completes the proof of the lemma.
\end{proof}
\end{Le}
Returning to the algorithm we choose in the interval \eqref{9-11} any $L$ (for example the maximal one) satisfying the inequality \eqref{9-13} and fix it. Now let $S_3:=S_3(L).$
  \item[Step 4]
For $\gamma\in S_3$ let $k(\gamma)$ be the length of the sequence $D(\gamma).$ Represent the set $S_3$ as the union of the sets $S_4(k),$ consisting of those matrices $\gamma\in S_3$ for which $k(\gamma)=k$ is fixed.
\begin{Le}\label{lemma9-3}
There exists $k$ for which
\begin{gather}\label{9-18}
|S_4(k)|\ge\frac{L^{2\delta}}{2^{18}A^5\log^4L}.
\end{gather}
\begin{proof}
Since for all $D\in V_{\A}(r)$ the inequality
$
\langle D\rangle\ge\langle \,\underbrace{1,1,\ldots,1}_{r}\,\rangle,
$
holds, then
\begin{equation*}
\langle D\rangle\ge\left(\frac{\sqrt{5}+1}{2}\right)^{r-1}
\end{equation*}
and consequently
\begin{gather}\label{9-19}
k\le\frac{\log\|\gamma\|}{\log\frac{\sqrt{5}+1}{2}}+1\le4\log\|\gamma\|\le4\log L.
\end{gather}
Hence, by the pigeonhole principle, there is a $k,$ for which
\begin{gather*}
|S_4(k)|\ge\frac{|S_3|}{4\log L}\ge\frac{L^{2\delta}}{(4\log L)2^{16}A^5\log^3L}=\frac{L^{2\delta}}{2^{18}A^5\log^4L}
\end{gather*}
by \eqref{9-13} and \eqref{9-19}. This completes the proof of the lemma.
\end{proof}
\end{Le}
Returning to the algorithm we fix any $k,$ satisfying the inequality \eqref{9-18} and write $S_4:=S_4(k).$
\textbf{Algorithm is completed.}
\end{description}
Now we write $\Xi(M):=S_4.$ Recall the properties of matrices $\gamma\in\Xi(M).$ For any $\gamma\in\Xi(M)$ we have from the construction:
\begin{description}
  \item[i] the first and the last $p$ elements of the sequence $D(\gamma)$ are equal to 1, where $p$ is defined by the inequality \eqref{9-4};
  \item[ii] $L(1-\log^{-1}L)\le\|\gamma\|\le L;$
  \item[iii] $k(\gamma)=const,$ that is, the length of $D(\gamma)$ is fixed fir all $\gamma\in\Xi(M).$
\end{description}
Besides, we have proved that
\begin{gather}\label{9-20}
\#\Xi(M)\ge\frac{L^{2\delta}}{2^{18}A^5\log^4L}.
\end{gather}
The first property allows us to prove an important lemma
\begin{Le}\label{pred9-1}
For every matrix $\gamma\in\Xi(M)$ of the form
$\gamma=
\begin{pmatrix}
a & b \\
c & d
\end{pmatrix}
$
the following inequalities hold
\begin{gather}\label{9-21}
\left|\frac{b}{d}-\frac{1}{\varphi}\right|\le\frac{2}{\log L},\quad
\left|\frac{c}{d}-\frac{1}{\varphi}\right|\le\frac{2}{\log L},
\end{gather}
where $\varphi$ is the golden ratio
\begin{gather}\label{9-22}
\varphi=1+[1,1,\ldots,1,\ldots]=\frac{\sqrt{5}+1}{2}.
\end{gather}
\begin{proof}
It follows from \eqref{continuant properties} and \eqref{matrix=continuant} that
$\frac{b}{d}=[D(\gamma)],\quad \frac{c}{d}=[\overleftarrow{D(\gamma)}].$ Hence the fraction
$\alpha=[\underbrace{1,1,\ldots,1}_p]$
is a convergent fraction to $\frac{b}{d}$ and to $\frac{c}{d}.$ The denominator of $\alpha$ is equal to $F_p$ and
it follows from the choice of parameters \eqref{9-4}, \eqref{9-11} that
$F_p\ge\log^{\frac{1}{2}}L.$
Hence,
\begin{gather}\label{9-24}
\left|\frac{b}{d}-\alpha\right|\le\frac{1}{\log L},\quad
\left|\frac{c}{d}-\alpha\right|\le\frac{1}{\log L}.
\end{gather}
But $\alpha$ is also a convergent fraction to $\frac{1}{\varphi},$ thus
$\left|\alpha-\frac{1}{\varphi}\right|\le\frac{1}{\log L}.$
Applying the triangle inequality we obtain the desired inequalities. This completes the proof of the lemma.
\end{proof}
\end{Le}
\section{Parameters and their properties}
Let $N\ge N_{min}=N_{min}(\epsilon_0,\A)$ and write
\begin{gather}\label{10-2}
J=J(N)=\left[\frac{\log\log N-4\log(10A)+2\log\epsilon_0}{-\log(1-\epsilon_0)}\right],
\end{gather}
where as usual $A\ge|\A|\ge2$ and require the following inequality $J(N_{min})\ge10$ to hold.
Using the definition \eqref{10-2}, one has
\begin{gather}\label{10-4}
\frac{10^4A^4}{\log N}\le\epsilon_0^2(1-\epsilon_0)^J\le\frac{10^5A^4}{\log N}.
\end{gather}
Now let define a finite sequence
\begin{gather}\label{10-6}
\left\{N_{-J-1},N_{-J},\ldots,N_{-1},N_0,N_1,\ldots,N_{J+1}\right\},
\end{gather}
having set $N_{J+1}=N$ and
\begin{gather}\label{10-8}
N_j=
\left\{
              \begin{array}{ll}
                N^{\frac{1}{2-\epsilon_0}(1-\epsilon_0)^{1-j}}, & \hbox{if $-1-J\le j\le1$;} \\
                N^{1-\frac{1}{2-\epsilon_0}(1-\epsilon_0)^{j}}, & \hbox{if $0\le j\le J$.}
              \end{array}
\right.
\end{gather}
It is obvious that the sequence is well-defined for $j=0$ and $j=1.$
\begin{Le}\label{lemma10-2}
\begin{enumerate}
  \item For $-J\le m\le J-1$ the following equation holds
  \begin{gather}\label{10-9}
N_{-m}N_{m+1}=N.
\end{gather}
  \item For $-J-1\le m\le J-1$ the following relations hold
\begin{gather}\label{10-10}
\frac{N_{m+1}}{N_m}=
\left\{
              \begin{array}{ll}
                N_{m+1}^{\epsilon_0}, & \hbox{if $m\le0$;} \\ \\
                \left(\frac{N}{N_m}\right)^{\epsilon_0}, & \hbox{if $m\ge0$,}
              \end{array}
\right.
\end{gather}
\begin{gather}\label{10-11}
\frac{N_{m+1}}{N_m}= N^{\frac{\epsilon_0}{2-\epsilon_0}(1-\epsilon_0)^{|m|}},
\end{gather}
\begin{gather}\label{10-12}
N_{m}\ge N_{m+1}^{1-\epsilon_0}.
\end{gather}
  \item Для $-1\le j<h\le J+1$ выполнено
\begin{gather}\label{10-13}
N_{h-J}^{(1-\epsilon_0)^{h-j}}=N_{j-J}.
\end{gather}
\end{enumerate}
\begin{proof}
All propositions follow directly from the definition  \eqref{10-8}. This completes the proof of the lemma.
\end{proof}
\end{Le}
\begin{Le}\label{lemma10-3}
For $-J\le m\le J-1$ the following estimate holds
\begin{gather}\label{10-16}
\frac{N_{m+1}}{N_{m}}\ge\exp\left(\frac{10^4A^4}{2\epsilon_0}\right);
\end{gather}
moreover
\begin{gather}\label{10-17}
\exp\left(\frac{10^4A^4}{2\epsilon_0^2}\right)\le\frac{N}{N_{J}}=N_{1-J}\le\exp\left(\frac{10^5A^4}{\epsilon_0^2}\right).
\end{gather}
\begin{proof}
The inequality \eqref{10-16} follows from \eqref{10-11} and the lower bound in \eqref{10-4}. Now let prove the inequality \eqref{10-17}. The equation $\frac{N}{N_{J}}=N_{1-J}$ follows from \eqref{10-9} with $m=-J.$  The estimate of $N_{1-J}$ follows from \eqref{10-8} and \eqref{10-4}. This completes the proof of the lemma.
\end{proof}
\end{Le}
\begin{Le}\label{lemma10-4}
For any $M,$ such that
\begin{gather}\label{10-19}
N_{1-J}\le M\le N_{J},
\end{gather}
there exist indexes  $j$ and $h,$ such that
\begin{gather}\label{10-20}
2\le j\le 2J,\quad 1\le h\le 2J-1,\quad h=2J-j+1,
\end{gather}
for which the following inequalities hold
\begin{gather}\label{10-21}
N_{j-J}^{1-\epsilon_0}\le M\le N_{j-J},
\end{gather}
\begin{gather}\label{10-22}
\left(\frac{N}{N_{h-J}}\right)^{1-\epsilon_0}\le M\le\frac{N}{N_{h-J}}.
\end{gather}
\begin{proof}
Since the sequence $\{N_j\}$ is increasing there exists the index $j$ in \eqref{10-20} such that
\begin{gather}\label{10-23}
N_{j-1-J}\le M\le N_{j-J},
\end{gather}
then \eqref{10-21} follows from \eqref{10-12}.
The inequality \eqref{10-22} can be obtained by substituting the equation \eqref{10-9} into \eqref{10-21}. This completes the proof of the lemma.
\end{proof}
\end{Le}
For a nonnegative integer number $n$ we write
\begin{gather}\label{10-26}
\widetilde{N}_{n-J}=
\left\{
              \begin{array}{ll}
                N_{n-J}, & \hbox{if $n\ge1$;} \\
                1, & \hbox{if $n=0.$}
              \end{array}
\right.
\end{gather}
Moreover for integers $j$ and $h$ such that
\begin{gather}\label{10-27}
0\le j<h\le 2J+1,
\end{gather}
we write
\begin{gather}\label{10-28}
j_0(j,h)=\min\left\{|n-J-1|\,\Bigl|\, j+1\le n\le h\right\}.
\end{gather}
Note that there are only three alternatives for the value of $j_0$
\begin{gather}\label{10-29}
j_0\in\left\{j-J,\,0,\,J+1-h\right\}.
\end{gather}
\begin{Le}\label{lemma10-5}
For integers $j$ and $h$ from \eqref{10-27} the following estimate holds
\begin{gather}\label{10-30}
\prod_{n=j+1}^{h}\left(2^9A^3\log\frac{N_{n-J}}{\widetilde{N}_{n-1-J}}\right)\le
\left((1-\epsilon_0)^{j_0}\log N\right)^{\frac{7}{4}(h-j)},
\end{gather}
where
$j_0=j_0(j,h).$
\begin{proof}
Consider two cases depending on the value of $j.$
\begin{description}
  \item[1)\,$j>0.$]
Using \eqref{10-16} if $n<2J+1$ and \eqref{10-17} if $n=2J+1,$ for any $n$ in the segment $j+1\le n\le h$ we obtain
\begin{gather}\label{10-31}
2^9A^3\le\frac{2^9}{10^3}(2\epsilon_0)^{3/4}\left(\log\frac{N_{n-J}}{N_{n-1-J}}\right)^{3/4}
\le\left(\log\frac{N_{n-J}}{N_{n-1-J}}\right)^{3/4}.
\end{gather}
Hence, since for $j>0$ one has $\widetilde{N}_{n-1-J}=N_{n-1-J},$ we obtain
\begin{gather}\label{10-32}
\prod_{n=j+1}^{h}\left(2^9A^3\log\frac{N_{n-J}}{\widetilde{N}_{n-1-J}}\right)\le
\prod_{n=j+1}^{h}\left(\log\frac{N_{n-J}}{N_{n-1-J}}\right)^{7/4}.
\end{gather}
Making allowance for $\epsilon_0\in\left(0,\frac{1}{2500}\right),$ it follows from \eqref{10-11} if  $n<2J+1$
and from \eqref{10-8} if $n=2J+1,$ that
\begin{gather}\label{10-33}
\log\frac{N_{n-J}}{N_{n-1-J}}\le
\frac{1}{2-\epsilon_0}(1-\epsilon_0)^{|n-J-1|}\log N\le(1-\epsilon_0)^{j_0}\log N
\end{gather}
Substituting the estimate \eqref{10-33} into \eqref{10-32} we obtain the statement of the theorem in case $j>0.$
  \item[2)\,$j=0.$]
Using the lower bound from \eqref{10-17} we obtain
\begin{gather}\label{10-34}
2^9A^3\le\frac{2^9}{10^3}(2\epsilon_0^2)^{3/4}\left(\log N_{1-J}\right)^{3/4}
\le\left(\log N_{1-J}\right)^{3/4}.
\end{gather}
It follows from the definition \eqref{10-26} and the result of the previous item that
\begin{gather}\notag
\prod_{n=j+1}^{h}\left(2^9A^3\log\frac{N_{n-J}}{\widetilde{N}_{n-1-J}}\right)=2^9A^3\log N_{1-J}
\prod_{n=2}^{h}\left(2^9A^3\log\frac{N_{n-J}}{\widetilde{N}_{n-1-J}}\right)\le\\\le
(\log N_{1-J})^{7/4}\left((1-\epsilon_0)^{j_0(1,h)}\log N\right)^{\frac{7}{4}(h-j-1)}.\label{10-35}
\end{gather}
We obtain by the definition \eqref{10-8} that
$$\log N_{1-J}=\frac{1}{2-\epsilon_0}(1-\epsilon_0)^{J}\log N\le(1-\epsilon_0)^{J}\log N.$$
Substituting this estimate into \eqref{10-35} we obtain
\begin{gather}\label{10-36}
\prod_{n=j+1}^{h}\left(2^9A^3\log\frac{N_{n-J}}{\widetilde{N}_{n-1-J}}\right)\le
\left((1-\epsilon_0)^{j_0(1,h)}\log N\right)^{\frac{7}{4}(h-j)}(1-\epsilon_0)^{\frac{7}{4}(J-j_0(1,h))}.
\end{gather}
Taking account of$j_0\le J$ and $1-\epsilon_0<1,$ we obtain
\begin{gather}\label{10-37}
\prod_{n=j+1}^{h}\left(2^9A^3\log\frac{N_{n-J}}{\widetilde{N}_{n-1-J}}\right)\le
\left((1-\epsilon_0)^{j_0(1,h)}\log N\right)^{\frac{7}{4}(h-j)}.
\end{gather}
The fact that $j_0(1,h)=j_0(0,h)$ completes the proof.
\end{description}
Lemma is proved.
\end{proof}
\end{Le}
\section{The ensemble: constructing the set $\Omega_N.$ }\label{postoenie}
In this section we construct a set $\Omega_N.$ It will be proved in §\ref{ansambl-property} that this set is $(\epsilon_0,N)$~--ensemble. We construct the set by the inductive algorithm with the steps numbered by indexes $1,2,\ldots,2J+1.$
\begin{enumerate}
  \item The first (the starting) step.\par
We set
\begin{gather}\label{11-1}
M=M_1=N_{1-J}.
\end{gather}
Because of the lower bound in \eqref{10-17} the condition \eqref{9-1} obviously holds. So we can run the algorithm of §\ref{predansambl} to generate the set $\Xi(M).$ During the algorithm we also obtain the numbers $L=L(M),p=p(M),k=k(M).$ By the construction the number $L$ belongs to the segment $\left[\frac{M}{64A^2}, M\right],$ so we can set
\begin{gather}\label{11-2}
L=\alpha_1M_1=\alpha_1 N_{1-J},
\end{gather}
where $\alpha_1$ is a number from
\begin{gather}\label{11-3}
\left[\frac{1}{64A^2}, 1\right].
\end{gather}
Let rename the returned pre-ensemble $\Xi(M)$ and numbers $L,p$ and $k$ to
\begin{gather*}
\Xi(M)=\Xi_1,\,L=L_1,\,p=p_1,\,k=k_1.
\end{gather*}
For the next step of the algorithm we define the number $M_2$
\begin{gather}\label{11-4}
M_2=\frac{N_{2-J}}{(1+\varphi^{-2})\alpha_1 N_{1-J}},
\end{gather}
where $\varphi$ id from \eqref{9-22}.
  \item The step with the number $j,$ where $2\le j\le2J+1$ (the inductive step).\par
Write $M=M_j.$ According to the inductive assumption the number $M_j$ has been defined on the previous step by the formula
\begin{gather}\label{11-5}
M_j=\frac{N_{j-J}}{(1+\varphi^{-2})\alpha_{j-1} N_{j-1-J}},
\end{gather}
where $\alpha_{j-1}$ is a number from \eqref{11-3}. To verify for such $M$ the condition \eqref{9-1} it is sufficient to apply bounds of Lemma \ref{lemma10-3} having put $m=j-J.$
Hence we can run the algorithm of §\ref{predansambl} to generate $\Xi(M).$ Besides there exists a number $\alpha_j$ from the interval \eqref{11-3} such that for the parameter $L$ the following equation holds
\begin{gather}\label{11-6}
L=\alpha_jM=\frac{\alpha_j N_{j-J}}{(1+\varphi^{-2})\alpha_{j-1} N_{j-1-J}}.
\end{gather}
We rename $\Xi(M)$ to $\Xi_j,$ the number $L$ to $L_j,$ the quantity $p$ to $p_j,$ the length $k$ to $k_j.$  If $j\le 2J,$ then the number $M_{j+1},$ which will be used in the next step, is defined by the equation
\begin{gather*}
M_{j+1}=\frac{N_{j+1-J}}{(1+\varphi^{-2})\alpha_{j} N_{j-J}}.
\end{gather*}
If $j=2J+1,$ then the notation $M_{j+1}$ is of no use, as the algorithm is completed.
\end{enumerate}
We now define the ensemble $\Omega_N$ writing all the sets generated in the algorithm for one another
\begin{gather*}
\Omega_N=\Xi_1\Xi_2\ldots\Xi_{2J}\Xi_{2J+1}.
\end{gather*}
It means that the set $\Omega_N$ consists of all possible products of the form

$\gamma_1\gamma_2\ldots\gamma_{2J}\gamma_{2J+1},$ with $\gamma_1\in\Xi_1,\,
\gamma_2\in\Xi_2,\ldots,
\gamma_{2J+1}\in\Xi_{2J+1}.$\par
To prove that $\Omega_N$ is really an ensemble we need two technical lemmas concerning quantities $L_j.$
We will use the following notation
\begin{gather}\label{12-4}
f(x)=O_1(g(x)),\,\mbox{if}\,|f(x)|\le  g(x).
\end{gather}

\begin{Le}\label{lemma12-2}
The following inequality holds
\begin{gather}\label{12-7}
\sum_{n=1}^{2J+1}\frac{1}{\log L_n}\le\frac{1}{16000}.
\end{gather}
\begin{proof}
Let prove that numbers  $L_j$ satisfy the following inequality
\begin{gather}\label{12-8}
L_j\ge\frac{1}{64A^2(1+\varphi^{-2})}\frac{N_{j-J}}{N_{j-J-1}}\ge\frac{N_{j-J}}{100A^2N_{j-J-1}}.
\end{gather}
Actually, for $j>1$ the inequality \eqref{12-8} follows from the definition \eqref{11-6}.
For $j=1$ to deduce the same inequality \eqref{12-8} from \eqref{11-2} it is sufficient to know that $N_{-J}\ge1.$ The last inequality follows from inequalities \eqref{10-12} and \eqref{10-17} with $m=-J:$
\begin{gather*}
N_{-J}\ge N_{1-J}^{1-\epsilon_0}\ge\exp\left(\frac{10^4A^4}{2\epsilon_0^2}(1-\epsilon_0)\right)\ge1.
\end{gather*}
Thus, the inequality \eqref{12-8} is proved. It follows from the bound \eqref{10-16} that
\begin{gather*}
\left(\frac{N_m}{N_{m-1}}\right)^{1/2}\ge\exp\left(\frac{10^4A^4}{4}\right)>100A^2,\quad -J\le m\le J-1.
\end{gather*}
In that case, if $j\le2J,$ then using \eqref{10-11} the estimate \eqref{12-8} can be resumed
\begin{gather}\label{12-9}
L_j\ge\left(\frac{N_{j-J}}{N_{j-J-1}}\right)^{\frac{1}{2}}\ge N^{\frac{1}{4}\epsilon_0(1-\epsilon_0)^{|j-J-1|}},
\end{gather}
Hence
\begin{gather}\label{12-10}
\log L_j\ge\frac{1}{4}\epsilon_0(1-\epsilon_0)^{|j-J-1|}\log N,\quad j\le2J.
\end{gather}
If $j=2J+1,$ then from the lower bound in \eqref{10-17} we obtain in a similar way
\begin{gather*}
L_{2J+1}\ge\left(\frac{N_{J+1}}{N_{J}}\right)^{\frac{1}{2}}\ge \exp\left(\frac{10^4A^4}{4\epsilon_0^2}\right),
\end{gather*}
whence it follows that
\begin{gather}\label{12-11}
\log L_j\ge\frac{10^4A^4}{4\epsilon_0^2}.
\end{gather}
Substituting the estimates \eqref{12-10} and \eqref{12-11} into the sum in \eqref{12-7} one has
\begin{gather}\label{12-12}
\sum_{n=1}^{2J+1}\frac{1}{\log L_n}\le\frac{4}{\epsilon_0}\sum_{n=0}^{2J}\frac{1}{(1-\epsilon_0)^{|n-J|}\log N}+\frac{4\epsilon_0^2}{10^4A^4}\le\frac{8}{\epsilon_0}\sum_{n=0}^{J}\frac{1}{(1-\epsilon_0)^{n}\log N}+\frac{4\epsilon_0^2}{10^4A^4}.
\end{gather}
Let estimate the geometric progression from  \eqref{12-12}:
\begin{gather*}
\sum_{n=0}^{J}\frac{1}{\epsilon_0(1-\epsilon_0)^{n}}\le\frac{1}{\epsilon_0(1-\epsilon_0)^{J}}\sum_{n=0}^{\infty}(1-\epsilon_0)^{n}\le
\frac{1}{\epsilon_0^2(1-\epsilon_0)^{J}}\le\frac{\log N}{10^4A^4},
\end{gather*}
since \eqref{10-4}. Substituting this bound into \eqref{12-12} we obtain
\begin{gather*}
\sum_{n=1}^{2J+1}\frac{1}{\log L_n}\le\frac{8}{\log N}\frac{\log N}{10^4A^4}+\frac{4\epsilon_0^2}{10^4A^4}=
\frac{4(2+\epsilon_0^2)}{10^4A^4}<\frac{1}{10^3A^4}\le\frac{1}{16000}.
\end{gather*}
This completes the proof of the lemma.
\end{proof}
\end{Le}
To state the following lemma we suppose the real numbers
\begin{gather*}
\Pi_1,\Pi_2,\ldots,\Pi_{2J+1}
\end{gather*}
to satisfy the relations
\begin{gather}\label{12-13}
\Pi_j=\left(1+2O_1(\log^{-1}L_j)\right)^2\prod_{n=1}^{j-1}\left(1+2O_1(\log^{-1}L_n)\right)^3,
\end{gather}
where the product over the empty set is regarded to be equal to one.
\begin{Le}\label{lemma12-3}
For any $j=1,2,\ldots,2J+1$ the following bound holds
\begin{gather}\label{12-14}
\exp(-10^{-3})\le\Pi_j\le\exp(10^{-3}).
\end{gather}
\begin{proof}
Taking the logarithm of the equation \eqref{12-13} and bounding from above the absolute value of the sum by the sum of absolute values we obtain
\begin{gather}\notag
|\log\Pi_j|\le2|\log\left(1+2O_1(\log^{-1}L_j)\right)|+3\sum_{n=1}^{j-1}|\log\left(1+2O_1(\log^{-1}L_n)\right)|\le\\\le
3\sum_{n=1}^{2J+1}|\log\left(1+2O_1(\log^{-1}L_n)\right)|.\label{12-15}
\end{gather}
In view of Lemma \ref{lemma12-2}, every number $\log^{-1}L_n$ for $n=1,2,\ldots,2J+1$ is less than $\frac{1}{16000};$ and in particular every number $2O_1(\log^{-1}L_n)$ belongs to the segment $\left[-\frac{1}{2},\frac{1}{2}\right].$ But for any $z$ in the segment $-\frac{1}{2}\le z\le\frac{1}{2}$ the inequality $|\log(1+z)|\le |z|\log4$ holds. Then by \eqref{12-15} we obtain
\begin{gather}\label{12-21}
|\log\Pi_j|\le3\sum_{n=1}^{2J+1}\left|2O_1(\log^{-1}L_n)\right|\log4<12\sum_{n=1}^{2J+1}\log^{-1}L_n.
\end{gather}
Using Lemma \ref{lemma12-2} we obtain
\begin{gather}\label{12-22}
|\log\Pi_j|\le\frac{12}{16000}<10^{-3}.
\end{gather}
The inequality \eqref{12-14} follows immediately from  \eqref{12-22}. This completes the proof of the lemma.
\end{proof}
\end{Le}
\section{Properties of $\Omega_N.$ It is really an ensemble!}\label{ansambl-property}
In this section we prove that the constructed set $\Omega_N$ is an ensemble, that is, it satisfies the definition of ensemble in §\ref{section idea BK}. Unique expansion is the easiest property to verify. Actually, if
$$\Omega^{(1)}=\Xi_1\Xi_2\ldots\Xi_{j},$$
$$\Omega^{(2)}=\Xi_{j+1}\Xi_{j+2}\ldots\Xi_{2J+1},$$
then firstly $\Omega_N=\Omega^{(1)}\Omega^{(2)}.$ Secondly, as the representation of a matrix in the form \eqref{matrix} is unique then the implication $\eqref{8-19}\Rightarrow\eqref{8-20}$ holds since the length $D(\gamma)=k_j$ is fixed for all $\gamma\in\Xi_{j}$ (the property (iii) in §\ref{section idea BK}), for each $j=1,2,\ldots,2J+1.$ \par
The next purpose is to prove that $\Omega_N$ is a pre-ensemble.
\begin{Le}\label{lemma13-1}
For any $j$ in the segment
\begin{gather}\label{13-1}
1\le j\le2J+1,
\end{gather}
for any  collection of matrices
\begin{gather}\label{13-2}
\xi_1\in\Xi_1,\,\xi_2\in\Xi_2,\ldots,\xi_j\in\Xi_j,
\end{gather}
one can find a number $\Pi_j,$ satisfying the equality \eqref{12-13}, such that
\begin{gather}\label{13-3}
\|\xi_1\xi_2\ldots\xi_j\|=\alpha_jN_{j-J}\Pi_j.
\end{gather}
\begin{proof}
Let first $j=1.$ Then, by the construction of the pre-ensemble $\Xi_1$ (§\ref{predansambl}) and by the equation \eqref{11-2}, the following equation holds
\begin{gather}\label{13-4}
\|\xi_1\|=\alpha_1N_{1-J}(1+O_1(\log^{-1}L_1)).
\end{gather}
Since
\begin{gather*}
1+O_1(\log^{-1}L_1)=(1+2O_1(\log^{-1}L_1))^2=\Pi_1,
\end{gather*}
then substituting the last equation into \eqref{13-4} one has
\begin{gather}\label{13-5}
\|\xi_1\|=\alpha_1N_{1-J}\Pi_1,
\end{gather}
and in the case $j=1$ lemma is proved.\par
We now assume that lemma is proved for some $j,$ such that $1\le j\le2J,$ and prove that it holds for $j+1.$ It follows from
\eqref{continuant identity} that
\begin{gather}\label{13-6}
\|\xi_1\xi_2\ldots\xi_j\xi_{j+1}\|=\|\xi_1\xi_2\ldots\xi_j\|\|\xi_{j+1}\|
(1+[\overleftarrow{D}(\xi_j),\overleftarrow{D}(\xi_{j-1}),\ldots,\overleftarrow{D}(\xi_1)][D(\xi_{j+1})]),
\end{gather}
where $D(\gamma),$ as usual, denotes the sequence $D(\gamma)=\{d_1,d_2,\ldots,d_k\},$ where
\begin{equation*}
\gamma=
\begin{pmatrix}
0 & 1 \\
1 & d_1
\end{pmatrix}
\begin{pmatrix}
0 & 1 \\
1 & d_2
\end{pmatrix}\ldots
\begin{pmatrix}
0 & 1 \\
1 & d_k
\end{pmatrix}.
\end{equation*}
It follows immediately from Lemma \ref{pred9-1} that
\begin{gather}\label{13-7}
[D(\xi_{j+1})]=\varphi^{-1}+2O_1(\log^{-1}L_{j+1}),\,
[\overleftarrow{D}(\xi_j)\ldots,\overleftarrow{D}(\xi_1)]=\varphi^{-1}+2O_1(\log^{-1}L_{j}).
\end{gather}
Substituting \eqref{13-7} into \eqref{13-6}, we obtain
\begin{gather}\label{13-9}
\|\xi_1\xi_2\ldots\xi_j\xi_{j+1}\|=\|\xi_1\xi_2\ldots\xi_j\|\|\xi_{j+1}\|
\left(1+\varphi^{-2}\right)\left(1+2O_1(\log^{-1}L_j)\right)\left(1+2O_1(\log^{-1}L_{j+1})\right).
\end{gather}
By the inductive hypothesis we have
\begin{gather}\label{13-10}
\|\xi_1\xi_2\ldots\xi_j\|=\alpha_jN_{j-J}\Pi_j.
\end{gather}
And by the construction of the ensemble $\Omega_N$ (§\ref{postoenie}, "The step with the number $j+1.$" ) the following equation holds
\begin{gather}\label{13-11}
\|\xi_{j+1}\|=L_{j+1}(1+O_1(\log^{-1}L_{j+1})=
\frac{\alpha_{j+1} N_{j+1-J}}{(1+\varphi^{-2})\alpha_{j} N_{j-J}}(1+O_1(\log^{-1}L_{j+1})).
\end{gather}
Substituting \eqref{13-10} and \eqref{13-11} into \eqref{13-9} and making cancelations, we obtain
\begin{gather}\label{13-12}
\|\xi_1\xi_2\ldots\xi_j\xi_{j+1}\|=
\alpha_{j+1} N_{j+1-J}\widetilde{\Pi}_{j+1},
\end{gather}
where
\begin{gather}\label{13-13}
\widetilde{\Pi}_{j+1}=(1+O_1(\log^{-1}L_{j+1}))\left(1+2O_1(\log^{-1}L_j)\right)\left(1+2O_1(\log^{-1}L_{j+1})\right)\Pi_j.
\end{gather}
Using the definition of $\Pi_j$ by the equation \eqref{12-13} we obtain
\begin{gather*}
\widetilde{\Pi}_{j+1}=(1+2O_1(\log^{-1}L_{j+1}))^2\left(1+2O_1(\log^{-1}L_j)\right)\Pi_j=\Pi_{j+1}
\end{gather*}
and hence
\begin{gather*}
\|\xi_1\xi_2\ldots\xi_j\xi_{j+1}\|=
\alpha_{j+1} N_{j+1-J}\Pi_{j+1}.
\end{gather*}
The lemma is proved.
\end{proof}
\end{Le}
\begin{Le}\label{lemma13-2}
For any  collection of matrices \eqref{13-2}, for any numbers $j,\,h$ in the interval
\begin{gather}\label{13-14}
1\le j\le2J+1,\quad j<h\le 2J+1
\end{gather}
the following inequalities hold
\begin{gather}\label{13-15}
\frac{1}{70A^2}N_{j-J}\le\|\xi_1\xi_2\ldots\xi_j\|\le1,01N_{j-J},
\end{gather}
\begin{gather}\label{13-16}
\frac{1}{70A^2}N\le\|\xi_1\xi_2\ldots\xi_{2J+1}\|\le1,01N,
\end{gather}
\begin{gather}\label{13-17}
\frac{1}{150A^2}\frac{N_{h-J}}{N_{j-J}}\le\|\xi_{j+1}\xi_{j+2}\ldots\xi_h\|\le73A^2\frac{N_{h-J}}{N_{j-J}};
\end{gather}
moreover, for $j\le2J$ one has
\begin{gather}\label{13-18}
\frac{1}{150A^2}\frac{N}{N_{j-J}}\le\|\xi_{j+1}\xi_{j+2}\ldots\xi_{2J+1}\|\le73A^2\frac{N}{N_{j-J}}.
\end{gather}
\begin{proof}
First we prove the inequality  \eqref{13-15}. Recall that by the construction of the set $\Omega_N$ the following inequality holds
\begin{gather}\label{13-19}
\frac{1}{64A^2}\le\alpha_{j}\le1,
\end{gather}
and by Lemma \ref{lemma12-3} we also have
\begin{gather}\label{13-20}
\exp(-10^{-3})\le\Pi_j\le\exp(10^{-3}).
\end{gather}
Substituting \eqref{13-19} and \eqref{13-20} into \eqref{13-3}, we obtain \eqref{13-15}. In particular, as $N_{J+1}=N,$ then by using \eqref{13-15} for $j=2J+1$ we obtain \eqref{13-16}.\par
Now we prove the inequality \eqref{13-17}. To do this we denote
\begin{gather*}
W(j,h)=\|\xi_{j}\xi_{j+1}\ldots\xi_h\|
\end{gather*}
and rewrite the inequality \eqref{continuant inequality} in the form
\begin{gather*}
W(1,j)W(j+1,h)\le W(1,h)\le2W(1,j)W(j+1,h).
\end{gather*}
Hence, applying the inequality \eqref{13-15} twice we obtain
\begin{gather*}
W(j+1,h)\ge\frac{W(1,h)}{2W(1,j)}\ge\frac{N_{h-J}/(70A^2)}{2,02N_{j-J}}\ge\frac{1}{150A^2}\frac{N_{h-J}}{N_{j-J}},
\end{gather*}
and in the same way
\begin{gather*}
W(j+1,h)\le\frac{W(1,h)}{W(1,j)}\le\frac{1,01N_{h-J}}{N_{j-J}/(70A^2)}\le73A^2\frac{N_{h-J}}{N_{j-J}}.
\end{gather*}
These prove the inequality \eqref{13-17}. Putting $h=2J+1$ in it we obtain \eqref{13-18}.
The lemma is proved.
\end{proof}
\end{Le}
For integers $j$ and $h,$ such that
\begin{gather}\label{13-21}
0\le j<h\le 2J+1,
\end{gather}
we put
\begin{gather}\label{13-22}
\Omega(j,h)=\Xi_{j+1}\Xi_{j+2}\ldots\Xi_{h}.
\end{gather}
\begin{Le}\label{lemma13-3}
The following estimate holds
\begin{gather}\label{13-23}
\left|\Omega(0,j)\right|\le9N_{j-J}^{2\delta}.
\end{gather}
\begin{proof}
By definition,  for $\gamma\in\Omega(0,j)$ one has $\gamma=\xi_{1}\xi_{2}\ldots\xi_{j}$ for a collection of matrices \eqref{13-2}. So, it follows from the inequality  \eqref{13-15} that
\begin{gather}\label{13-24}
\|\gamma\|\le1,01N_{j-J}.
\end{gather}
The number of matrices $\gamma,$ satisfying the inequality \eqref{13-24} can be bounded by Theorem \ref{Kan-Hen}. Estimating the result from above we obtain \eqref{13-23}. This completes the proof of the lemma.
\end{proof}
\end{Le}
Recall that parameters $\widetilde{N}_{n-J}$ and $j_0(j,h)$ were introduced by formulae \eqref{10-26} and \eqref{10-28}.  Note that the restrictions \eqref{13-21} on $j$ and $h$ coincide with the restrictions \eqref{10-27}.
\begin{Le}\label{lemma13-4}
For $j$ and $h$ in \eqref{13-21} the following bound holds
\begin{gather}\label{13-25}
\left|\Omega(j,h)\right|\ge\frac{1}{\left((1-\epsilon_0)^{j_0}\log N\right)^{7(h-j)}}
\left(\frac{N_{h-J}}{\widetilde{N}_{j-J}}\right)^{2\delta},
\end{gather}
where $j_0=j_0(j,h).$
\begin{proof}
Multiplying the lower bounds \eqref{9-20} we obtain
\begin{gather}\label{13-26}
\left|\Omega(j,h)\right|\ge\prod_{n=j+1}^{h}\left|\Xi_n\right|\ge\prod_{n=j+1}^{h}\frac{L_n^{2\delta}}{2^{18}A^5\log^4L_n}.
\end{gather}
It follows from formulae \eqref{11-2} and \eqref{11-6} that
\begin{gather}\label{13-27}
\frac{N_{n-J}}{c_1\widetilde{N}_{n-1-J}}\le L_n\le\frac{c_1N_{n-J}}{\widetilde{N}_{n-1-J}},
\end{gather}
where
\begin{gather}\label{13-28}
c_1=64A^2(1+\varphi^{-2})\le 2^7A^2.
\end{gather}
Let estimate the product of the numerators in \eqref{13-26}. Applying \eqref{13-27} and \eqref{13-28} we obtain
\begin{gather*}
\prod_{n=j+1}^{h}L_n^{2\delta}\ge\prod_{n=j+1}^{h}\left(\frac{N_{n-J}}{2^7A^2\widetilde{N}_{n-1-J}}\right)^{2\delta}
\end{gather*}
After the cancelations we obtain
\begin{gather}\label{13-29}
\prod_{n=j+1}^{h}L_n^{2\delta}\ge
\left(\frac{N_{h-J}}{\widetilde{N}_{j-J}}\right)^{2\delta}\prod_{n=j+1}^{h}(2^7A^2)^{-2\delta}.
\end{gather}
So the estimate \eqref{13-26} can be resumed in such a way
\begin{gather}\notag
\left|\Omega(j,h)\right|\ge
\left(\frac{N_{h-J}}{\widetilde{N}_{j-J}}\right)^{2\delta}
\prod_{n=j+1}^{h}\frac{(2^7A^2)^{-2\delta}}{2^{18}A^5\log^4L_n}\ge\\\ge
\left(\frac{N_{h-J}}{\widetilde{N}_{j-J}}\right)^{2\delta}
\left(\prod_{n=j+1}^{h}(2^8A^3\log L_n)\right)^{-4}.\label{13-30}
\end{gather}
The last product in  \eqref{13-30} will be estimated separately. In view of the upper bound in \eqref{13-27} we have
\begin{gather}\label{13-31}
\prod_{n=j+1}^{h}(2^8A^3\log L_n)\le\prod_{n=j+1}^{h}\left(2^8A^3\left(
\log(2^7A^2)+\log\frac{N_{n-J}}{\widetilde{N}_{n-1-J}}\right)\right).
\end{gather}
Applying Lemma \ref{lemma10-3}, we obtain
\begin{gather}\label{13-31}
\log\frac{N_{n-J}}{\widetilde{N}_{n-1-J}}\ge\frac{10^4A^4}{2\epsilon_0}\ge\log(2^7A^2),
\end{gather}
hence, applying Lemma \ref{lemma10-5}, we obtain
\begin{gather}\label{13-32}
\prod_{n=j+1}^{h}(2^8A^3\log L_n)\le\prod_{n=j+1}^{h}\left(2^9A^3\log\frac{N_{n-J}}{\widetilde{N}_{n-1-J}}\right)\le
\left((1-\epsilon_0)^{j_0}\log N\right)^{\frac{7}{4}(h-j)}.
\end{gather}
Substituting the estimate \eqref{13-32} into \eqref{13-30}, we obtain \eqref{13-25}. This completes the proof of the lemma.
\end{proof}
\end{Le}
\begin{Th}\label{theorem13-1}
For $j$ and $h$ in \eqref{13-21} the following estimate holds
\begin{gather}\label{13-34}
\left|\Omega(j,h)\right|\ge
\left(\frac{N_{h-J}}{\widetilde{N}_{j-J}}\right)^{2\delta}
\exp\left(-\left(\frac{\log\log N}{\log(1-\epsilon_0)}+j_0\right)^2\right),
\end{gather}
where $j_0=j_0(j,h).$
\begin{proof}
It follows from \eqref{13-25} that it is enough to prove the inequality
\begin{gather*}
\exp\left(\left(\frac{\log\log N}{\log(1-\epsilon_0)}+j_0\right)^2\right)\ge
\exp\left(7(h-j)\left(\log\log N+j_0\log(1-\epsilon_0)\right)\right).
\end{gather*}
Hence, it is sufficient to prove
\begin{gather*}
7(h-j)\left(\log\log N+j_0\log(1-\epsilon_0)\right)\le
\left(\frac{\log\log N+j_0\log(1-\epsilon_0)}{\log(1-\epsilon_0)}\right)^2.
\end{gather*}
One can readily obtain from \eqref{10-2} that
\begin{gather}\label{13-35}
J\le\frac{\log\log N}{-\log(1-\epsilon_0)}-1,
\end{gather}
and since $j_0\le J,$ so one has
\begin{gather*}
\log\log N+j_0\log(1-\epsilon_0)>0.
\end{gather*}
Thus, it is sufficient to prove that
\begin{gather}\label{13-36}
7(h-j)\le
\frac{\log\log N+j_0\log(1-\epsilon_0)}{\log^2(1-\epsilon_0)}=
\frac{1}{-\log(1-\epsilon_0)}\left(\frac{\log\log N}{-\log(1-\epsilon_0)}-j_0\right).
\end{gather}
We observe that as $\epsilon_0\in\left(0,\frac{1}{2500}\right),$ so $$-7\log(1-\epsilon_0)\le\frac{1}{7}.$$ It follows from this, \eqref{13-35} and \eqref{13-36} that it is sufficient to prove
\begin{gather}\label{13-37}
\frac{1}{7}(h-j)\le J+1-j_0.
\end{gather}
For $j_0=0$ the inequality \eqref{13-37} follows from the trivial bound $h-j\le2J+1.$ For $j_0\neq0$ $\left( \mbox{hence,} j>J\,\mbox{or}\,h\le J\right)$ it follows from \eqref{10-28} and \eqref{10-29} that
\begin{gather}\label{13-38}
h-j\le
\left\{
              \begin{array}{ll}
                2J+1-j=J+1-j_0, & \hbox{если $j>J$;} \\
                h=J+1-j_0, & \hbox{если $h\le J$.}
              \end{array}
\right.
\end{gather}
Hence the inequality \eqref{13-37} follows. This completes the proof of the theorem.
\end{proof}
\end{Th}
\begin{Sl}\label{conseq13-1}
For $N\ge\exp(\epsilon_0^{-5})$ the following estimate holds
\begin{gather}\label{13-39}
N^{2\delta-\epsilon_0}\le N^{2\delta}
\exp\left(-\left(\frac{\log\log N}{\epsilon_0}\right)^2\right)\le
\#\Omega_N\le9N^{2\delta}.
\end{gather}
\begin{proof}
To prove the upper bound we apply Lemma \ref{lemma13-3} with $j=0,\,h=2J+1.$ To prove the lower bound we use Theorem \ref{theorem13-1} putting $j=0,\,h=2J+1$ in it (hence,  $j_0=0.$)
As $\epsilon_0\in\left(0,\frac{1}{2500}\right),$ so we have $\log^2(1-\epsilon_0)>\epsilon_0^2$ and obtain the lower bound in \eqref{13-39}. The corollary is proved.
\end{proof}
\end{Sl}
\begin{Sl}\label{conseq13-2}
The set $\Omega_N$ is a $N$~--pre-ensemble.
\begin{proof}
It follows from the inequality \eqref{13-16} that the property $\|\gamma\|\asymp N$ holds for each matrix $\gamma\in\Omega_N.$ By the Corollary \ref{conseq13-1} we have $\#\Omega_N\ge_{\epsilon}N^{2\delta-\epsilon}.$  To put it another way, both items of the definition of pre-ensemble hold. The corollary is proved.
\end{proof}
\end{Sl}
We now verify the last property of ensemble related to the relations  \eqref{8-21} and \eqref{8-22}.
\begin{Le}\label{lemma13-5}
For any $M$ in the interval \eqref{10-19} there exist $j$ and $h$ in the intervals \eqref{10-20}, such that for any collection of matrices \eqref{13-2} the following inequalities hold
\begin{gather}\label{13-40}
0,99\|\xi_1\xi_2\ldots\xi_j\|^{1-\epsilon_0}\le M\le
70A^2\|\xi_1\xi_2\ldots\xi_j\|,
\end{gather}
\begin{gather}\label{13-41}
\frac{1}{73A^2}\|\xi_{h+1}\xi_{h+2}\ldots\xi_{2J+1}\|^{1-\epsilon_0}\le M\le
150A^2\|\xi_{h+1}\xi_{h+2}\ldots\xi_{2J+1}\|.
\end{gather}
\begin{proof}
Let $M$ be fixed in the interval \eqref{10-19}. Then using Lemma \ref{lemma10-4} we find $j$ and $h$ in \eqref{10-20}, such that the inequalities  \eqref{10-21} and \eqref{10-22} hold. We note that a bilateral bound on the number $N_{j-J}$ in terms of $\|\xi_1\xi_2\ldots\xi_j\|$ follows from \eqref{13-15}. To obtain the inequality \eqref{13-40} one should substitute this bilateral bound into \eqref{10-21}. To obtain the inequality \eqref{13-41} one should substitute a bilateral bound on $N/N_{h-J},$ following from the inequality \eqref{13-18} with $j=h,$ into \eqref{10-22}. This completes the proof of the lemma.
\end{proof}
\end{Le}
\begin{Sl}\label{conseq13-3}
For any $M,$ satisfying the inequality
\begin{gather}\label{13-42}
\exp\left(\frac{10^5A^4}{\epsilon_0^2}\right)\le M\le N\exp\left(-\frac{10^5A^4}{\epsilon_0^2}\right),
\end{gather}
there exist indexes $j$ and $h$ in the intervals \eqref{10-20}, such that for any collection of matrices \eqref{13-2}the inequalities \eqref{13-40} and \eqref{13-41} hold.
\begin{proof}
By applying Lemma \ref{lemma10-3}, it follows from the inequality \eqref{13-42} that the inequality \eqref{10-19} holds. It is sufficient now to apply Lemma \ref{lemma13-5}. The corollary is proved.
\end{proof}
\end{Sl}
\begin{Th}\label{theorem13-2}
The set $\Omega_N$ is a bilateral $(\epsilon_0,N)$~--ensemble.
\begin{proof}
The unique expansion of the products, which are equal to $\Omega_N,$ has been proved in the beginning of §\ref{ansambl-property}. The property of $\Omega_N$ to be a $N$~--pre-ensemble has been proved in the Corollary \ref{conseq13-2}. The right $(\epsilon_0,N)$~--ensemble property is proved by the inequality \eqref{13-41}, the left~--by \eqref{13-40}, since it follows from the Corollary \ref{conseq13-3} that these inequalities hold for any $M,$ satisfying the inequality \eqref{13-42}. This completes the proof of the theorem.
\end{proof}
\end{Th}
Thus the main purpose of the section is achieved. However, we formulate a few more properties of a bilateral ensemble. These properties will be of use while estimating exponential sums.\par
For $M,$ satisfying the inequality \eqref{13-42}, we denote by $\jj$ and $\h$ the numbers $j$ and $h$ from Lemma \ref{lemma10-4} corresponding to $M$. For brevity we will write that numbers $\jj,\,\h$ corresponds to $M.$ In the following theorem $\jj^{(1)}$ corresponds to $M^{(1)},$ and $\h^{(3)}$ corresponds to $M^{(3)}.$
\begin{Le}\label{lemma13-6}
Let the inequality
\begin{gather}\label{13-43}
M^{(1)}M^{(3)}< N^{1-\epsilon_0}.
\end{gather}
holds for $M^{(1)}$ and $M^{(3)}$ in the interval \eqref{13-42}. Then $\jj^{(1)}<\h^{(3)},$ and for $M=M^{(1)}$ the inequality \eqref{13-40} holds and for $M=M^{(3)}.$ the inequality \eqref{13-41} holds.
\begin{proof}
It is sufficient to verify the condition $\jj^{(1)}<\h^{(3)}.$ Then the statement of Lemma \ref{lemma13-6} will follow from Lemma \ref{lemma13-5} and Corollary \ref{conseq13-3}.\par
We recall that $\h^{(3)}=2J-\jj^{(3)}+1.$ Hence, it is sufficient to prove that $$\jj^{(1)}+\jj^{(3)}<2J+1.$$ Assume the contrary.
It follows from \eqref{10-23} that
\begin{gather*}
N_{\jj^{(1)}-1-J}N_{\jj^{(3)}-1-J}\le M^{(1)}M^{(3)}.
\end{gather*}
Since $\jj^{(1)}+\jj^{(3)}\ge2J+1,$ one can consider the inequality $\jj^{(1)}\ge J+1$ to hold.  It follows from the relation $\jj^{(3)}-1-J\ge-\jj^{(1)}+J$ and the increasing of the sequence $\left\{N_j\right\}_{j=-J-1}^{J+1},$ that
\begin{gather*}
N_{\jj^{(1)}-1-J}N_{-\jj^{(1)}+J}\le M^{(1)}M^{(3)},\quad \jj^{(1)}\ge J+1.
\end{gather*}
We put $m=\jj^{(1)}-J-1,$  then
\begin{gather*}
N_{m}N_{-m-1}\le M^{(1)}M^{(3)},\quad m\ge0.
\end{gather*}
Because of \eqref{10-8}, for $ m\ge0$ we obtain
\begin{gather*}
N_{m}N_{-m-1}=N^{1-\frac{1}{2-\epsilon_0}(1-\epsilon_0)^{m}+\frac{1}{2-\epsilon_0}(1-\epsilon_0)^{m+2}}=
N^{1-\epsilon_0(1-\epsilon_0)^{m}}\ge N^{1-\epsilon_0}.
\end{gather*}
Hence, the inequality $M^{(1)}M^{(3)}\ge N^{1-\epsilon_0}$ holds, contrary to \eqref{13-43}. This completes the proof of the lemma.
\end{proof}
\end{Le}

Recall that $\Omega(j,h)$ is defined in \eqref{13-22} and let
\begin{gather}\label{13-44}
\Omega_2(M)=\Omega(\h,2J+1),\quad\Omega_1(M)=\Omega(0,\jj),
\end{gather}
\begin{Th}\label{theorem13-3}
For any $M,$ satisfying the inequality \eqref{13-42}, the ensemble $\Omega_N$ can be represented in the following way
\begin{gather}\label{13-46-1}
\Omega_N=\Omega^{(1)}\Omega^{(3)},
\end{gather}
where
\begin{gather}\label{13-46-2}
\Omega^{(1)}=\Omega_1(M^{(1)})=\Xi_{1}\Xi_{2}\ldots\Xi_{\jj},\quad
\Omega^{(3)}=\Xi_{\jj+1}\Xi_{j+2}\ldots\Xi_{2J+1}
\end{gather}
and for any $\gamma_1\in\Omega^{(1)},\,\gamma_3\in\Omega^{(3)}$ the following inequalities hold
\begin{gather}\label{13-47}
\frac{M^{(1)}}{70A^2}\le\|\gamma_1\|\le1,03(M^{(1)})^{1+2\epsilon_0},
\end{gather}
\begin{gather}\label{13-47-1}
\frac{N}{140A^2H_1(M^{(1)})}\le\|\gamma_3\|\le\frac{73A^2N}{M^{(1)}},\quad\mbox{где}\quad H_1(M^{(1)})=1,03(M^{(1)})^{1+2\epsilon_0},
\end{gather}
\begin{proof}
In view of the Corollary \ref{conseq13-3} the inequality \eqref{13-40} holds for the matrices $\gamma_1=\xi_1\xi_2\ldots\xi_{\jj}.$ Using this inequality we obtain \eqref{13-47}. Moreover, it follows from \eqref{continuant inequality} and \eqref{13-16} that
\begin{gather}\label{13-48-1}
\frac{N}{140A^2}\le\frac{\|\xi_1\xi_2\ldots\xi_{2J+1}\|}{2}\le\|\gamma_1\|\|\gamma_2\|\le\|\xi_1\xi_2\ldots\xi_{2J+1}\|\le1,01N.
\end{gather}
Substituting the bound \eqref{13-47} into \eqref{13-48-1} we obtain the inequality \eqref{13-47-1}. This completes the proof of the theorem.
\end{proof}
\end{Th}

\begin{Th}\label{theorem13-6}
Let $M^{(4)}\ge(M^{(1)})^{2\epsilon_0}$
and let the inequality \eqref{13-42} holds for $M=M^{(1)}$ and $M=M^{(1)}M^{(4)}.$ Let  $M^{(1)}M^{(4)}$ and $M^{(6)}$ satisfy the hypotheses of Lemma \ref{lemma13-6}. Then the ensemble $\Omega_N$ can be represented in the following form
\begin{gather*}
\Omega_N=\Omega^{(1)}\Omega^{(3)}=\Omega^{(1)}\Omega^{(4)}\Omega^{(5)}\Omega^{(6)},
\end{gather*}
In doing so one has \eqref{13-46-2} and
\begin{gather*}
\Omega^{(4)}=\Xi_{\jj_1+1}\Xi_{\jj_1+2}\ldots\Xi_{\jj_4},\quad
\Omega^{(1)}\Omega^{(4)}=\Omega_1(M^{(1)}M^{(4)}),\\
\Omega^{(5)}=\Xi_{\jj_4+1}\Xi_{\jj_4+2}\ldots\Xi_{\h_6},\quad
\Omega^{(6)}=\Omega_2(M^{(6)})=\Xi_{\h_6+1}\Xi_{\h_6+2}\ldots\Xi_{2J+1},
\end{gather*}
where $\jj_1$ corresponds to $M^{(1)},$  $\jj_4$ corresponds to $M^{(1)}M^{(4)},$ and $\h_6$ to $M^{(6)}.$ Moreover for any $\gamma_{i}\in\Omega^{(i)},\,i=1,\ldots,6$ one has \eqref{13-47}, \eqref{13-47-1} and
\begin{gather}\label{13-61}
\frac{M^{(1)}M^{(4)}}{70A^2}\le\|\gamma_1\gamma_4\|\le1,03(M^{(1)}M^{(4)})^{1+2\epsilon_0},
\end{gather}
\begin{gather}\label{13-62}
\frac{M^{(4)}}{150A^2(M^{(1)})^{2\epsilon_0}}\le\|\gamma_4\|\le73A^2\frac{(M^{(4)})^{1+2\epsilon_0}}{(M^{(1)})^{2\epsilon_0}},
\end{gather}
\begin{gather}\label{13-63}
\frac{M^{(6)}}{150A^2}\le\|\gamma_6\|\le80A^{2,1}(M^{(6)})^{1+2\epsilon_0},
\end{gather}
\begin{gather}\label{13-64}
\frac{N}{25000A^5(M^{(1)}M^{(6)})^{1+2\epsilon_0}}\le\|\gamma_4\gamma_5\|\le11000A^4\frac{N}{M^{(1)}M^{(6)}}.
\end{gather}
\begin{proof}
To ensure that the partition of the ensemble is well-defined, it is enough to verify that $\jj_1<\jj_4<\h_6.$ The second inequality follows from Lemma \ref{lemma13-6} ($\jj_4=j^{(1)},\,\h_6=h^{(3)}$).For proving  $\jj_1<\jj_4$ it is sufficient to show that for any $M^{(1)}$ in the interval $N_{\jj_1-1-J}\le M^{(1)}\le N_{\jj_1-J}$ the inequality $N_{\jj_1-J}\le M^{(1)}M^{(4)}$ holds. It follows from the conditions of the theorem that
\begin{gather}\label{13-48-7}
M^{(4)}\ge(M^{(1)})^{2\epsilon_0}\ge(M^{(1)})^{\frac{1}{1-\epsilon_0}-1}\ge\frac{N_{\jj_1-J}}{M^{(1)}},
\end{gather}
where the last inequality holds because of $M^{(1)}\ge N_{\jj-J}^{1-\epsilon_0}$ see \eqref{10-21}. Thus we have proved that the partition of the ensemble is well-defined.  The bound \eqref{13-61} can be proved in the same way as \eqref{13-47}, and \eqref{13-62} follows from \eqref{13-47},\, \eqref{13-61} and \eqref{continuant inequality}. It also follows from Lemma \ref{lemma13-6} that the inequality \eqref{13-41} holds for $\gamma_6\in\Omega^{(6)}.$ Using this inequality we obtain the estimate \eqref{13-63}. Next, in the same way as \eqref{13-48-1} we obtain
\begin{gather}\label{13-66}
\frac{N}{280A^2}\le\|\gamma_1\|\|\gamma_4\gamma_5\|\|\gamma_6\|\le1,01N.
\end{gather}
Substituting \eqref{13-47} and \eqref{13-63} into \eqref{13-66}, we obtain \eqref{13-64}. This completes the proof of the theorem.
\end{proof}
\end{Th}
We write
$j_0(M)=j_0(0,\jj)=j_0(\h,2J+1),$
where $j_0(j,h)$ is defined in \eqref{10-28}. Let verify that $j_0(M)$ is well-defined. Actually, since $\h=2J-\jj+1,$ so we are to prove that
$$j_0(0,\jj)=j_0(2J-\jj+1,2J+1)\,\mbox{for }2\le\jj\le2J.$$
If $\jj\le J+1,$ then
$$j_0(0,\jj)=\left|\jj-J-1\right|=J+1-\jj=\left|2J+2-\jj-J-1\right|=j_0(2J-\jj+1,2J+1).$$
If $\jj>J+1,$ then
$$j_0(0,\jj)=j_0(2J-\jj+1,2J+1)=0.$$
Hence, $j_0(M)$ is well-defined.
\begin{Le}\label{lemma13-7}
For $M$ in the interval \eqref{13-42} the following inequality holds
\begin{gather}\label{13-51}
\frac{\log\log M^2}{\log(1-\epsilon_0)}-1\le j_0(M)+\frac{\log\log N}{\log(1-\epsilon_0)}.
\end{gather}
\begin{proof}
We consider two cases.
\begin{enumerate}
  \item Let
$\exp\left(\frac{10^5A^4}{\epsilon_0^2}\right)\le M\le N_1=N^{\frac{1}{2-\epsilon_0}}.$\par
Then it follows from Lemma \ref{lemma10-4} and \eqref{10-23} that firstly  $2\le\jj\le J+1$ and secondly
\begin{gather*}
N_{\jj-J}^{1-\epsilon_0}\le M\le N_{\jj-J}.
\end{gather*}
Hence $j_0(M)=J+1-\jj$ and applying \eqref{10-8} we obtain
\begin{gather}\label{13-52}
M\ge N_{\jj-J}^{1-\epsilon_0}=N^{\frac{1}{2-\epsilon_0}(1-\epsilon_0)^{2-\jj+J}}=
N^{\frac{1}{2-\epsilon_0}(1-\epsilon_0)^{1+j_0(M)}}\ge N^{\frac{1}{2}(1-\epsilon_0)^{1+j_0(M)}}.
\end{gather}
Taking a logarithm twice we obtain
\begin{gather*}
\log\log M^2\ge(1+j_0(M))\log(1-\epsilon_0)+\log\log N.
\end{gather*}
From this the inequality \eqref{13-51} follows.
  \item Let
$N^{\frac{1}{2-\epsilon_0}}=N_1\le M\le N\exp\left(-\frac{10^5A^4}{\epsilon_0^2}\right).$\par
Then $\jj>J+1$ and, hence, $j_0(M)=0.$ In view of Lemma \ref{lemma10-3} and \eqref{10-8} we obtain
\begin{gather}\notag
M\ge  N_{\jj-J}^{1-\epsilon_0}=N^{\left(1-\frac{1}{2-\epsilon_0}(1-\epsilon_0)^{\jj-J}\right)(1-\epsilon_0)}\ge
N^{\left(1-\frac{1}{2-\epsilon_0}(1-\epsilon_0)^{2}\right)(1-\epsilon_0)}\\\ge
N^{\frac{1}{2}(1-\epsilon_0)^{1+j_0(M)}},\label{13-53}
\end{gather}
since $1-\frac{1}{2-\epsilon_0}(1-\epsilon_0)^{2}\ge\frac{1}{2}.$ The inequality  \eqref{13-53} coincides with the bound \eqref{13-52}. Thus we obtain \eqref{13-51} in the same way.
\end{enumerate}
The lemma is proved.
\end{proof}
\end{Le}
For $M$ in the interval \eqref{13-42} we define the following function
\begin{gather}\label{13-54}
h(M)=\exp\left(-\left(\frac{\log\log M^2}{\log(1-\epsilon_0)}-1\right)^2\right).
\end{gather}
We observe that for $\epsilon_0\in(0,\frac{1}{2500})$ one has
\begin{gather}\label{13-54-1}
M^{-\epsilon_0}\le h(M), \quad\mbox{if}\quad M\ge \exp(\epsilon_0^{-5}).
\end{gather}
\begin{Th}\label{theorem13-5}
Let $M\ge \exp(\epsilon_0^{-5})$ and belongs to the interval \eqref{13-42}. Then the following bounds on the cardinality of the sets $|\Omega_1(M)|$ and $|\Omega_2(M)|$ hold
\begin{gather}\label{13-55}
M^{2\delta-\epsilon_0}\le M^{2\delta}h(M)\le\left|\Omega_1(M)\right|\le 9M^{2\delta+4\epsilon_0},
\end{gather}
\begin{gather}\label{13-56}
\left|\Omega_2(M)\right|\ge M^{2\delta}h(M)\ge M^{2\delta-\epsilon_0}.
\end{gather}
\begin{proof}
Using the definition of the set $|\Omega_1(M)|$ and the inequality
\begin{gather*}
N_{\jj-J}\le M^{\frac{1}{1-\epsilon_0}}\le M^{1+2\epsilon_0}
\end{gather*}
we obtain that the upper bound in \eqref{13-55} follows immediately  from Lemma \ref{lemma13-3}.
\begin{enumerate}
  \item We estimate the cardinality of the set $|\Omega_1(M)|$ from below.
Taking into account that ${M\le N_{\jj-J}}$ and $j_0(M)=j_0(0,\jj),$ it follows from the Theorem \ref{theorem13-1} that
\begin{gather}\label{13-57}
\left|\Omega_1(M)\right|\ge M^{2\delta}\exp\left(-\left(\frac{\log\log N}{\log(1-\epsilon_0)}+j_0(M)\right)^2\right).
\end{gather}
Using \eqref{13-35} and \eqref{13-51}, we obtain, in view of $j_0(M)\le J,$ that
\begin{gather*}
\frac{\log\log M^2}{\log(1-\epsilon_0)}-1\le j_0(M)+\frac{\log\log N}{\log(1-\epsilon_0)}\le0.
\end{gather*}
Hence,
\begin{gather*}
h(M)\le \exp\left(-\left(\frac{\log\log N}{\log(1-\epsilon_0)}+j_0(M)\right)^2\right).
\end{gather*}
Substituting the estimate into \eqref{13-57}, we obtain the lower bound in \eqref{13-55}.
  \item We estimate the cardinality of the set $|\Omega_2(M)|$ from below.
Taking into account that ${M\le \frac{N}{N_{\h-J}}}$ and $j_0(M)=j_0(\h,2J+1),$ it follows from the Theorem \ref{theorem13-1} that
\begin{gather}\label{13-58}
\left|\Omega_2(M)\right|\ge M^{2\delta}\exp\left(-\left(\frac{\log\log N}{\log(1-\epsilon_0)}+j_0(M)\right)^2\right).
\end{gather}
From this the estimate \eqref{13-56} follows in the same way.
\end{enumerate}
The theorem is proved.
\end{proof}
\end{Th}
\part{Estimates of exponential sums and integrals. A generalization of the Bourgain-Kontorovich's method.}
\section{General estimates of exponential sums over ensemble.}\label{Trigsums1}
Recall that to prove our main theorem \ref{uslov} we should obtain the maximum accurate bound on the integral
\begin{gather}\label{14-1}
\int_0^1\left|S_N(\theta)\right|^2d\theta=\int_0^1\left|\sum_{\gamma\in\Omega_{N} }e(\theta\|\gamma\|)\right|^2d\theta,
\end{gather}
where $N$ is a sufficiently large integer and for
\begin{equation}\label{14-2}
\gamma=
\begin{pmatrix}
a & b \\
c & d
\end{pmatrix}=
\begin{pmatrix}
0 & 1 \\
1 & d_1
\end{pmatrix}
\begin{pmatrix}
0 & 1 \\
1 & d_2
\end{pmatrix}\ldots
\begin{pmatrix}
0 & 1 \\
1 & d_k
\end{pmatrix}
\end{equation}
the norm $\|\gamma\|$  is defined by
\begin{gather}\label{14-3}
\|\gamma\|=d=<d_1,d_2,\ldots,d_k>=(0,1)
\begin{pmatrix}
a & b \\
c & d
\end{pmatrix}
\begin{pmatrix}
0   \\
1
\end{pmatrix}.
\end{gather}
Then
\begin{gather}\label{14-4}
S_N(\theta)=\sum_{\gamma\in\Omega_{N} }e(\theta\|\gamma\|)=\sum_{\gamma\in\Omega_{N} }e((0,1)\gamma
\begin{pmatrix}
0   \\
1
\end{pmatrix}
\theta).
\end{gather}
To estimate the sum \eqref{14-4} different methods were used in ~\cite{BK} depending on the value of $\theta.$  Suppose that a partition of the ensemble as in Theorem \ref{theorem13-3} is given. Then $\|\gamma_1\|\le H_1(M^{(1)})=H_1.$
For $n\in\{1,3\}$ we write
\begin{gather}\label{14-11}
\widetilde{\Omega}^{(n)}=
\left\{
              \begin{array}{ll}
                (0,1)\Omega^{(1)}=\left\{(0,1)g_1\Bigl| g_1\in\Omega^{(1)}\right\}, & \hbox{if $n=1$,} \\
              \Omega^{(3)}\begin{pmatrix}0  \\1 \end{pmatrix}=\left\{g_3\begin{pmatrix}0  \\1 \end{pmatrix}\Bigl| g_3\in\Omega^{(3)}\right\}, & \hbox{if $n=3$.}
                 \end{array}
\right.,
\end{gather}
Let define the following functions
\begin{gather}\label{14-12}
T(x)=\max\left\{0,1-|x|\right\},\quad S(x)=\left(\frac{\sin \pi x}{\pi x}\right)^2.
\end{gather}
It is  common knowledge ~\cite[(4.83)]{IK} that $\hat{S}(x)=T(x),$ where $\hat{f}(x)=\int_{-\infty}^{\infty}f(t)e(-xt)dt$ is the Fourier transform of the function $f(x).$ We consider $S_2(x)=3S(\frac{x}{2}).$ It is obvious that $S_2(x)$ is a nonnegative function and
\begin{gather}\label{14-13}
S_2(x)>1\, \mbox{при}\, x\in[-1,1]
\end{gather}
Since $\hat{S}_2(x)=6T(2x),$ we have $\hat{S}_2(x)\neq 0$ only for $|x|<\frac{1}{2}.$  We are to estimate the sum of the form
\begin{gather}\label{14-14}
\sigma_{N,Z}=\sum_{\theta\in Z}\left|S_N(\theta)\right|,
\end{gather}
where $Z$ is a finite subset of the interval $[0,1].$
\begin{Le}\label{lemma-14-1}
The following estimate holds
\begin{gather}\label{14-15}
\sigma_{N,Z}\le\left|\Omega^{(1)}\right|^{1/2}\Biggl(
\sum_{g_1\in\Z^2}\s\left(\frac{g_1}{H_1}\right)
\left|\sum_{\theta\in Z}\xi(\theta)\sum_{g_{3}\in\Om^{(3)}}e(g_1g_3\theta)\right|^2\Biggr)^{1/2},
\end{gather}
where $\s(x,y)=S_2(x)S_2(y),$ and $\xi(\theta)$ is a complex number with $|\xi(\theta)|=1.$
\begin{proof}
The numbers $\xi(\theta)$ are defined by the relation $\left|S_N(\theta)\right|=\xi(\theta)S_N(\theta).$ Then we obtain from \eqref{14-4} and the definition \eqref{14-11} that
\begin{gather}\label{14-16}
\sigma_{N,Z}=\sum_{\theta\in Z}\xi(\theta)S_N(\theta)=\sum_{\theta\in Z}\xi(\theta)\sum_{g_{1}\in\Om^{(1)}}
\sum_{g_{3}\in\Om^{(3)}}e(g_1g_3\theta),
\end{gather}
where $g_1,g_3$ are already vectors in $\Z^2.$ It follows from \eqref{14-16} that
\begin{gather}\label{14-17}
\sigma_{N,Z}\le\sum\limits_{g_1\in\Z^2\atop |g_1|\le H_1}\1_{g_{1}\in\Om^{(1)}}\left|
\sum_{\theta\in Z}\xi(\theta)\sum_{g_{3}\in\Om^{(3)}}e(g_1g_3\theta)\right|.
\end{gather}
We note that  the condition $|g_1|\le H_1$ does not impose any extra restrictions. Applying the Cauchy~-Schwarz inequality in \eqref{14-17} we obtain
\begin{gather}\label{14-18}
\sigma_{N,Z}\le
\left(\sum\limits_{g_1\in\Z^2\atop |g_1|\le H_1}\1_{g_{1}\in\Om^{(1)}}\right)^{1/2}
\left(\sum\limits_{g_1\in\Z^2\atop |g_1|\le H_1}
\left|\sum_{\theta\in Z}\xi(\theta)\sum_{g_{3}\in\Om^{(3)}}e(g_1g_3\theta)\right|^2\right)^{1/2}.
\end{gather}
Considering that
$$\left(\sum\limits_{g_1\in\Z^2\atop |g_1|\le H_1}\1_{g_{1}\in\Om^{(1)}}\right)^{1/2}=\left|\Omega^{(1)}\right|^{1/2}$$
and taking into account that the function $\s(x,y)>1$ for {$(x,y)\in[-1,1]^2$} and is nonnegative, we obtain \eqref{14-15}. This completes the proof of the lemma.
\end{proof}
\end{Le}
\begin{Le}\label{lemma-14-2}
Under the hypotheses of Lemma \ref{lemma-14-1} the following bound holds
\begin{gather}\label{14-19}
\sigma_{N,Z}\le
10H_1\left|\Omega^{(1)}\right|^{1/2}
\Biggl(
\sum\limits_{g^{(1)}_{3},g^{(2)}_{3}\in\Om^{(3)}\atop\theta^{(1)},\theta^{(2)}\in Z}
\1_{\{\|z\|_{1,2}\le \frac{1}{2H_1}\}}
\Biggr)^{1/2},
\end{gather}
where $\|x\|_{1,2}=\max\{\|x_1\|,\|x_2\|\}$ for $x=(x_1,x_2)\in\rr^2,$  and
\begin{gather}\label{14-20}
z=g_3^{(1)}\theta^{(1)}-g_3^{(2)}\theta^{(2)}
\end{gather}
\begin{proof}
Applying the relation $|x|^2=x\overline{x}$ and reversing orders we easily obtain that
\begin{gather}\label{14-22}
\sum_{g_1\in\Z^2}\s\left(\frac{g_1}{H_1}\right)
\left|\sum_{\theta\in Z}\xi(\theta)\sum_{g_{3}\in\Om^{(3)}}e(g_1g_3\theta)\right|^2\le
\sum_{g^{(1)}_{3},g^{(2)}_{3}\in\Om^{(3)}}\sum_{\theta^{(1)},\theta^{(2)}\in Z}
\left|\sum_{g_1\in\Z^2}\s\left(\frac{g_1}{H_1}\right)e(g_1z)\right|.
\end{gather}
By application of the Poisson summation formula ~\cite[§4.3.]{IK}:
\begin{gather*}
\sum_{n\in\Z^2}f(n)=\sum_{k\in\Z^2}\hat{f}(k),
\end{gather*}
and writing $f(n)=\s\left(\frac{n}{H_1}\right)e(nz),$ we transform the inner sum in the right side of \eqref{14-22}:
\begin{gather}\label{14-23}
\sum_{g_1\in\Z^2}\s\left(\frac{g_1}{H_1}\right)e(g_1z)=
\sum_{k\in\Z^2}\int_{x\in\rr^2}\s\left(\frac{x}{H_1}\right)e(x(z-k))dx=H_1^2\sum_{k\in\Z^2}\hat{\s}((k-z)H_1).
\end{gather}
We note that the relation \eqref{14-23} can be obtained directly from ~\cite[(4.25)]{IK}. As $\hat{\s}(x,y)\neq0$ only if  $|x|\le\frac{1}{2}$ and $|y|\le\frac{1}{2},$ so the sum in the right side of \eqref{14-23} consists of at most one summand, hence,
\begin{gather}\label{14-24}
\left|\sum_{g_1\in\Z^2}\s\left(\frac{g_1}{H_1}\right)e(g_1z)\right|\le
36H_1^2\1_{\{\|z\|_{1,2}\le \frac{1}{2H_1}\}}.
\end{gather}
Substituting \eqref{14-24} into \eqref{14-22}, we obtain
\begin{gather}\label{14-25}
\sum_{g_1\in\Z^2}\s\left(\frac{g_1}{H_1}\right)
\left|\sum_{\theta\in Z}\xi(\theta)\sum_{g_{3}\in\Om^{(3)}}e(g_1g_3\theta)\right|^2\le
36H_1^2
\sum\limits_{g^{(1)}_{3},g^{(2)}_{3}\in\Om^{(3)}\atop\theta^{(1)},\theta^{(2)}\in Z}
\1_{\{\|z\|_{1,2}\le \frac{1}{2H_1}\}}
\end{gather}
Substituting \eqref{14-25} into \eqref{14-15}, we obtain \eqref{14-19}. This completes the proof of the lemma.
\end{proof}
\end{Le}
To transform the right side of \eqref{14-19}, we are to specify the set $Z.$ It follows from the Dirichlet theorem that for any
$\theta\in[0,1]$ there exist $a,q\in\N\cup\{0\}$ and $\beta\in\rr,$ such that
\begin{gather}\label{14-26}
\theta=\frac{a}{q}+\beta,\;(a,q)=1,\; 0\le a\le q\le\frac{N^{1/2}}{10A},\;\beta=\frac{K}{N},\; |K|\le\frac{10AN^{1/2}}{q},
\end{gather}
and $a=0$ or $a=q$ only if $q=1.$ We denote
\begin{gather}\label{14-27}
P_{Q_1,Q}^{(\beta)}=\left\{
\theta=\frac{a}{q}+\beta\;\Bigl|\;(a,q)=1,\;  0\le a\le q,\;Q_1\le q\le Q
\right\}.
\end{gather}
In what follows we always have $Z\subseteq P_{Q_1,Q}^{(\beta)}$ for some $Q_1,Q,\beta.$ We will write numbers $\theta^{(1)},\theta^{(2)}\in P_{Q_1,Q}^{(\beta)}$ in the following way
\begin{gather}\label{14-28}
\theta^{(1)}=\frac{a^{(1)}}{q^{(1)}}+\beta,\quad \theta^{(2)}=\frac{a^{(2)}}{q^{(2)}}+\beta.
\end{gather}
Let
\begin{gather}\label{14-29}
\NN_0=\left\{ (g^{(1)}_{3},g^{(2)}_{3},\theta^{(1)},\theta^{(2)})\in
\Om^{(3)}\times\Om^{(3)}\times Z^2\Bigl|\,
\left\|g^{(1)}_3\theta^{(1)}-g_3^{(2)}\theta^{(2)}\right\|_{1,2}\le \frac{1}{2H_1},
\right\}
\end{gather}
\begin{gather}\label{14-30}
\NN=\left\{ (g^{(1)}_{3},g^{(2)}_{3},\theta^{(1)},\theta^{(2)})\in
\Om^{(3)}\times\Om^{(3)}\times Z^2\Bigl|\,
\eqref{14-31}\, \mbox{and}\, \eqref{14-32}\, \mbox{hold}
\right\},
\end{gather}
where
\begin{gather}\label{14-31}
\|g^{(1)}_3\frac{a^{(1)}}{q^{(1)}}-g_3^{(2)}\frac{a^{(2)}}{q^{(2)}}\|_{1,2}\le \frac{74A^2\KK}{M^{(1)}},
\end{gather}
\begin{gather}\label{14-32}
|g^{(1)}_3-g_3^{(2)}|_{1,2}\le\min\left\{
\frac{73A^2N}{M^{(1)}},\; \frac{73A^2N}{M^{(1)}\overline{K}}+\frac{N}{\overline{K}}
\left\|g^{(1)}_3\frac{a^{(1)}}{q^{(1)}}-g_3^{(2)}\frac{a^{(2)}}{q^{(2)}}\right\|_{1,2}
\right\},
\end{gather}
and $\KK=\max\{1,|K|\}.$ We note that it follows from Lemma \ref{lemma-14-2} and the definition \eqref{14-29} that
\begin{gather}\label{14-33}
\sigma_{N,Z}\le
10H_1\left|\Omega^{(1)}\right|^{1/2}
\left|\NN_0\right|^{1/2}.
\end{gather}
\begin{Le}\label{lemma-14-3}
For $Z\subseteq P_{Q_1,Q}^{(\beta)}$ and for $M^{(1)}$ in \eqref{13-42} such that $M^{(1)}>146A^2\KK,$ the following inequality holds
\begin{gather}\label{14-34}
\sigma_{N,Z}\le
10H_1\left|\Omega^{(1)}\right|^{1/2}
\left|\NN\right|^{1/2}.
\end{gather}
\begin{proof}
In view of \eqref{14-33}, it is sufficient to prove that $\NN_0\subseteq \NN.$ Using \eqref{13-47-1}, one has
\begin{gather}\label{14-35}
|g^{(1)}_3-g_3^{(2)}|_{1,2}\le\frac{73A^2N}{M^{(1)}}.
\end{gather}
Hence,
\begin{gather}\label{14-36}
|(g_3^{(1)}-g_3^{(2)})\beta|_{1,2}\le\frac{73A^2N}{M^{(1)}}\frac{\KK}{N}\le\frac{1}{2},
\end{gather}
so
\begin{gather}\label{14-37}
|(g_3^{(1)}-g_3^{(2)})\beta|_{1,2}=\|(g_3^{(1)}-g_3^{(2)})\beta\|_{1,2},
\end{gather}
and
\begin{gather}\label{14-38}
\|(g_3^{(1)}-g_3^{(2)})\beta\|_{1,2}\le\frac{73A^2\KK}{M^{(1)}}.
\end{gather}
It follows from the definition \eqref{14-29} and the bound \eqref{14-38} that
\begin{gather}\label{14-39}
\|g_3^{(1)}\frac{a^{(1)}}{q^{(1)}}-g_3^{(2)}\frac{a^{(2)}}{q^{(2)}}\|_{1,2}\le \frac{73A^2\KK}{M^{(1)}}+\frac{1}{2H_1}<\frac{74A^2\KK}{M^{(1)}},
\end{gather}
that is, for $(g^{(1)}_{3},g^{(2)}_{3},\theta^{(1)},\theta^{(2)})\in\NN_0$ the inequality \eqref{14-31} holds. Using \eqref{14-37} and \eqref{14-29}, we obtain
\begin{gather}\label{14-40}
|(g_3^{(1)}-g_3^{(2)})\frac{K}{N}|_{1,2}\le\frac{1}{H_1}+
\|g_3^{(1)}\frac{a^{(1)}}{q^{(1)}}-g_3^{(2)}\frac{a^{(2)}}{q^{(2)}}\|_{1,2},
\end{gather}
whence
\begin{gather}\label{14-41}
|g_3^{(1)}-g_3^{(2)}|_{1,2}\le
\frac{N}{M^{(1)}\K}+\frac{N}{\K}
\left\|g_3^{(1)}\frac{a^{(1)}}{q^{(1)}}-g_3^{(2)}\frac{a^{(2)}}{q^{(2)}}\right\|_{1,2}.
\end{gather}
The inequality  \eqref{14-32} follows from the estimates \eqref{14-41} and \eqref{14-35}. This completes the proof of the lemma.
\end{proof}
\end{Le}
We denote
\begin{gather}\label{14-42}
\M=\left\{ (g^{(1)}_{3},g^{(2)}_{3},\theta^{(1)},\theta^{(2)})\in
\Om^{(3)}\times\Om^{(3)}\times Z^2\Bigl|\,
\eqref{14-43}\, \mbox{и}\, \eqref{14-44}\, \mbox{hold}
\right\},
\end{gather}
where
\begin{gather}\label{14-43}
|g_3^{(1)}-g_3^{(2)}|_{1,2}\le \frac{73A^2N}{M^{(1)}\KK},
\end{gather}
\begin{gather}\label{14-44}
\|g_3^{(1)}\frac{a^{(1)}}{q^{(1)}}-g_3^{(2)}\frac{a^{(2)}}{q^{(2)}}\|_{1,2}=0.
\end{gather}
Now we transform the equation \eqref{14-44}. Let $\q=[q^{(1)},q^{(2)}],$ then \eqref{14-44} can be written as
\begin{gather}\label{15-4}
\left(g_3^{(1)}\frac{a^{(1)}q^{(2)}}{(q^{(1)},q^{(2)})}-g_3^{(2)}\frac{a^{(2)}q^{(1)}}{(q^{(1)},q^{(2)})}\right)_{1,2}\equiv 0\pmod{\q}.
\end{gather}
By setting $g_3^{(1)}=(x_1,x_2)^{t},\,g_3^{(2)}=(y_1,y_2)^{t}$ in \eqref{15-4} we obtain the congruence
\begin{gather}\label{15-5}
x_1\frac{a^{(1)}q^{(2)}}{(q^{(1)},q^{(2)})}\equiv y_1\frac{a^{(2)}q^{(1)}}{(q^{(1)},q^{(2)})}\pmod{\q}
\end{gather}
and, the same one for $x_2,y_2.$ But $(a^{(1)},\frac{q^{(1)}}{(q^{(1)},q^{(2)})})\le(a^{(1)},q^{(1)})=1$ and, therefore,
\begin{gather*}
x_1\equiv 0 \pmod{\frac{q^{(1)}}{(q^{(1)},q^{(2)})}},\quad
x_2\equiv 0 \pmod{\frac{q^{(1)}}{(q^{(1)},q^{(2)})}}
\end{gather*}
and, the same one for $y_1,y_2.$ At the same time $(x_1,x_2)=(y_1,y_2)=1$ as the component of the vectors $g_3^{(1)},\, g_3^{(2)},$ thus $q^{(1)}=(q^{(1)},q^{(2)})=q^{(2)}=\q.$ So
\begin{gather}\label{14-42-1}
\M=\left\{ (g^{(1)}_{3},g^{(2)}_{3},\theta^{(1)},\theta^{(2)})\in
\Om^{(3)}\times\Om^{(3)}\times Z^2\Bigl|\,
\eqref{14-43}\, \mbox{и}\, \eqref{14-44-1}\, \mbox{hold}
\right\},
\end{gather}
where
\begin{gather}\label{14-44-1}
(a^{(1)}g_3^{(1)}-a^{(2)}g_3^{(2)})_{1,2}\equiv 0\pmod{\q}\quad\mbox{and}\quad
\q=q^{(1)}=q^{(2)}.
\end{gather}
\begin{Le}\label{lemma-14-4}
Let $M^{(1)}$ be such that for any $\theta^{(1)},\theta^{(2)}\in Z$ the following inequality holds
\begin{gather}\label{14-45}
[q^{(1)},q^{(2)}]<\frac{M^{(1)}}{74A^2\KK}.
\end{gather}
Then the following bound holds
\begin{gather}\label{14-46}
\sigma_{N,Z}\le
10H_1\left|\Omega^{(1)}\right|^{1/2}
\left|\M\right|^{1/2}.
\end{gather}
\begin{proof}
In view of \eqref{14-34} it is sufficient to prove that $\NN\subseteq \M.$ We note that to prove this it is sufficient to obtain that under the hypotheses of Lemma \ref{lemma-14-4} and $(g^{(1)}_{3},g^{(2)}_{3},\theta^{(1)},\theta^{(2)})\in\NN$ the relation \eqref{14-44} holds. It follows from \eqref{14-31} and \eqref{14-45} that
\begin{gather}\label{14-47}
\|g_3^{(1)}\frac{a^{(1)}}{q^{(1)}}-g_3^{(2)}\frac{a^{(2)}}{q^{(2)}}\|_{1,2}\le \frac{74A^2\KK}{M^{(1)}}<
\frac{1}{[q^{(1)},q^{(2)}]},
\end{gather}
this implies \eqref{14-44}. This completes the proof of the lemma.
\end{proof}
\end{Le}

Thus we reduced the problem of estimating $\sigma_{N,Z}$ to the evaluation the cardinality of one of the sets $\NN,\,\M.$ Let state one more lemma of a general nature. A similar statement was used by S.V.\,Konyagin in ~\cite[ 17]{Konyagin}.
\begin{Le}\label{lemma-14-5}
Let $W$ be a finite subset of the interval $[0,1]$ and let $|W|>10.$ Let $f:W\rightarrow \rr_{+}$ be a function such that, for any subset $Z\subseteq W$ the following bound holds
\begin{gather*}
\sum_{\theta\in Z}f(\theta)\le C_1|Z|^{1/2},
\end{gather*}
where $C_1$ is a non-negative constant not depending on the set $Z.$ Then the following estimate holds
\begin{gather}\label{14-49}
\sum_{\theta\in W}f^2(\theta)\ll C_1^2\log|W|
\end{gather}
with the absolute constant in Vinogradov symbol.
\begin{proof}
Let number the value of $f(\theta)$ in the decreasing order
\begin{gather*}
f_1\ge\ldots\ge f_{|W|}>0.
\end{gather*}
Then for any $k$ such that $1\le k\le|W|,$ one has
\begin{gather*}
kf_k\le\sum_{n=1}^{k}f_n\le 2C_1k^{1/2}
\end{gather*}
and, hence $f_k\le 2C_1k^{-1/2}.$ Thus
\begin{gather*}
\sum_{\theta\in W}f^2(\theta)=\sum_{n=1}^{|W|}f_n^2\le 8C_1^2\log|W|.
\end{gather*}
This completes the proof of the lemma.
\end{proof}
\end{Le}

\section{Consequences of general estimates of exponential sum.}\label{mu=3}
Let $Q_1,Q,\beta$ be given. For any $q$ in $Q_1\le q\le Q$ we define by any means the number $a_q,$ such that $(a_q,q)=1,\,  0\le a_q\le q.$ Let denote
\begin{gather}\label{15-0}
Z^{*}=\left\{
\theta=\frac{a_q}{q}+\beta\;\Bigl|\;Q_1\le q\le Q
\right\},
\end{gather}
then $|Z^{*}|\le Q.$ The following trivial bound holds
\begin{gather}\label{15-16-1}
\sum_{\theta\in P_{Q_1,Q}^{(\beta)}}\left|S_N(\theta)\right|^2\le
Q\sum_{Q_1\le q\le Q}\max\limits_{1\le a\le q, (a,q)=1}\left|S_N(\frac{a}{q}+\frac{K}{N})\right|^2=Q\sum_{\theta\in Z^*}\left|S_N(\theta)\right|^2,
\end{gather}
where as $a_q$ we have chosen numerators for which the maximum is achieved. We write $\gamma=1-\delta$ and
\begin{gather*}
Q_0=\max\left\{\exp\left(\frac{10^5A^4}{\epsilon_0^2}\right),\exp(\epsilon_0^{-5})\right\}.
\end{gather*}

\begin{Le}\label{lemma-15-2-1}
If
$
\KK^{5/2}Q^{3}\le N^{1-2\epsilon_0},\; \KK Q\ge Q_0,
$
then for any $Z\subseteq Z^{*}$ then the following bound holds
\begin{gather}\label{15-100}
\sigma_{N,Z}\ll
\frac{|\Omega_N||Z|}{|Z|^{1/2}\KK Q_1}(\KK^{5/2}Q^3)^{\gamma}(\KK^{12}Q^{15})^{\epsilon_0}.
\end{gather}
\begin{proof}
We put
\begin{gather}\label{15-101}
M^{(1)}=76A^2\KK Q^2,\quad
M^{(4)}=\KK^{1/2}(M^{(1)})^{2\epsilon_0},\quad
M^{(6)}=\KK Q.
\end{gather}
It follows from the statement of the lemma that all conditions of the Theorem \ref{theorem13-6} and Lemma \ref{lemma-14-4} hold. We prove that
\begin{gather}\label{15-106}
\left|\M\right|\ll|Z|\left|\Omega^{(3)}\right|(M^{(4)}M^{(1)})^{4\epsilon_0}\left|\Omega^{(5)}\right|\left(
\frac{(M^{(6)})^{1+2\epsilon_0}}{Q_1}+1\right)\left(
\frac{(M^{(1)}M^{(6)})^{1+2\epsilon_0}}{M^{(1)}\KK Q_1}+1\right),
\end{gather}
where $\M$ is defined in \eqref{14-42-1}. Let fix $\q$, for which there are $|Z|$ choices, this gives the first factor in \eqref{15-106}. Then it follows from the conditions on the set $Z$ that $a^{(1)}=a^{(2)}$ and the congruence \eqref{14-44-1} can be simplified to
\begin{gather}\label{15-107}
(g_3^{(1)}-g_3^{(2)})_{1,2}\equiv 0\pmod{\q}.
\end{gather}
We choose and fix the vector $g_3^{(2)},$ for which there are $\left|\Omega^{(3)}\right|$ choices. This gives the second factor in \eqref{15-106}. Using Theorem \ref{theorem13-6} we can write the second vector in the following form $g_3^{(2)}=\gamma_{4}\gamma_{5}\gamma_{6}(0,1)^{t},$ where $\gamma_{i}\in\Omega^{(i)}, i=4,5,6.$ Let $g_3^{(1)}=(x_1,x_2)^{t},\,g_3^{(2)}=(y_1,y_2)^{t}.$ It follows from \eqref{14-43} that
\begin{gather}\label{15-108}
\left|\frac{x_1}{x_2}-\frac{y_1}{y_2}\right|\ll\frac{(M^{(1)})^{2\epsilon_0}}{\KK}
\end{gather}
Let
$
\gamma_4=
\begin{pmatrix}
a & b \\
c & d
\end{pmatrix}
$
then
\begin{gather}\label{15-109}
\left|\frac{y_1}{y_2}-\frac{b}{d}\right|\ll\frac{1}{d^2}
\end{gather}
It follows from \eqref{13-62} that
\begin{gather}\label{15-150}
\frac{1}{d^2}=\frac{1}{\|\gamma_4\|^2}\ll\frac{(M^{(1)})^{4\epsilon_0}}{(M^{(4)})^2}.
\end{gather}
Using \eqref{15-108}, \eqref{15-109}, \eqref{15-150}, we obtain $\left|\frac{x_1}{x_2}-\frac{b}{d}\right|\ll\frac{(M^{(1)})^{2\epsilon_0}}{\KK}+\frac{(M^{(1)})^{4\epsilon_0}}{(M^{(4)})^2}.$
It follows from \eqref{13-62} that
\begin{gather*}\label{15-110}
\left|\frac{b_1}{d_1}-\frac{b_2}{d_2}\right|\ge\frac{1}{d_1d_2}\gg\frac{(M^{(1)})^{4\epsilon_0}}{(M^{(4)})^{2+4\epsilon_0}}
\end{gather*}
Hence, the number of different $\frac{b}{d},$ that is the number of matrices $\gamma_4$ is less than
\begin{gather}\label{15-120}
\frac{(M^{(4)})^{2+4\epsilon_0}}{(M^{(1)})^{4\epsilon_0}}\left(
\frac{(M^{(1)})^{2\epsilon_0}}{\KK}+\frac{(M^{(1)})^{4\epsilon_0}}{(M^{(4)})^2}\right)=
(M^{(4)})^{4\epsilon_0}\left(1+\frac{(M^{(4)})^2}{(M^{(1)})^{2\epsilon_0}\KK}\right)\ll
(M^{(4)})^{4\epsilon_0}(M^{(1)})^{3\epsilon_0}
\end{gather}
this is the third factor in \eqref{15-106}. We choose and fix the matrix $\gamma_{5}$ for which there are $\left|\Omega^{(5)}\right|$ choices. This gives the fourth factor in \eqref{15-106}. Let $\gamma_{6}(0,1)^{t}=(x,y)^{t},$ then \eqref{14-43}, \eqref{15-107} can be written in the following form
\begin{gather*}\label{15-110}
|x_1-(ax+by)|\le \frac{73A^2N}{M^{(1)}\KK},\quad
|x_2-(cx+dy)|\le \frac{73A^2N}{M^{(1)}\KK},
\end{gather*}
\begin{gather}\label{15-121}
x_1\equiv(ax+by)\pmod{\q},\quad
x_2\equiv(cx+dy)\pmod{\q},
\end{gather}
where $g_3^{(1)}=(x_1,x_2)^{t},$ and $
\gamma_4\gamma_5=
\begin{pmatrix}
a & b \\
c & d
\end{pmatrix}
$
are fixed. As $\det(\gamma_4\gamma_5)=\pm1,$ we obtain from \eqref{15-121} that
\begin{gather}\label{15-122}
x\equiv\pm(dx_1-bx_2)\pmod{\q},\quad
y\equiv\pm(ax_2-cx_1)\pmod{\q},
\end{gather}
Hence, applying \eqref{13-63} and \eqref{15-121}, we obtain that for $x$ there are less than $\left(\frac{(M^{(6)})^{1+2\epsilon_0}}{\q}+1\right)$ choices. This gives the fifth factor in \eqref{15-106}. For $y$ there are less than $\left(\frac{73A^2N}{M^{(1)}\KK\q d}+1\right)$ choices. Using \eqref{13-64}, we obtain that it is less than
\begin{gather*}\label{15-110}
\frac{73A^2N}{M^{(1)}\KK Q_1}\frac{25000A^5(M^{(1)}M^{(6)})^{1+2\epsilon_0}}{N}+1\ll
\frac{(M^{(1)}M^{(6)})^{1+2\epsilon_0}}{M^{(1)}\KK Q_1}+1
\end{gather*}
This gives the sixth factor in \eqref{15-106} Thus the bound \eqref{15-106} is proved. Substituting \eqref{15-106} into \eqref{14-46} and applying \eqref{15-101}, we obtain
\begin{gather*}
\sigma_{N,Z}\ll|Z|^{1/2}
\left|\Omega^{(1)}\right|^{1/2}
\left|\Omega^{(3)}\right|^{1/2}\left|\Omega^{(5)}\right|^{1/2}
(M^{(1)})^{1+4\epsilon_0}(M^{(4)})^{2\epsilon_0}
\frac{(M^{(1)})^{1/2+\epsilon_0}(M^{(6)})^{1+2\epsilon_0}}{(M^{(1)})^{1/2}\KK^{1/2} Q_1},
\end{gather*}
Since
\begin{gather}\label{15-112}
\left|\Omega^{(1)}\right|^{1/2}
\left|\Omega^{(3)}\right|^{1/2}\left|\Omega^{(5)}\right|^{1/2}=
\left|\Omega^{(1)}\right|^{1/2}
\left|\Omega^{(4)}\right|^{1/2}\left|\Omega^{(5)}\right|\left|\Omega^{(6)}\right|^{1/2}=
\frac{|\Omega_N|}{\left|\Omega^{(1)}\right|^{1/2}
\left|\Omega^{(4)}\right|^{1/2}\left|\Omega^{(6)}\right|^{1/2}},
\end{gather}
and applying Theorem \ref{theorem13-5} we obtain
\begin{gather}\label{15-151}
\sigma_{N,Z}\ll|Z|^{1/2}|\Omega_N|
\frac{(M^{(1)})^{1+5\epsilon_0}(M^{(4)})^{2\epsilon_0}(M^{(6)})^{1+2\epsilon_0}}{
(M^{(1)}M^{(4)}M^{(6)})^{\delta-\epsilon_0/2}
\KK^{1/2} Q_1}.
\end{gather}
Substituting $M^{(i)}$ from \eqref{15-101} into \eqref{15-151}, one has
\begin{gather*}
\sigma_{N,Z}\ll|Z|^{1/2}|\Omega_N|
\KK^{5\gamma/2-1+12\epsilon_0}Q^{3\gamma-1+15\epsilon_0}\frac{Q}{Q_1}.
\end{gather*}
This completes the proof of the lemma.
\end{proof}
\end{Le}
\begin{Le}\label{lemma-15-2-2}
If
$
\KK^{5/2}Q^{3}\le N^{1-2\epsilon_0}, \KK Q\ge Q_0,
$
then the following estimate holds
\begin{gather}\label{15-132}
\sum_{\theta\in P_{Q_1,Q}^{(\beta)}}\left|S_N(\theta)\right|^2\ll
\frac{|\Omega_N|^2Q}{\KK^2Q_1^2}(\KK^{5}Q^6)^{\gamma}(\KK^{24}Q^{30})^{\epsilon_0}.
\end{gather}
\begin{proof}
Applying Lemma \ref{lemma-14-5} with $W=Z^{*},$ we obtain from \eqref{15-100} that
\begin{gather*}
\sum_{\theta\in Z^*}\left|S_N(\theta)\right|^2\ll
|\Omega_N|^2\KK^{5\gamma-2+24\epsilon_0}Q^{6\gamma-2+30\epsilon_0}\frac{Q^2}{Q^2_1}.
\end{gather*}
Using the trivial bound \eqref{15-16-1}, we obtain the desired estimate \eqref{15-132}. Lemma is proved.
\end{proof}
\end{Le}

\begin{Le}\label{lemma-15-2-3}
Let
\begin{gather*}
1\le a_{i}\le q_{i}\le Q,\; (a_i,q_i)=1,\; i=1,2;\quad   \frac{a_1}{q_1}\neq\frac{a_2}{q_2},
\end{gather*}
then for $R>\max\left\{\frac{Q^2}{Y_1},4Y,2Q\right\}$ there are no coprime numbers $y_1,y_2,$ satisfying the following conditions
\begin{gather}\label{15-113}
\|y_i\left(\frac{a_1}{q_1}-\frac{a_2}{q_2}\right)\|<\frac{1}{R},\quad
0<Y_1\le y_i\le Y,\;i=1,2.
\end{gather}
\begin{proof}
We define $\alpha=\left|\frac{a_1}{q_1}-\frac{a_2}{q_2}\right|.$ Assume the contrary, then there exists integers $n_1,n_2$ such that
\begin{gather}\label{15-114}
|y_i\alpha-n_i|<\frac{1}{R},\quad i=1,2.
\end{gather}
First we prove that $n_1,n_2\neq0.$ If $n_1=0,$ then $\alpha<\frac{1}{Ry_1}\le\frac{1}{RY_1}.$ On the other side, $\alpha\ge\frac{1}{q_1q_2}\ge\frac{1}{Q^2}.$ And we obtain $\frac{1}{Q^2}<\frac{1}{RY_1},$ but this is contradictory to the conditions of the lemma.
It follows from \eqref{15-114} that there exist $-1<\theta_i<1, i=1,2$ such that
\begin{gather}\label{15-115}
y_i\alpha=n_i+\frac{\theta_i}{R},\quad i=1,2.
\end{gather}
Hence
$
(y_1n_2-y_2n_1)\alpha=\frac{\theta_1n_2-\theta_2n_1}{R}
$
and so one has
\begin{gather}\label{15-116}
|y_1n_2-y_2n_1|<\frac{|n_2|+|n_1|}{\alpha R}<
\frac{y_2\alpha+y_1\alpha}{\alpha R}+\frac{2}{\alpha R^2}<\frac{2Y}{R}+\frac{2Q^2}{R^2}<1
\end{gather}
Thus we obtain that $y_1n_2=y_2n_1.$ But $(y_1,y_2)=1$ and so one has $n_1=y_1,\;n_2=y_2.$ Substituting the last equalities into  \eqref{15-115} we obtain
\begin{gather}\label{15-116}
|1-\alpha|<\frac{1}{RY_1}<\frac{1}{Q^2}.
\end{gather}
But this is impossible. Lemma is proved.
\end{proof}
\end{Le}

\begin{Le}\label{lemma-15-3}
If $\KK q\ge Q_0,$ then the following bound holds
\begin{gather}\label{15-22}
|S_N(\theta)|\ll\frac{|\Omega_N|}{\KK q}N^{\gamma}(\KK^2q^2N)^{\epsilon_0}.
\end{gather}
\begin{proof}
We use the partition of the ensemble $\Omega_N$ given by  Theorem \ref{theorem13-3}. We put $Z=\{\theta\}$ and
\begin{gather}\label{15-23}
M^{(1)}=76A^2\KK q.
\end{gather}
Then conditions of Lemma \ref{lemma-14-4} hold and therefore the estimate \eqref{14-46} is valid. Thus we obtain
\begin{gather}\notag
\left|\M\right|\le
\sum_{g_3^{(2)}\in\Omega^{(3)}}\sum_{g_3^{(1)}\in\Z^2}\1_{\{
g_3^{(1)}\equiv g_3^{(2)} \pmod{q},\,|g_3^{(1)}-g_3^{(2)}|_{1,2}\le \frac{N}{\KK^2q}
\}}\le\\\le|\Omega^{(3)}|\left(1+\frac{N}{\KK^2q^2}\right)^2\ll\frac{|\Omega^{(3)}|N^2}{\KK^4q^4},\label{15-25}
\end{gather}
as, in view of \eqref{14-26} we have $\KK q\le N^{1/2}.$ Substituting \eqref{15-25} into \eqref{14-46}, we obtain
\begin{gather}\label{15-26}
|S_N(\theta)|\ll H_1\left|\Omega^{(1)}\right|^{1/2}|\Omega^{(3)}|^{1/2}\frac{N}{\KK^2q^2}\ll
|\Omega_N|\frac{(M^{(1)})^{1+2\epsilon_0}}{|\Omega_N|^{1/2}}\frac{N}{\KK^2q^2}.
\end{gather}
Using \eqref{15-23} and the lower bound of \eqref{13-39}, we have
\begin{gather*}
|S_N(\theta)|\ll|\Omega_N|\frac{(\KK q)^{1+2\epsilon_0}N^{1-\delta+\epsilon_0}}{\KK^2q^2}.
\end{gather*}
This completes the proof of the lemma.
\end{proof}
\end{Le}
\begin{Le}\label{lemma-15-4}
Let the following inequalities hold
$$N^{\epsilon_0/2}\le Q^{1/2}\le Q_1\le Q\le\frac{N^{1/2}}{10A},\,\KK Q\le N^{1/2+\epsilon_0}.$$
Then for any $Z\subseteq P_{Q_1,Q}^{(\beta)}$ the following bound holds
\begin{gather}\label{15-27-2}
\sum_{\theta\in Z}|S_N(\theta)|\ll
|Z|^{1/2}|\Omega_N|\left(
\frac{N^{\gamma+4\epsilon_0}}{(\KK Q_1)^{1/2}}+
\frac{N^{5\gamma/4-1/4+5\epsilon_0}}{\KK^{1/2}}+N^{\gamma-1/4+3\epsilon_0}\right).
\end{gather}
\begin{proof}
We use the partition of the ensemble $\Omega_N$ given by Theorem \ref{theorem13-3}. We put $Z\subseteq P_{Q_1,Q}^{(\beta)}$ and
\begin{gather}\label{15-28}
M^{(1)}=150A^2N^{3/4+\epsilon_0}.
\end{gather}
Then $M^{(1)}\ge150A^2N^{1/2+\epsilon_0}\ge150A^2\KK$ and conditions of Lemma \ref{lemma-14-3} hold. Therefore, the inequality \eqref{14-34} holds.
Let $g_3^{(1)}=(x_1,x_2)^{t},\,g_3^{(2)}=(y_1,y_2)^{t},$ then it follows from Theorem \ref{theorem13-3} that
\begin{gather}\label{15-31}
\max\{y_1,y_2\}=y_2\le\frac{73A^2N}{M^{(1)}}.
\end{gather}
We denote
\begin{gather}\label{15-32}
Y=
\begin{pmatrix}
x_1 & y_1 \\
x_2 & y_2
\end{pmatrix},\quad
\mathcal{Y}=det(Y)=x_1y_2-y_1x_2.
\end{gather}
Then it follows from the triangle inequality that
\begin{gather}\label{15-33}
\|\mathcal{Y}\frac{a^{(1)}}{q^{(1)}}\|\le
\|y_2(x_1\frac{a^{(1)}}{q^{(1)}}-y_1\frac{a^{(2)}}{q^{(2)}})\|+\|y_1(y_2\frac{a^{(2)}}{q^{(2)}}-x_2\frac{a^{(1)}}{q^{(1)}})\|.
\end{gather}
Applying \eqref{14-31},\,\eqref{15-31} and  \eqref{15-28} we similarly estimate both summands in the right side of \eqref{15-33}
\begin{gather}\label{15-34}
\|\mathcal{Y}\frac{a^{(1)}}{q^{(1)}}\|\le
2\frac{74A^2\KK}{M^{(1)}}\frac{73A^2N}{M^{(1)}}<\frac{1}{2Q}<\frac{1}{q^{(1)}}.
\end{gather}
Hence, $\mathcal{Y}\equiv 0 \pmod{q^{(1)}}$ and similarly $\mathcal{Y}\equiv 0 \pmod{q^{(2)}}.$ So one has
$\mathcal{Y}\equiv 0 \pmod{\q},$ where $\q=[q^{(1)},q^{(2)}]$ and so $Q_1\le\q\le Q^2.$
The set $\NN$ can be represented as a union of the sets $\M_1,\M_2.$ For the first set $\mathcal{Y}=0,$ for the second one $\mathcal{Y}\neq0.$ To prove the estimate \eqref{15-27-2} we use the following bounds
\begin{gather}\label{15-35}
|\M_1|\ll|Z||\Omega^{(3)}|,
\end{gather}
\begin{gather}\label{15-36}
|\M_2|\ll
|Z|\left(\frac{N}{M^{(1)}}\right)^2\frac{1}{\KK Q_1}
\left(|\Omega^{(3)}|N^{2\epsilon_0}+Q^{1+\epsilon_0}\right).,
\end{gather}
which will be proved below in Lemma \ref{lemma-15-5} and \ref{lemma-15-6} respectively. Hence,
\begin{gather}\label{15-36-9}
|\NN|^{1/2}\ll|Z|^{1/2}\left(
\frac{N^{1+\epsilon_0}}{M^{(1)}}\frac{|\Omega^{(3)}|^{1/2}}{\KK^{1/2}Q^{1/2}_1}+
\frac{NQ^{\epsilon}}{M^{(1)}\KK^{1/2}}+|\Omega^{(3)}|^{1/2}\right).
\end{gather}
Substituting \eqref{15-36-9} into \eqref{14-34} we obtain
\begin{gather}\label{15-36-10}
\sigma_{N,Z}\ll
(M^{(1)})^{1+2\epsilon_0}\left|\Omega^{(1)}\right|^{1/2}|Z|^{1/2}\left(
\frac{N^{1+\epsilon_0}}{M^{(1)}}\frac{|\Omega^{(3)}|^{1/2}}{\KK^{1/2}Q^{1/2}_1}+
\frac{NQ^{\epsilon_0}}{M^{(1)}\KK^{1/2}}+|\Omega^{(3)}|^{1/2}
\right).
\end{gather}
Using the bounds \eqref{13-55} and \eqref{13-39}, we have
\begin{gather}\label{15-36-4}
|\Omega^{(1)}|^{1/2}|\Omega^{(3)}|^{1/2}=|\Omega_N|\frac{1}{|\Omega_N|^{1/2}}\le
|\Omega_N|\frac{1}{N^{\delta-\epsilon_0/2}}.
\end{gather}
\begin{gather}\label{15-36-11}
\left|\Omega^{(1)}\right|^{1/2}=|\Omega_N|\frac{|\Omega^{(1)}|^{1/2}}{|\Omega_N|}\le
|\Omega_N|\frac{(M^{(1)})^{\delta+2\epsilon_0}}{N^{2\delta-\epsilon_0}}.
\end{gather}
Substituting \eqref{15-36-4} and \eqref{15-36-11} into \eqref{15-36-10}, we obtain
\begin{gather*}\label{15-36-12}
\sigma_{N,Z}\ll
|\Omega_N||Z|^{1/2}\Bigl(
\frac{N^{1-\delta+2\epsilon_0}(M^{(1)})^{2\epsilon_0}}{\KK^{1/2}Q^{1/2}_1}+
N^{1-2\delta+\epsilon_0}(M^{(1)})^{\delta+4\epsilon_0}\frac{Q^{\epsilon_0}}{\KK^{1/2}}+
\frac{(M^{(1)})^{1+2\epsilon_0}}{N^{\delta-\epsilon_0/2}}\Bigr).
\end{gather*}
Substituting $M^{(1)}$ from \eqref{15-28} one has \eqref{15-27-2}. This completes the proof of the lemma.
\end{proof}
\end{Le}
\begin{Le}\label{lemma-15-5}
Under the hypotheses of Lemma \ref{lemma-15-4} one has
\begin{gather}\label{15-37}
|\M_1|\ll|Z||\Omega^{(3)}|.
\end{gather}
\begin{proof}
To simplify we put $R=(M^{(1)})^{1+2\epsilon_0}.$
Recall that $M^{(1)}$ is defined in \eqref{15-28}. It follows from $\mathcal{Y}=0$ that
$\frac{x_1}{x_2}=\frac{y_1}{y_2}.$ Since
\begin{gather*}
g_3^{(1)}=(x_1,x_2)^{t},\,g_3^{(2)}=(y_1,y_2)^{t},
\end{gather*}
one has $(x_1,x_2)=1,(y_1,y_2)=1$ and hence $x_1=y_1,x_2=y_2.$ In particular we have ${|\beta(g_3^{(1)}-g_3^{(2)})|_{1,2}=0.}$ Then it follows from  \eqref{14-29} that
\begin{gather}\label{15-38}
\|y_1\left(\frac{a^{(1)}}{q^{(1)}}-\frac{a^{(2)}}{q^{(2)}}\right)\|\le
\|y_1\left(\theta^{(1)}-\theta^{(2)}\right)\|+0\le\frac{1}{H_1}<\frac{1}{R}
\end{gather}
and similarly $\|y_2\left(\frac{a^{(1)}}{q^{(1)}}-\frac{a^{(2)}}{q^{(2)}}\right)\|\le\frac{1}{H_1}<\frac{1}{R}.$ Thus we have
\begin{gather}\label{15-39}
|\M_1|\le \sum_{g_3^{(2)}\in\Om^{(3)}}\sum_{\frac{a^{(1)}}{q^{(1)}}\in Z}\sum_{\frac{a^{(2)}}{q^{(2)}}\in Z}
\1_{\{\|y_{1,2}\left(\frac{a^{(1)}}{q^{(1)}}-\frac{a^{(2)}}{q^{(2)}}\right)\|<\frac{1}{R}\}}.
\end{gather}
Let prove that all conditions of Lemma \ref{lemma-15-2-3} hold. It follows from Theorem \ref{theorem13-3} that
\begin{gather*}
Y_1\le y_i\le Y, \mbox{where}\quad
Y_1=\frac{N}{150A^2(M^{(1)})^{1+2\epsilon_0}},\;Y=\frac{73A^2N}{M^{(1)}}.
\end{gather*}
Under the hypotheses of Lemma \ref{lemma-15-4} one has
$R>\max\left\{\frac{Q^2}{Y_1},4Y,2Q\right\},$ and so we can apply Lemma \ref{lemma-15-2-3}. We obtain that in \eqref{15-39} only summands corresponding to $\frac{a^{(1)}}{q^{(1)}}=\frac{a^{(2)}}{q^{(2)}}$ are nonzero. Hence
$|\M_1|\ll|Z||\Omega^{(3)}|.$ This completes the proof of the lemma.
\end{proof}
\end{Le}
\begin{Le}\label{lemma-15-6}
Under the hypotheses of Lemma \ref{lemma-15-4} one has
\begin{gather}\label{15-46}
|\M_2|\ll
|Z|\left(\frac{73A^2N}{M^{(1)}}\right)^2\frac{1}{\KK Q_1}
\left(|\Omega^{(3)}|N^{2\epsilon_0}+Q^{1+\epsilon_0}\right).
\end{gather}
\begin{proof}
Since $|\mathcal{Y}|\le\max\{x_1y_2,x_2y_1\}\le x_2y_2,$ then applying the bound \eqref{15-31} to each factor we obtain
\begin{gather}\label{15-48}
|\mathcal{Y}|\le\left(\frac{73A^2N}{M^{(1)}}\right)^2.
\end{gather}
As $\q|\mathcal{Y}$ and $\mathcal{Y}\neq0,$ one has
\begin{gather}\label{15-49}
\q\le\min\{Q^2,|\mathcal{Y}|\}\le\min\{Q^2,\left(\frac{73A^2N}{M^{(1)}}\right)^2\}\le
\sqrt{Q^2\left(\frac{73A^2N}{M^{(1)}}\right)^2}=Q\frac{73A^2N}{M^{(1)}}.
\end{gather}
It follows from the conditions of the lemma and from \eqref{15-28} that
\begin{gather}\label{15-49-1}
\q\frac{74A^2\KK}{M^{(1)}}\le Q\frac{73A^2N}{M^{(1)}}\frac{74A^2\KK}{M^{(1)}}<1.
\end{gather}
By \eqref{15-49-1} we obtain that the right side of \eqref{14-31} is less than $\frac{1}{\q},$ whence
\begin{gather}\label{15-50}
\|g_3^{(1)}\frac{a^{(1)}}{q^{(1)}}-g_3^{(2)}\frac{a^{(2)}}{q^{(2)}}\|_{1,2}=0.
\end{gather}
It follows from the equation \eqref{15-50} that the congruence \eqref{15-4} holds, from which we deduced \eqref{14-44-1}, that is
\begin{gather}\label{15-51}
q^{(1)}=q^{(2)}=\q,\quad (g_3^{(1)}a^{(1)}-g_3^{(2)}a^{(2)})_{1,2}\equiv 0\pmod{\q}.
\end{gather}
Then by \eqref{14-32} we obtain
\begin{gather}\label{15-51-1}
|g_3^{(1)}-g_3^{(2)}|_{1,2}\le\frac{73A^2N}{M^{(1)}\KK}.
\end{gather}
To simplify we put $U=\frac{73A^2N}{M^{(1)}}.$ Note that if $g_3^{(1)},\,g_3^{(2)},\,a^{(1)}$ are fixed then $a^{(2)}$ is determined uniquely by \eqref{15-51}. Actually, let $$g_3^{(1)}=(x_1,x_2)^{t},\,g_3^{(2)}=(y_1,y_2)^{t},$$ then it follows from \eqref{15-51} that
\begin{gather}\label{15-57}
x_1a^{(1)}\equiv y_1a^{(2)}\pmod{\q},\,x_2a^{(1)}\equiv y_2a^{(2)}\pmod{\q}.
\end{gather}
Since $(a^{(1)},\q)=(a^{(1)},\q)=1,$ we have $$\delta_1=(x_1,\q)=(y_1,\q),\,\delta_2=(x_2,\q)=(y_2,\q)$$ and moreover $(\delta_1,\delta_2)=1,$ as $(x_1,x_2)=1.$ Then one can obtain from \eqref{15-57} that
\begin{gather}\label{15-58}
a^{(2)}\equiv A_1a^{(1)}\pmod{\frac{\q}{\delta_1}},\,a^{(2)}\equiv A_2a^{(2)}\pmod{\frac{\q}{\delta_2}}.
\end{gather}
These two congruences are equivalent to the congruence modulo $[\frac{\q}{\delta_1},\frac{\q}{\delta_2}]=\q,$ and hence $a^{(2)}$ is uniquely determined. Since $\q|\mathcal{Y},$ one has $x_1y_2\equiv x_2y_1\pmod{\q}$ and thus
\begin{gather}\label{15-59}
|\M_2|\le|Z|\sum_{g_3^{(1)}\in\Omega^{(3)}}\sum_{g_3^{(2)}\in\Omega^{(3)}\atop |g_3^{(1)}-g_3^{(2)}|_{1,2}\le\frac{U}{\KK}}\1_{\{x_1y_2\equiv x_2y_1\pmod{\q}\}}
\end{gather}
The change of variables $z_1=x_1-y_1,\, z_2=x_2-y_2$ leads to
\begin{gather}\label{15-60}
|\M_2|\le|Z|\sum_{g_3^{(1)}\in\Omega^{(3)}}\sum_{|z_{1,2}|\le\frac{U}{\KK}}\1_{\{x_1z_2\equiv x_2z_1\pmod{\q}\}}
\end{gather}
We consider three cases.
\begin{enumerate}
  \item Let $z_1>0,\,z_2>0.$ We fix the vector $g_3^{(1)}\in\Omega^{(3)},$ then $x_1z_2-x_2z_1=j\q.$ Let estimate the amount of $j.$ We have
$$x_1-x_2\frac{U}{\KK}\le j\q\le x_1\frac{U}{\KK}-x_2$$
and, hence, $\#j\le \frac{U^2}{\q\KK}+1.$ For a fixed $j$ the solution of the congruence is given by the formulae
\begin{gather*}
z_1=z_{1,0}+nx_1,\;z_2=z_{2,0}+nx_2
\end{gather*}
In view of $x_2\gg\frac{U}{(M^{(1)})^{2\epsilon_0}},\,|z_{1,2}|\le\frac{U}{\KK}$ we have $\#n\ll\frac{(M^{(1)})^{2\epsilon_0}}{\KK}+1.$ Thus
\begin{gather}\label{15-61}
\sum_{g_3^{(1)}\in\Omega^{(3)}}\sum_{0<z_{1,2}\le\frac{T}{\KK}}\1_{\{x_1z_2\equiv x_2z_1\pmod{\q}\}}\ll
|\Omega^{(3)}|\left(\frac{U^2}{\q\KK}+1\right)\left(\frac{(M^{(1)})^{2\epsilon_0}}{\KK}+1\right)
\end{gather}
It follows from he conditions of Lemma \ref{lemma-15-4} that $U^2>\q\KK,$ so one has
\begin{gather}\label{15-62}
\sum_{g_3^{(1)}\in\Omega^{(3)}}\sum_{0<z_{1,2}\le\frac{U}{\KK}}\1_{\{x_1z_2\equiv x_2z_1\pmod{\q}\}}\ll
|\Omega^{(3)}|\frac{U^2}{\q\KK}N^{2\epsilon_0}.
\end{gather}
  \item Let $z_1>0,\,z_2<0.$ In the same way as in the previous case we obtain
\begin{gather}\label{15-63}
\sum_{g_3^{(1)}\in\Omega^{(3)}}\sum_{0<-z_2,z_{1}\le\frac{U}{\KK}}\1_{\{x_1z_2\equiv x_2z_1\pmod{\q}\}}\ll
|\Omega^{(3)}|\frac{U^2}{\q\KK}N^{2\epsilon_0}.
\end{gather}
  \item Let$z_1=0.$ One has
\begin{gather}\label{15-64}
\sum_{g_3^{(1)}\in\Omega^{(3)}}\sum_{|z_2|\le\frac{U}{\KK}}\1_{\{x_1z_2\equiv0\pmod{\q}\}}\le
\sum_{g_3^{(1)}\in\Omega^{(3)}}\left(\frac{U}{q\KK}(x_1,q)+1\right)\le\\\le
|\Omega^{(3)}|+\frac{U^2}{q\KK}\sum_{x_1\le T}(x_1,q)
\end{gather}
Next
\begin{gather*}
\sum_{x_1\le U}(x_1,q)\le\sum_{d|q}d\left(\frac{U}{d}+1\right)\ll_{\epsilon}Uq^{\epsilon}+q^{1+\epsilon}
\end{gather*}
and so
\begin{gather}\label{15-65}
\sum_{g_3^{(1)}\in\Omega^{(3)}}\sum_{|z_2|\le\frac{T}{\KK}}\1_{\{x_1z_2\equiv0\pmod{\q}\}}\ll_{\epsilon}
|\Omega^{(3)}|+\frac{U^2}{q\KK}\left(Uq^{\epsilon}+q^{1+\epsilon}\right).
\end{gather}
\end{enumerate}
Using \eqref{15-62},\,\eqref{15-63},\,\eqref{15-65} and putting $\epsilon=\epsilon_0$ we obtain
\begin{gather}\label{15-66}
|\M_2|\ll|Z|\frac{U^2}{q\KK}\left(Uq^{\epsilon_0}+q^{1+\epsilon_0}+|\Omega^{(3)}|N^{2\epsilon_0}\right)\le
|Z|\frac{U^2}{\KK Q_1}\left(Q^{1+\epsilon}+|\Omega^{(3)}|N^{2\epsilon_0}\right).
\end{gather}
Lemma is proved.
\end{proof}
\end{Le}

\begin{Sl}\label{sled15-1}
Under the hypotheses of Lemma \ref{lemma-15-4}, one has
\begin{gather}\label{15-74-0}
\sum_{\theta\in P_{Q_1,Q}^{(\beta)}}|S_N(\theta)|^2\ll
|\Omega_N|^2C_1^2Q^{\epsilon_0},
\end{gather}
where
\begin{gather}\label{15-75}
C_1=\frac{N^{\gamma+4\epsilon_0}}{(\KK Q_1)^{1/2}}+
\frac{N^{5\gamma/4-1/4+5\epsilon_0}}{\KK^{1/2}}+N^{\gamma-1/4+3\epsilon_0}
\end{gather}
\begin{proof}
It was proved in Lemma \ref{lemma-15-4} that for any $Z\subseteq P_{Q_1,Q}^{(\beta)}$ one has
$
\sum_{\theta\in Z}|S_N(\theta)|\ll
|\Omega_N||Z|^{1/2}C_1.
$
Applying Lemma \ref{lemma-14-5} with $W=P_{Q_1,Q}^{(\beta)},f(\theta)=\frac{|S_N(\theta)|}{|\Omega_N|},$ we obtain \eqref{15-74-0}. This completes the proof of the corollary.
\end{proof}
\end{Sl}

\section{Replacing integrals by sums.}
The purpose of the following reasonings is a slight modification of the results of §\ref{Trigsums1}. It follows from the statement of Lemma \ref{lemma-17-1} that we need to know how to estimate the following expression
\begin{gather}\label{12-1}
\frac{1}{N}\mathop{{\sum}^*}_{0\le a\le q\le X}\int\limits_{|K|\le Y}
\left|S_N(\frac{a}{q}+\frac{K}{N})\right|^2dK,
\end{gather}
where $Y$ may depend on $q.$ The following reasonings are similar to \cite[Lemma 26 p.145]{Korobov}. Let take a sufficiently large number $T,$ then
\begin{gather}\label{12-2}
\int\limits_{|K|\le Y}\left|S_N(\frac{a}{q}+\frac{K}{N})\right|^2dK\le
\sum_{|l|\le TY}\int\limits_{l/T}^{(l+1)/T}\left|S_N(\frac{a}{q}+\frac{K}{N})\right|^2dK.
\end{gather}
Hence, $K=\frac{l}{T}+\lambda, \theta=\frac{a}{q}+\frac{K}{N}=\frac{a}{q}+\frac{l}{TN}+\frac{\lambda}{N}$ and
\begin{gather}\label{12-3}
\left|S_N(\theta)-S_N(\frac{a}{q}+\frac{l}{TN})\right|\le\lambda|\Omega_{N}|\Rightarrow
\left|S_N(\theta)\right|^2\le\left|S_N(\frac{a}{q}+\frac{l}{TN})\right|^2+\lambda^2|\Omega_{N}|^2.
\end{gather}
So one has
\begin{gather}\label{12-3-1}
\frac{1}{N}\mathop{{\sum}^*}_{0\le a\le q\le X}\int\limits_{|K|\le Y}
\left|S_N(\frac{a}{q}+\frac{K}{N})\right|^2dK\le\frac{1}{N}
\mathop{{\sum}^*}_{0\le a\le q\le X}\sum_{|l|\le TY}\left(\frac{1}{T}\left|S_N(\frac{a}{q}+\frac{l}{TN})\right|^2+
\frac{|\Omega_{N}|^2}{T^3}\right)
\end{gather}
Choosing $T$ sufficiently large we obtain that the investigation of the expression of the form \eqref{12-1} reduce to the investigation of the quantity
\begin{gather}
\frac{1}{TN}
\mathop{{\sum}^*}_{0\le a\le q\le X}\sum_{|l|\le TY}\left|S_N(\frac{a}{q}+\frac{l}{TN})\right|^2.
\end{gather}
Our next purpose is to modify Lemma \ref{lemma-14-3} and Lemma \ref{lemma-14-4}. Let
\begin{gather}\label{12-4-2}
P_{Q_1,Q}^{\Ll_1,\Ll}=\left\{
\theta=\frac{a}{q}+\frac{l}{TN}\;\Bigl|\;(a,q)=1,\;  0\le a\le q,\;Q_1\le q\le Q,\; \Ll_1\le|l|\le\Ll
\right\},
\end{gather}
and $\LL_1=\frac{\Ll_1}{T},\;\LL=\frac{\Ll}{T}.$ We note that $\LL$ is similar to $\KK$ and that $\left|P_{Q_1,Q}^{\Ll_1,\Ll}\right|\ll\Ll Q^2=\LL TQ^2. $
Let $Z\subseteq P_{Q_1,Q}^{\Ll_1,\Ll},$ we denote
\begin{gather}\label{14-30-2}
\NNN=\left\{ (g^{(1)}_{3},g^{(2)}_{3},\theta^{(1)},\theta^{(2)})\in
\Om^{(3)}\times\Om^{(3)}\times Z^2\Bigl|\,
\eqref{14-31-2}\, \mbox{and}\, \eqref{14-32-2}\, \mbox{hold}
\right\},
\end{gather}
where
\begin{gather}\label{14-31-2}
\|g_3^{(1)}\frac{a^{(1)}}{q^{(1)}}-g_3^{(2)}\frac{a^{(2)}}{q^{(2)}}\|_{1,2}\le \frac{74A^2\LL}{M^{(1)}},
\end{gather}
\begin{gather}\label{14-32-2}
|g_3^{(1)}\frac{l_1}{TN}-g_3^{(2)}\frac{l_2}{TN}|_{1,2}\le\min\left\{
\frac{73A^2\LL}{M^{(1)}},\; \frac{73A^2}{M^{(1)}}+
\left\|g_3^{(1)}\frac{a^{(1)}}{q^{(1)}}-g_3^{(2)}\frac{a^{(2)}}{q^{(2)}}\right\|_{1,2}
\right\},
\end{gather}

\begin{gather}\label{14-42-2}
\MM=\left\{ (g^{(1)}_{3},g^{(2)}_{3},\theta^{(1)},\theta^{(2)})\in
\Om^{(3)}\times\Om^{(3)}\times Z^2\Bigl|\,
\eqref{14-43-2}\, \mbox{and}\, \eqref{14-44-2}\, \mbox{hold}
\right\},
\end{gather}
where
\begin{gather}\label{14-43-2}
|g_3^{(1)}\frac{l_1}{TN}-g_3^{(2)}\frac{l_2}{TN}|_{1,2}\le \frac{1}{M^{(1)}},
\end{gather}
\begin{gather}\label{14-44-2}
\|g_3^{(1)}\frac{a^{(1)}}{q^{(1)}}-g_3^{(2)}\frac{a^{(2)}}{q^{(2)}}\|_{1,2}=0.
\end{gather}
The following lemma can de proved in the same manner as Lemma \ref{lemma-14-3}.
\begin{Le}\label{lemma-14-3-2}
For $Z\subseteq P_{Q_1,Q}^{\Ll_1,\Ll}$ and for $M^{(1)}$ in \eqref{13-42} such that $M^{(1)}>146A^2\LL,$ the following inequality holds
\begin{gather}\label{14-34-2}
\sigma_{N,Z}\ll
H_1\left|\Omega^{(1)}\right|^{1/2}
\left|\NNN\right|^{1/2}.
\end{gather}
\end{Le}

The following lemma can de proved in the same manner as Lemma \ref{lemma-14-4}.
\begin{Le}\label{lemma-14-4-2}
Let $Z\subseteq P_{Q_1,Q}^{\Ll_1,\Ll}$ and let ensemble $\Omega_N$ satisfy the hypotheses of Lemma \ref{lemma-14-4}.
Let $M^{(1)}$ be such that for any $\theta^{(1)},\theta^{(2)}\in Z$ the following inequality holds
\begin{gather}\label{14-45-2}
M^{(1)}>76A^2[q^{(1)},q^{(2)}]\LL.
\end{gather}
Then the following bound holds
\begin{gather}\label{14-46-2}
\sigma_{N,Z}\ll H_1\left|\Omega^{(1)}\right|^{1/2}
\left|\MM\right|^{1/2}.
\end{gather}
\end{Le}
For any $q$ in $Q_1\le q\le Q$ we define by any means the number $a_q,$ such that $(a_q,q)=1,\,  0\le a_q\le q.$ Let denote
\begin{gather}\label{15-0}
Z^{*}=
\left\{
\theta=\frac{a_q}{q}+\frac{l}{TN}\;\Bigl|\;Q_1\le q\le Q,\; \Ll_1\le|l|\le\Ll
\right\},
\end{gather}
then $|Z^{*}|\le\LL TQ.$ The following statement is similar to Lemma \ref{lemma-15-2-1}.
\begin{Le}\label{lemma-16-2-1}
If
$
\LL^{3/2}Q^{3}\le N^{1-2\epsilon_0},\;\LL\ge Q^{10\epsilon_0},\; \LL Q\ge Q_0,
$
then for any $Z\subseteq Z^{*}$ then the following bound holds
\begin{gather}\label{16-100}
\sigma_{N,Z}\ll|Z|^{1/2}T^{1/2}|\Omega_N|
\LL^{3\gamma/2-1/2+10\epsilon_0}Q^{3\gamma-1+12\epsilon_0}\frac{Q\LL^{\delta/2}}{Q_1\LL_1^{\delta/2}}.
\end{gather}
\begin{proof}
We put
\begin{gather}\label{16-101}
M^{(1)}=76A^2\LL Q^2,\quad
M^{(4)}=\LL_1^{1/2}(M^{(1)})^{\epsilon_0},\quad
M^{(6)}=Q.
\end{gather}
It follows from the statement of the lemma that all conditions of the Theorem \ref{theorem13-6} and Lemma \ref{lemma-14-4-2} hold. We prove that
\begin{gather}\label{16-106}
\left|\M\right|\ll|Z|T(M^{(1)})^{2\epsilon_0}
\left|\Omega^{(3)}\right|(M^{(4)})^{4\epsilon_0}\left|\Omega^{(5)}\right|\left(
\frac{(M^{(6)})^{1+2\epsilon_0}}{Q_1}+1\right)^2.
\end{gather}
In the same manner as in Lemma \ref{lemma-15-2-1} one can prove that $q^{(1)}=(q^{(1)},q^{(2)})=q^{(2)}=\q.$
We choose and fix $\theta_1$ (thus $q_1=\q,\;l_1$ are fixed) for which there are $|Z|$ choices, this gives the first factor in \eqref{16-106}. Then it follows from the conditions on the set $Z$ that $a^{(1)}=a^{(2)}$ and the congruence \eqref{14-44-2} can be simplified to
\begin{gather}\label{16-107}
(g_3^{(1)}-g_3^{(2)})_{1,2}\equiv 0\pmod{\q}.
\end{gather}
Using the first inequality in \eqref{14-43-2} and the lower bound \eqref{13-47-1} we obtain that there are less than $T(M^{(1)})^{2\epsilon_0}$ choices for $l_2$ (if $g_3^{(1)}\,g_3^{(2)},\,l_1$ have been already fixed).
This gives the second factor in \eqref{16-106}. We choose and fix the vector $g_3^{(1)}$ for which there are $\left|\Omega^{(3)}\right|$ choices. This gives the third factor in \eqref{16-106}. Using Theorem \ref{theorem13-6} we can write the second vector in the following form $g_3^{(2)}=\gamma_{4}\gamma_{5}\gamma_{6}(0,1)^{t},$ where $\gamma_{i}\in\Omega^{(i)}, i=4,5,6.$ Let $g_3^{(1)}=(x_1,x_2)^{t},\,g_3^{(2)}=(y_1,y_2)^{t}.$ It follows from \eqref{14-43-2} that
\begin{gather}\label{16-108}
\left|\frac{x_1}{x_2}-\frac{y_1}{y_2}\right|\ll\frac{(M^{(1)})^{2\epsilon_0}}{\LL_1}
\end{gather}
In the same manner as in Lemma \ref{lemma-15-2-1} we obtain that the number of matrices $\gamma_4$ is less then
$(M^{(4)})^{4\epsilon_0}.$ This is the fourth factor in \eqref{16-106}. We choose and fix the matrix $\gamma_{5}$ for which there are $\left|\Omega^{(5)}\right|$ choices. This gives the fifth factor in \eqref{16-106}. Let $\gamma_{6}(0,1)^{t}=(x,y)^{t},$ then\eqref{16-107} can be written in the form
\begin{gather*}\label{16-110}
x_1\equiv(ax+by)\pmod{\q},\quad
x_2\equiv(cx+dy)\pmod{\q},
\end{gather*}
where $g_3^{(1)}=(x_1,x_2)^{t},$ and
$\gamma_4\gamma_5=
\begin{pmatrix}
a & b \\
c & d
\end{pmatrix}
$
have been already fixed. Using \eqref{13-63} we obtain that for $x$ and for $y$ there are less than $\left(\frac{(M^{(6)})^{1+2\epsilon_0}}{\q}+1\right)$ choices. This gives the sixth factor in \eqref{16-106}. Thus the estimate \eqref{16-106} is proved. Substituting \eqref{16-106} into \eqref{14-46-2} and applying \eqref{16-101}, we obtain
\begin{gather*}
\sigma_{N,Z}\ll|Z|^{1/2}T^{1/2}
\left|\Omega^{(1)}\right|^{1/2}
\left|\Omega^{(3)}\right|^{1/2}\left|\Omega^{(5)}\right|^{1/2}
(M^{(1)})^{1+2\epsilon_0}(M^{(4)})^{2\epsilon_0}
\frac{(M^{(6)})^{1+2\epsilon_0}}{Q_1}.
\end{gather*}
Using \eqref{15-112} and Theorem \ref{theorem13-5} one has
\begin{gather}\label{16-111}
\sigma_{N,Z}\ll|Z|^{1/2}T^{1/2}|\Omega_N|
\frac{(M^{(1)}M^{(6)})^{1+2\epsilon_0}(M^{(4)})^{2\epsilon_0}}{(M^{(1)}M^{(4)}M^{(6)})^{\delta-\epsilon_0/2}Q_1}
\end{gather}
Substituting $M^{(i)}$ from \eqref{16-101} into \eqref{16-111}, we have
\begin{gather*}\label{16-111}
\sigma_{N,Z}\ll|Z|^{1/2}T^{1/2}|\Omega_N|
\LL^{3\gamma/2-1/2+10\epsilon_0}Q^{3\gamma-1+12\epsilon_0}\frac{Q\LL^{\delta}}{Q_1\LL_1^{\delta}}.
\end{gather*}
Lemma is proved.
\end{proof}
\end{Le}
\begin{Le}\label{lemma-16-2-2}
If
$
\LL^{3/2}Q^{3}\le N^{1-2\epsilon_0},\,\LL\ge Q^{10\epsilon_0}, \LL Q\ge Q_0,
$
then the following estimate holds
\begin{gather}\label{16-112}
\sum_{\theta\in P_{Q_1,Q}^{\Ll_1,\Ll}}\left|S_N(\theta)\right|^2\ll
|\Omega_N|^2T
\LL^{3\gamma-1+20\epsilon_0}Q^{6\gamma-1+24\epsilon_0}\frac{Q^2\LL^{\delta}}{Q_1^2\LL_1^{\delta}}.
\end{gather}
\begin{proof}
Applying Lemma \ref{lemma-14-5} with $W=Z^{*},$ we obtain from \eqref{16-100} that
\begin{gather*}
\sum_{\theta\in Z^*}\left|S_N(\theta)\right|^2\ll
|\Omega_N|^2T\LL^{3\gamma-1+30\epsilon_0}Q^{6\gamma-2+24\epsilon_0}\frac{Q^2\LL^{\delta}}{Q_1^2\LL_1^{\delta}}.
\end{gather*}
Using the trivial bound \eqref{15-16-1}, we obtain the desired estimate \eqref{16-112}. Lemma is proved.
\end{proof}
\end{Le}
The following statement is similar to Lemma \ref{lemma-15-4}.
\begin{Le}\label{lemma-16-4}
Let the following inequalities hold
$N^{\epsilon_0}\le Q\le N^{1/2-10\epsilon_0},\,\LL Q\le N^{1/2+\epsilon_0}.$ Then for any $Z\subseteq P_{Q_1,Q}^{\Ll_1,\Ll}$ the following bound holds
\begin{gather}\notag
\sum_{\theta\in Z}|S_N(\theta)|\ll
|\Omega_N||Z|^{1/2}T^{1/2}\Bigl(
\frac{N^{\gamma+2\epsilon_0}(\LL Q)^{2\epsilon_0}}{(\LL_1Q_1)^{1/2}}+
N^{\gamma-1/2+2\epsilon_0}(\LL Q)^{1/2+2\epsilon_0}+\\+
N^{3\gamma/2-1/4+4\epsilon_0}\frac{(\LL Q)^{1/4-\gamma/2+3\epsilon_0}Q^{\epsilon_0}}{Q_1^{1/2}}+
N^{3\gamma/2-3/4+4\epsilon_0}(\LL Q)^{3/4-\gamma/2+3\epsilon_0}Q^{1/2+\epsilon_0}
\Bigr).\label{16-27-2}
\end{gather}
\begin{proof}
We use the partition of the ensemble $\Omega_N$ given by Theorem \ref{theorem13-3}. We put $Z\subseteq P_{Q_1,Q}^{\Ll_1,\Ll}$ and
\begin{gather}\label{16-28}
M^{(1)}=150A^2N^{1/2}(\LL Q)^{1/2}.
\end{gather}
Then $M^{(1)}\ge150A^2N^{1/2+\epsilon_0}\ge150A^2\KK$ and conditions of Lemma \ref{lemma-14-3-2} hold. Therefore, the inequality \eqref{14-34-2} holds.
Let $g_3^{(1)}=(x_1,x_2)^{t},\,g_3^{(2)}=(y_1,y_2)^{t},\,\mathcal{Y}=det(Y)=x_1y_2-y_1x_2.$ In the same way as in Lemma \ref{lemma-15-4} we obtain that $\mathcal{Y}\equiv 0 \pmod{\q},$ where $\q=[q^{(1)},q^{(2)}].$ The set $\NNN$ can be represented as a union of the sets $\M_1,\M_2.$ For the first set $\mathcal{Y}=0,$ for the second one $\mathcal{Y}\neq0.$ To prove the estimate \eqref{16-27-2} we use the following bounds
\begin{gather}\label{16-35}
|\M_1|\ll|Z|T(M^{(1)})^{2\epsilon_0}|\Omega^{(3)}|,
\end{gather}
\begin{gather}\label{16-36}
|\M_2|\ll_{\epsilon}|Z|T(M^{(1)})^{2\epsilon_0}\left(\frac{|\Omega^{(3)}|(M^{(3)})^2}{\LL_1Q_1}+
\frac{(M^{(3)})^3Q^{\epsilon}}{Q_1}+M^{(3)}Q^{1+\epsilon}\right),
\end{gather}
which will be proved below in Lemma \ref{lemma-16-5} and \ref{lemma-16-6} respectively. Hence,
\begin{gather}\label{16-36-9}
|\NNN|^{1/2}\ll|Z|^{1/2}T^{1/2}(M^{(1)})^{\epsilon_0}\left(
\frac{|\Omega^{(3)}|^{1/2}M^{(3)}}{(\LL_1Q_1)^{1/2}}+|\Omega^{(3)}|^{1/2}+
\frac{(M^{(3)})^{3/2}Q^{\epsilon_0}}{Q_1^{1/2}}+(M^{(3)})^{1/2}Q^{1/2+\epsilon_0}
\right).
\end{gather}
Substituting \eqref{16-36-9} into \eqref{14-34-2}, we obtain
\begin{gather*}\label{15-36-10}
\sigma_{N,Z}\ll|Z|^{1/2}T^{1/2}
(M^{(1)})^{1+3\epsilon_0}\left|\Omega^{(1)}\right|^{1/2}\left(
\frac{|\Omega^{(3)}|^{1/2}M^{(3)}}{(\LL_1Q_1)^{1/2}}+|\Omega^{(3)}|^{1/2}+
\frac{(M^{(3)})^{3/2}Q^{\epsilon_0}}{Q_1^{1/2}}+(M^{(3)})^{1/2}Q^{1/2+\epsilon_0}
\right).
\end{gather*}
Using the bounds \eqref{15-36-4}, \eqref{15-36-11} and the definition of $M^{(1)}$ and $M^{(3)}$, we have
\begin{gather*}\label{15-36-12}
\sigma_{N,Z}\ll
|\Omega_N||Z|^{1/2}T^{1/2}\Bigl(
\frac{N^{\gamma+2\epsilon_0}(\LL Q)^{2\epsilon_0}}{(\LL_1Q_1)^{1/2}}+
N^{\gamma-1/2+2\epsilon_0}(\LL Q)^{1/2+2\epsilon_0}+\\
N^{3\gamma/2-1/4+4\epsilon_0}\frac{(\LL Q)^{1/4-\gamma/2+3\epsilon_0}Q^{\epsilon_0}}{Q_1^{1/2}}+
N^{3\gamma/2-3/4+4\epsilon_0}(\LL Q)^{3/4-\gamma/2+3\epsilon_0}Q^{1/2+\epsilon_0}
\Bigr).
\end{gather*}
Lemma is proved.
\end{proof}
\end{Le}
\begin{Le}\label{lemma-16-5}
Under the hypotheses of Lemma \ref{lemma-16-4} one has
\begin{gather}\label{16-37}
|\M_1|\ll|Z|T(M^{(1)})^{2\epsilon_0}|\Omega^{(3)}|.
\end{gather}
\begin{proof}
To simplify we put $$R=\frac{M^{(1)}}{75A^2\LL},M^{(3)}=\frac{73A^2N}{M^{(1)}}.$$
It can be proved in the same way as in Lemma \ref{lemma-15-5} that  $x_1=y_1,x_2=y_2.$ Then it follows from \eqref{14-31-2} that
\begin{gather}\label{16-38}
\|y_1\left(\frac{a^{(1)}}{q^{(1)}}-\frac{a^{(2)}}{q^{(2)}}\right)\|<\frac{1}{R}
\end{gather}
and similarly $\|y_2\left(\frac{a^{(1)}}{q^{(1)}}-\frac{a^{(2)}}{q^{(2)}}\right)\|<\frac{1}{R}.$ It follows from \eqref{14-32-2} that
\begin{gather*}\label{16-38}
|y_{1,2}\frac{l_1-l_2}{TN}|\le\frac{73A^2}{M^{(1)}}+
\|y_{1,2}\left(\frac{a^{(1)}}{q^{(1)}}-\frac{a^{(2)}}{q^{(2)}}\right)\|
\end{gather*}
Thus we obtain
\begin{gather}\label{16-39}
|\M_1|\le \sum_{g_3^{(2)}\in\Om^{(3)}}\sum_{\theta^{(1)}\in Z}\sum_{\theta^{(2)}\in Z}
\1_{\{\|y_{1,2}\left(\frac{a^{(1)}}{q^{(1)}}-\frac{a^{(2)}}{q^{(2)}}\right)\|<\frac{1}{R}\}}
\1_{
\{
|y_{1,2}\frac{l_1-l_2}{T}|\le M^{(3)}+
N|y_{1,2}\left(\frac{a^{(1)}}{q^{(1)}}-\frac{a^{(2)}}{q^{(2)}}\right)|
\}
}.
\end{gather}
Let prove that all conditions of Lemma \ref{lemma-15-2-3} hold. It follows from Theorem \ref{theorem13-3} that
\begin{gather*}
Y_1\le y_i\le Y,\; \mbox{где}
Y_1=\frac{N}{150A^2(M^{(1)})^{1+2\epsilon_0}},\;Y=\frac{73A^2N}{M^{(1)}}.
\end{gather*}
Under the hypotheses of Lemma \ref{lemma-16-4} one has
$R>\max\left\{\frac{Q^2}{Y_1},4Y,2Q\right\},$ and so we can apply Lemma \ref{lemma-15-2-3}. We obtain that in \eqref{16-39} only summands corresponding to $\frac{a^{(1)}}{q^{(1)}}=\frac{a^{(2)}}{q^{(2)}}$ are nonzero. We choose and fix $\theta^{(1)}$ for which there are $Z$ choices. Thus $\frac{a^{(1)}}{q^{(1)}}=\frac{a^{(2)}}{q^{(2)}},\, l_1$ are fixed. Using the inequality
$
|y_{1,2}\frac{l_1-l_2}{T}|\le M^{(3)}
$
and the lower bound from \eqref{13-47-1} we obtain that for $l_2$ there are less than $T(M^{(1)})^{2\epsilon_0}$ choices. Hence
$|\M_1|\ll|Z|T(M^{(1)})^{2\epsilon_0}|\Omega^{(3)}|.$ Lemma is proved.
\end{proof}
\end{Le}
\begin{Le}\label{lemma-16-6}
Under the hypotheses of Lemma \ref{lemma-16-4} one has
\begin{gather}\label{16-46}
|\M_2|\ll_{\epsilon}|Z|T(M^{(1)})^{2\epsilon_0}\left(\frac{|\Omega^{(3)}|(M^{(3)})^2}{\LL_1Q_1}+
\frac{(M^{(3)})^3Q^{\epsilon}}{Q_1}+M^{(3)}Q^{1+\epsilon}\right).
\end{gather}
\begin{proof}
In the same way as in Lemma \ref{lemma-15-6} we obtain
\begin{gather}\label{16-51}
q^{(1)}=q^{(2)}=\q,\quad (g_3^{(1)}a^{(1)}-g_3^{(2)}a^{(2)})_{1,2}\equiv 0\pmod{\q}.
\end{gather}
Then it follows from \eqref{14-32-2} that
\begin{gather}\label{16-51-1}
|g_3^{(1)}\frac{l_1}{TN}-g_3^{(2)}\frac{l_2}{TN}|_{1,2}\le\frac{73A^2N}{M^{(1)}}.
\end{gather}
We choose and fix $\theta_1$  (that is $q_1=\q,\;l_1$ are fixed) for which there are $|Z|$ choices. Using the first inequality  \eqref{16-51-1} and the lower bound \eqref{13-47-1} we obtain that  ($g_3^{(1)},\,g_3^{(2)},\,l_1$ are fixed) there are less than $T(M^{(1)})^{2\epsilon_0}$ choices for $l_2$. As in Lemma \ref{lemma-15-6} we obtain that $a^{(2)}$ is determined uniquely by \eqref{16-51}. Hence
\begin{gather}\label{16-59}
|\M_2|\le|Z|T(M^{(1)})^{2\epsilon_0}\sum_{g_3^{(1)}\in\Omega^{(3)}}\sum_{g_3^{(2)}\in\Omega^{(3)}\atop |x_2\frac{l_1}{TN}-y_2\frac{l_2}{TN}|\le\frac{73A^2}{M^{(1)}}}\1_{\{x_1y_2\equiv x_2y_1\pmod{\q}\}}
\end{gather}
To simplify we put $M^{(3)}=\frac{73A^2N}{M^{(1)}}.$ Let
$$\delta_1=(x_1,q),\,\delta_2=(x_2,q),\,
x_3=\frac{x_1}{\delta_1},\,x_4=\frac{x_2}{\delta_2},\,y_3=\frac{y_1}{\delta_1},\,y_4=\frac{y_2}{\delta_2},\,
p=\frac{\q}{\delta_1\delta_2},$$
then
\begin{gather}\label{16-69}
|\M_2|\le|Z|T(M^{(1)})^{2\epsilon_0}\sum_{g_3^{(1)}\in\Omega^{(3)}}\sum_{y_3\le\frac{M^{(3)}}{\delta_1}}
\sum_{y_4\le\frac{M^{(3)}}{\delta_2}\atop |x_4\frac{l_1}{T}-y_4\frac{l_2}{T}|\le\frac{M^{(3)}}{\delta_2}}
\1_{\{x_3y_4\equiv x_4y_3\pmod{p}\}}.
\end{gather}
Since $(x_3,p)=(x_4,p)=1$ the congruence \eqref{16-69} can be written in the form  $y_4\equiv cy_3\pmod{p},$ where $c\equiv x_3^{-1}x_4\pmod{p}.$ Then using the function $\delta_p(a)$
\begin{gather*}\label{deltasum}
\delta_p(a)=\frac{1}{p}\sum_{x=1}^p e\left(\frac{ax}{p}\right)=
\left\{
              \begin{array}{ll}
                1, & \hbox{if $a\equiv 0 \pmod{p}$;} \\
                0, & \hbox{else,}
              \end{array}
\right.
\end{gather*}
we obtain
\begin{gather*}\label{16-70}
|\M_2|\le|Z|T(M^{(1)})^{2\epsilon_0}\sum_{g_3^{(1)}\in\Omega^{(3)}}\frac{1}{p}\sum_{k=1}^p
\sum_{y_3\le\frac{M^{(3)}}{\delta_1}}
\sum_{y_4\le\frac{M^{(3)}}{\delta_2}\atop |x_4\frac{l_1}{T}-y_4\frac{l_2}{T}|\le\frac{M^{(3)}}{\delta_2}}
e\left(k\frac{x_3y_4-x_4y_3}{p}\right)
\end{gather*}
In the same way as in~\cite[p.18]{Korobov} we obtain
\begin{gather}\label{16-70}
\frac{1}{p}\sum_{k=1}^p
\sum_{y_3\le\frac{M^{(3)}}{\delta_1}}
\sum_{y_4\le\frac{M^{(3)}}{\delta_2}\atop |x_4\frac{l_1}{T}-y_4\frac{l_2}{T}|\le\frac{M^{(3)}}{\delta_2}}
e\left(k\frac{x_3y_4-x_4y_3}{p}\right)\le
\frac{1}{p}\left(\frac{M^{(3)}}{\delta_1}+1\right)\left(\frac{M^{(3)}}{\delta_2\LL_1}+1\right)+
O(s(\frac{c}{p})\log^2p),
\end{gather}
where $s(\alpha)=\sum\limits_{1\le i \le s}a_i$ is the sum of partial quotients of the number $\alpha=[0;a_1,\ldots,a_s].$ Substituting \eqref{16-70} into \eqref{16-69} we have
\begin{gather}\notag
|\M_2|\le|Z|T(M^{(1)})^{2\epsilon_0}\Bigl(\sum_{g_3^{(1)}\in\Omega^{(3)}}
\frac{1}{q}\left(
\frac{(M^{(3)})^2}{\LL_1}+M^{(3)}\delta_2+
\frac{M^{(3)}\delta_1}{\LL_1}+\delta_1\delta_2\right)+\\+
\log^2Q\sum_{x_1\le M^{(3)}}\sum_{x_2\le M^{(3)}}s(\frac{x_3^{-1}x_2}{q/\delta_1})\Bigr).\label{16-71}
\end{gather}
Using the following result of Knuth and Yao ~\cite{KnuthYao},
$$\sum_{a\le b}s(a/b)\ll b\log^2 b,$$
we obtain
\begin{gather}\label{16-72}
\sum_{x_1\le M^{(3)}}\sum_{x_2\le M^{(3)}}s(\frac{x_3^{-1}x_2}{q/\delta_1})\le\sum_{x_1\le M^{(3)}}
\left(\frac{M^{(3)}}{q/\delta_1}+1\right)\frac{q}{\delta_1}\log^2 q\le((M^{(3)})^2+M^{(3)}q)\log^2 q.
\end{gather}
One has
\begin{gather}\label{16-73}
\sum_{x_1\le M^{(3)}}\delta_1\ll
\left(M^{(3)}+q\right)q^{\epsilon_0}.
\end{gather}
Substituting \eqref{16-72},\eqref{16-73} into \eqref{16-71} we have
\begin{gather*}\label{15-73}
|\M_2|\ll_{\epsilon}|Z|T(M^{(1)})^{2\epsilon_0}\left(\frac{|\Omega^{(3)}|(M^{(3)})^2}{\LL_1Q_1}+
\frac{(M^{(3)})^3Q^{\epsilon}}{Q_1}+M^{(3)}Q^{1+\epsilon}\right).
\end{gather*}
Lemma is proved.
\end{proof}
\end{Le}

Applying Lemma \ref{lemma-14-5} we obtain
\begin{Le}\label{lemma-16-7}
Let
$N^{\epsilon_0}\le Q\le N^{1/2-10\epsilon_0},\,\LL Q\le N^{\frac{1}{2}+\epsilon_0},$
then the following bound holds
\begin{gather}\notag
\sum_{\theta\in P_{Q_1,Q}^{\Ll_1,\Ll}}\left|S_N(\theta)\right|^2\ll
|\Omega_N|^2T\Bigl(
\frac{N^{2\gamma+4\epsilon_0}(\LL Q)^{4\epsilon_0}}{\LL_1Q_1}+
N^{2\gamma-1+4\epsilon_0}(\LL Q)^{1+4\epsilon_0}+\\+
N^{3\gamma-1/2+8\epsilon_0}\frac{(\LL Q)^{1/2-\gamma+6\epsilon_0}Q^{2\epsilon_0}}{Q_1}+
N^{3\gamma-3/2+8\epsilon_0}(\LL Q)^{3/2-\gamma+6\epsilon_0}Q^{1+2\epsilon_0}
\Bigr).\label{16-74}
\end{gather}
\end{Le}


\section{Estimates for integrals of $|S_N(\theta)|^2$.}
We recall that the sequence $\{N_j\}_{-J-1}^{J+1}$ was defined in \eqref{10-8}. It will be convenient to amplify this sequence by the element $N_{1/2}=\frac{N^{1/2}}{10A}.$ We note that $N_0<N_{1/2}<N_{1}.$ Let $N^{*}_{1/2}=\frac{N}{N_{1/2}}.$
\begin{Le}\label{lemma-17-1}
The following inequality holds
\begin{gather}\label{17-1}
\int_0^1\left|S_N(\theta)\right|^2d\theta\le \frac{1}{N}
\mathop{{\sum}^*}_{0\le a\le q\le\QN}\int\limits_{|K|\le\frac{\KN}{q}}
\left|S_N(\frac{a}{q}+\frac{K}{N})\right|^2dK,
\end{gather}
where $\mathop{{\sum}^*}$ means that the sum is taken over $a$ and $q$ being coprime for $q\ge1,$ and $a=0,1$ for $q=1.$
\begin{proof}
It follows from the Dirichlet theorem that for any $\theta\in[0,1]$ there exist $a,q\in\N$ and $\beta\in\rr$ such that
\begin{gather*}
\theta=\frac{a}{q}+\beta,\;(a,q)=1,\; 0\le a\le q\le\QN,\;|\beta|\le\frac{1}{q\QN},
\end{gather*}
so
\begin{gather}\label{17-2}
\int_0^1\left|S_N(\theta)\right|^2d\theta=\int_0^1
\left|S_N(\frac{a}{q}+\beta)\right|^2d(\frac{a}{q}+\beta)\le
\mathop{{\sum}^*}_{0\le a\le q\le\QN}\int\limits_{|\beta|\le\frac{1}{q\QN}}
\left|S_N(\frac{a}{q}+\beta)\right|^2d\beta.
\end{gather}
The change of variables $K=N\beta$ in \eqref{17-2} leads to the inequality \eqref{17-1}. This completes the proof of the lemma.
\end{proof}
\end{Le}
Recall that
\begin{gather*}
Q_0=\max\left\{\exp\left(\frac{10^5A^4}{\epsilon_0^2}\right),\exp(\epsilon_0^{-5})\right\}.
\end{gather*}
\begin{Le}\label{lemma-17-2}
The following inequality holds
\begin{gather}\notag
\int_0^1\left|S_N(\theta)\right|^2d\theta\le
2Q_0^2\frac{|\Omega_N|^2}{N}+
\frac{1}{N}\mathop{{\sum}^*}_{0\le a\le q\le\QN\atop q>Q_0}
\int\limits_{|K|\le\frac{\KN}{q}}
\left|S_N(\frac{a}{q}+\frac{K}{N})\right|^2dK+\\
\frac{1}{N}\mathop{{\sum}^*}_{0\le a\le q\le Q_0}\int\limits_{\frac{Q_0}{q}\le|K|\le\frac{\KN}{q}}
\left|S_N(\frac{a}{q}+\frac{K}{N})\right|^2dK
\label{17-4}
\end{gather}
\begin{proof}
To simplify we write $f(K)=\left|S_N(\frac{a}{q}+\frac{K}{N})\right|^2,$ then
\begin{gather}\notag
\mathop{{\sum}^*}_{0\le a\le q\le\QN}\int\limits_{|K|\le\frac{\KN}{q}}f(K)dK=
\mathop{{\sum}^*}_{0\le a\le q\le\QN\atop q>Q_0}\int\limits_{|K|\le\frac{\KN}{q}}f(K)dK+\\
\mathop{{\sum}^*}_{0\le a\le q\le Q_0}\int\limits_{\frac{Q_0}{q}<|K|\le\frac{\KN}{q}}f(K)dK+
\mathop{{\sum}^*}_{1\le a\le q\le Q_0}\int\limits_{|K|\le\frac{Q_0}{q}}f(K)dK.\label{17-5}
\end{gather}
We estimate the fourth integral trivially
\begin{gather}\label{17-6}
\sum_{q\le Q_0}\mathop{{\sum}^*}_{0\le a\le q}\int\limits_{|K|\le\frac{Q_0}{q}}f(K)dK\le 2Q_0^2|\Omega_N|^2.
\end{gather}
Substituting \eqref{17-6} into \eqref{17-5} and using \eqref{17-1}, we obtain \eqref{17-4}. This completes the proof of the lemma.
\end{proof}
\end{Le}
First we estimate the third integral in the right side of \eqref{17-4}. It is convenient to use the following notation
\begin{gather}\label{17-7}
\gamma=1-\delta,\quad\xi_1=N^{2\gamma+10\epsilon_0}.
\end{gather}

\begin{Le}\label{lemma-17-3-2}
For $\gamma<\frac{1}{5}-5\epsilon_0$ and $\epsilon_0\in(0,\frac{1}{2500})$ the following inequality holds
\begin{gather}\label{17-7-8}
\frac{1}{N}\mathop{{\sum}^*}_{0\le a\le q\le Q_0}\int\limits_{\frac{Q_0}{q}\le|K|\le\xi_1}
\left|S_N(\frac{a}{q}+\frac{K}{N})\right|^2dK\ll\frac{|\Omega_N|^2}{N}.
\end{gather}
\begin{proof}
We denote by $I$ the integral in the left side of \eqref{17-7-8}. Applying Lemma \ref{lemma-15-2-2} we obtain
\begin{gather}\label{17-7-9}
I\ll|\Omega_N|^2
Q_0^{6\gamma+1+30\epsilon_0}\int\limits_{\frac{Q_0}{q}\le|K|\le\xi_1}\KK^{5\gamma-2+24\epsilon_0}dK\ll
|\Omega_N|^2,
\end{gather}
by the choice of the parameter $\gamma.$ Summing \eqref{17-7-9} over $0\le a\le q\le Q_0$ we obtain \eqref{17-7-8}. Lemma is proved.
\end{proof}
\end{Le}
The integral
\begin{gather}\label{17-7-8-1}
\frac{1}{N}\mathop{{\sum}^*}_{0\le a\le q\le Q_0}\int\limits_{\xi_1\le|K|\le\frac{\KN}{q}}
\left|S_N(\frac{a}{q}+\frac{K}{N})\right|^2dK.
\end{gather}
will be estimated later in Lemma \ref{lemma-17-3}.

It remains to estimate the first integral in the right side of \eqref{17-4}, that is,
\begin{gather}\label{17-7-10}
\frac{1}{N}\mathop{{\sum}^*}_{0\le a\le q\le\QN\atop q>Q_0}\int\limits_{\frac{Q_0}{q}\le|K|\le\frac{\KN}{q}}
\left|S_N(\frac{a}{q}+\frac{K}{N})\right|^2dK.
\end{gather}
The following lemmas will be devoted to this. We partition the range of summation and integration over $q,\,K$ into five subareas:

\begin{center}
  \includegraphics[width=450pt,height=350pt]{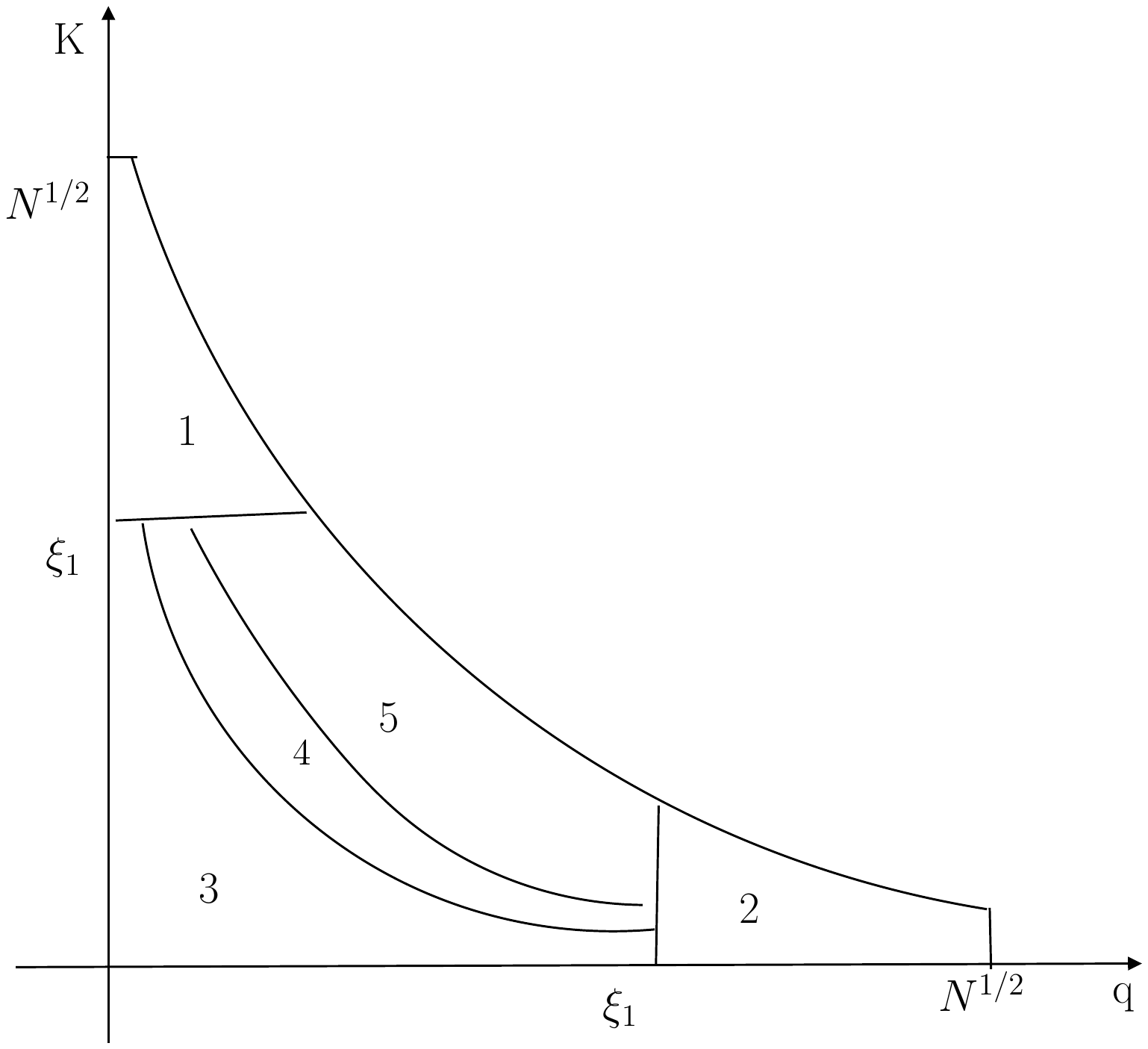}
\end{center}
Lemma \ref{lemma-17-3} corresponds to the domain 1,
Lemma \ref{lemma-17-6} corresponds to the domain 2,
Lemma \ref{lemma-17-7} corresponds to the domain 3,
Lemma \ref{lemma-17-8} corresponds to the domain 4,
Lemma \ref{lemma-17-9} corresponds to the domain 5.


Let
\begin{gather*}
c_1=c_1(N),\,c_2=c_2(N),\,Q_0\le c_1<c_2\le  N^{1/2},\\
f_1=f_1(N,q),\,f_2=f_2(N,q),\,f_1<f_2\le  \frac{\KN}{q},\\
m_1=\min\{f_1(N,N_j),f_1(N,N_{j+1})\},\,m_2=\max\{f_2(N,N_j),f_2(N,N_{j+1})\}.
\end{gather*}
\begin{Le}\label{lemma-17-5}
If the functions $f_1(N,q),f_2(N,q)$ are monotonic for $q,$ then the following inequality holds
\begin{gather}\notag
\mathop{{\sum}^*}_{c_1\le q\le c_2\atop 1\le a\le q }\,\int\limits_{f_1\le|K|\le f_2}
\left|S_N(\frac{a}{q}+\frac{K}{N})\right|^2dK\le\\\le
\sum_{j:\,c_1^{1-\epsilon_0}\le N_j\le c_2\,}\int\limits_{m_1\le|K|\le m_2}
\mathop{{\sum}^*}_{N_j\le q\le N_{j+1}\atop 1\le a\le q }\left|S_N(\frac{a}{q}+\frac{K}{N})\right|^2dK.\label{17-7-15}
\end{gather}
\begin{proof}
The interval $[c_1,c_2]$ can be covered  by the intervals $[N_j,N_{j+1}],$ then
\begin{gather*}
\sum_{c_1\le q\le c_2}\le\sum_{j:\,c_1^{1-\epsilon_0}\le N_j\le c_2}\sum_{N_j\le q\le N_{j+1}}.
\end{gather*}
Interchanging the order of summation over $q$ and integration over $K,$ we obtain \eqref{17-7-15}. This completes the proof of the lemma.
\end{proof}
\end{Le}
Using \eqref{12-3-1} we in the same way as in Lemma \ref{lemma-17-5} obtain the following statement.
\begin{Le}\label{lemma-17-5-2}
If the functions $f_1(N,q),f_2(N,q)$ are monotonic for $q$ and и $m_1\ge Q_0,$ then the following inequality holds
\begin{gather}\notag
\mathop{{\sum}^*}_{c_1\le q\le c_2\atop 1\le a\le q }\,\int\limits_{f_1\le|K|\le f_2}
\left|S_N(\frac{a}{q}+\frac{K}{N})\right|^2dK\ll\\\ll\frac{1}{T}
\sum_{j:\,c_1^{1-\epsilon_0}\le N_j\le c_2\,}
\sum_{i:\,m_1^{1-\epsilon_0}\le N_i\le m_2\,}
\mathop{{\sum}^*}_{N_j\le q\le N_{j+1}\atop 1\le a\le q }
\sum_{TN_i\le|l|\le TN_{i+1}}
\left|S_N(\frac{a}{q}+\frac{l}{TN})\right|^2.\label{17-7-15-2}
\end{gather}
\end{Le}
To simplify we denote
$$Q_1=N_j,\;Q=N_{j+1},\;\Ll_1=TN_i,\;\Ll=TN_{i+1},\;\LL_1=N_i,\;\LL=N_{i+1}.$$
Using the relation of Lemma \ref{lemma10-2}, we obtain
\begin{gather*}
\frac{Q}{Q_1}\le Q^{\epsilon_0}\le N^{\epsilon_0/2},\quad c_1^{1-\epsilon_0}\le Q_1\le c_2,\quad c_1\le Q\le c_2^{1+2\epsilon_0},\\
\frac{\LL}{\LL_1}\le \LL^{\epsilon_0}\le N^{\epsilon_0/2},\quad m_1^{1-\epsilon_0}\le\LL_1\le m_2,\quad m_1\le\LL\le m_2^{1+2\epsilon_0}
\end{gather*}
Further, we will use this bounds without reference to them.
\begin{Le}\label{lemma-17-3}
For $\epsilon_0\in(0,\frac{1}{2500})$ the following inequality holds
\begin{gather}\label{17-7-11}
\frac{1}{N}\mathop{{\sum}^*}_{1\le a\le q\le\QN}\int\limits_{\xi_1\le|K|\le\frac{\KN}{q}}
\left|S_N(\frac{a}{q}+\frac{K}{N})\right|^2dK\ll\frac{|\Omega_N|^2}{N}.
\end{gather}
\begin{proof}
We denote by $I$ the integral in the left side of \eqref{17-7-11}.  Applying Lemma \ref{lemma-15-3}, we obtain
\begin{gather}\label{17-7-12}
I\ll|\Omega_N|^2\frac{N^{2\gamma+2\epsilon_0}}{q^{2-4\epsilon_0}}
\int\limits_{\xi_1\le|K|\le\frac{\KN}{q}}\frac{dK}{\KK^{2-4\epsilon_0}}\ll
|\Omega_N|^2\frac{N^{2\gamma+2\epsilon_0}}{q^{2-4\epsilon_0}\xi_1^{1-4\epsilon_0}}.
\end{gather}
Substituting \eqref{17-7-12} into the left side of \eqref{17-7-11}, we have
\begin{gather}\label{17-7-13}
\frac{1}{N}\mathop{{\sum}^*}_{1\le a\le q\le\QN}\int\limits_{\xi_1\le|K|\le\frac{\KN}{q}}
\left|S_N(\frac{a}{q}+\frac{K}{N})\right|^2dK\ll
\frac{|\Omega_N|^2}{N}
\sum_{1\le q\le N_{1/2}}\frac{N^{2\gamma+2\epsilon_0}}{q^{1-4\epsilon_0}\xi_1^{1-4\epsilon_0}}.
\end{gather}
By the choice of the parameter $\xi_1,$ we obtain
\begin{gather}\label{17-7-14}
\sum_{1\le q\le N_{1/2}}\frac{N^{2\gamma+2\epsilon_0}}{q^{1-4\epsilon_0}\xi_1^{1-4\epsilon_0}}\ll
\frac{N^{2\gamma+4\epsilon_0}}{\xi_1^{1-4\epsilon_0}}\ll N^{-\epsilon_0/2}\ll1.
\end{gather}
Substituting \eqref{17-7-14} into \eqref{17-7-13}, we obtain \eqref{17-7-11}. This completes the proof of the lemma.
\end{proof}
\end{Le}
\begin{Le}\label{lemma-17-6}
For $\gamma\le\frac{1}{5}-4\epsilon_0,\,\epsilon_0\in(0,\frac{1}{2500})$ the following inequality holds
\begin{gather}\label{17-8}
\frac{1}{N}\mathop{{\sum}^*}_{1\le a\le q\le\QN\atop q>\xi_1}
\int\limits_{|K|\le\frac{\KN}{q}}
\left|S_N(\frac{a}{q}+\frac{K}{N})\right|^2dK\ll\frac{|\Omega_N|^2}{N}.
\end{gather}
\begin{proof}
We use Lemma \ref{lemma-17-5} with
\begin{gather*}
c_1=\xi_1,\,c_2=\QN,
f_1=0,\,f_2=\frac{\KN}{q},\,
m_1=\frac{Q_0}{Q},\,m_2=\frac{\KN}{Q_1}.
\end{gather*}
We note that $KQ\le\KN\frac{Q}{Q_1}\le N^{1/2+\epsilon_0/2},$ thus, applying Corollary \ref{sled15-1}, we obtain
\begin{gather}\label{17-10}
\mathop{{\sum}^*}_{N_j\le q\le N_{j+1}\atop 1\le a\le q }\left|S_N(\frac{a}{q}+\frac{K}{N})\right|^2\ll
|\Omega_N|^2\left(
\frac{N^{2\gamma+8\epsilon_0}}{\KK Q_1}+
\frac{N^{5/2\gamma-1/2+9\epsilon_0}}{\KK}+N^{2\gamma-1/2+5\epsilon_0}\right).
\end{gather}
Integrating over $K$ we have
\begin{gather*}
\int\limits_{m_1\le|K|\le m_2}\mathop{{\sum}^*}_{N_j\le q\le N_{j+1}\atop 1\le a\le q }
\left|S_N(\frac{a}{q}+\frac{K}{N})\right|^2dK\ll
|\Omega_N|^2\Bigl(
\frac{N^{2\gamma+9\epsilon_0}}{Q_1}+N^{5/2\gamma-1/2+10\epsilon_0}+
\frac{N^{2\gamma+6\epsilon_0}}{Q_1}\Bigr).
\end{gather*}
Since $Q_1\ge\xi_1^{1-2\epsilon_0}$ and $\gamma<\frac{1}{5}-4\epsilon_0,$ we obtain
\begin{gather}\label{17-11}
\int\limits_{m_1\le|K|\le m_2}\mathop{{\sum}^*}_{N_j\le q\le N_{j+1}\atop 1\le a\le q }
\left|S_N(\frac{a}{q}+\frac{K}{N})\right|^2dK\ll
|\Omega_N|^2 N^{-0,1\epsilon_0}.
\end{gather}
Since the number of summands in the sum over $j$ is less than $\log\log N,$ one has
\begin{gather*}
\sum_{j:\,c_1^{1-\epsilon_0}\le N_j\le c_2\,}\int\limits_{m_1\le|K|\le m_2}
\mathop{{\sum}^*}_{N_j\le q\le N_{j+1}\atop 1\le a\le q }\left|S_N(\frac{a}{q}+\frac{K}{N})\right|^2dK\ll
|\Omega_N|^2 .
\end{gather*}
This completes the proof of the lemma.
\end{proof}
\end{Le}
\begin{Le}\label{lemma-17-7}
For $\gamma\le\frac{1}{6}-10\epsilon_0,\,\epsilon_0\in(0,\frac{1}{2500})$  the following inequality holds
\begin{gather}\label{17-12}
\frac{1}{N}\mathop{{\sum}^*}_{1\le a\le q\le\xi_1\atop q>Q_0}
\int\limits_{|K|\le\frac{N^{2/5-2\epsilon_0}}{q^{6/5}}}
\left|S_N(\frac{a}{q}+\frac{K}{N})\right|^2dK\ll\frac{|\Omega_N|^2}{N}.
\end{gather}
\begin{proof}
We use Lemma \ref{lemma-17-5} with
\begin{gather*}
c_1=Q_0,\,c_2=\xi_1,\,
f_1=0,\,f_2=\frac{N^{2/5-2\epsilon_0}}{q^{6/5}},\,
m_1=\frac{Q_0}{Q},\,m_2=\frac{N^{2/5-2\epsilon_0}}{Q_1^{6/5}}.
\end{gather*}
Applying Lemma \ref{lemma-15-2-2} one has
\begin{gather}\label{17-13}
\mathop{{\sum}^*}_{N_j\le q\le N_{j+1}\atop 1\le a\le q }\left|S_N(\frac{a}{q}+\frac{K}{N})\right|^2\ll
|\Omega_N|^2
\KK^{5\gamma-2+24\epsilon_0}Q^{6\gamma-1+30\epsilon_0}\frac{Q^2}{Q^2_1}.
\end{gather}
Integrating over  $K$ we have
\begin{gather*}
\int\limits_{m_1\le|K|\le m_2}\mathop{{\sum}^*}_{N_j\le q\le N_{j+1}\atop 1\le a\le q }
\left|S_N(\frac{a}{q}+\frac{K}{N})\right|^2dK\ll
|\Omega_N|^2Q^{6\gamma-1+30\epsilon_0}\frac{Q^2}{Q^2_1}\ll|\Omega_N|^2Q^{6\gamma-1+32\epsilon_0}.
\end{gather*}
For the sum over $j$ to be bounded by a constant, it is sufficient to have $\gamma<\frac{1}{6}-10\epsilon_0.$ Lemma is proved.
\end{proof}
\end{Le}
\begin{Le}\label{lemma-17-8}
For $\gamma\le\frac{1}{6}-6\epsilon_0,\,\epsilon_0\in(0,\frac{1}{2500})$ the following inequality holds
\begin{gather}\label{17-14}
\frac{1}{N}\mathop{{\sum}^*}_{1\le a\le q\le\xi_1\atop q>N^{(1-5\gamma)/3-9\epsilon_0}}
\int\limits_{\frac{N^{2/5-2\epsilon_0}}{q^{6/5}}\le|K|\le\frac{N^{2/3-3\epsilon_0}}{q^{2}}}
\left|S_N(\frac{a}{q}+\frac{K}{N})\right|^2dK\ll\frac{|\Omega_N|^2}{N}.
\end{gather}
\begin{proof}
We use Lemma \ref{lemma-17-5-2} with
\begin{gather*}
c_1=N^{(1-5\gamma)/3-9\epsilon_0},\,c_2=\xi_1,
f_1=\frac{N^{2/5-2\epsilon_0}}{q^{6/5}},\,f_2=\frac{N^{2/3-3\epsilon_0}}{q^{2}},\,
m_1=\frac{N^{2/5-2\epsilon_0}}{Q^{6/5}},\,m_2=\frac{N^{2/3-3\epsilon_0}}{Q_1^{2}}.
\end{gather*}
Applying Lemma \ref{lemma-16-2-2} we obtain
\begin{gather*}
\mathop{{\sum}^*}_{N_j\le q\le N_{j+1}\atop 1\le a\le q }
\sum_{TN_i\le|l|\le TN_{i+1}}
\left|S_N(\frac{a}{q}+\frac{l}{TN})\right|^2\ll|\Omega_N|^2T
\LL^{3\gamma-1+20\epsilon_0}Q^{6\gamma-1+24\epsilon_0}\frac{Q^2\LL^{\delta}}{Q_1^2\LL_1^{\delta}}\ll\\\ll
|\Omega_N|^2
\LL^{3\gamma-1+22\epsilon_0}Q^{6\gamma-1+26\epsilon_0}.
\end{gather*}
By Lemma \ref{lemma-17-5-2} we have
\begin{gather*}
\mathop{{\sum}^*}_{c_1\le q\le c_2\atop 1\le a\le q }\,\int\limits_{f_1\le|K|\le f_2}
\left|S_N(\frac{a}{q}+\frac{K}{N})\right|^2dK\ll\\\ll\frac{1}{T}
\sum_{j:\,c_1^{1-\epsilon_0}\le N_j\le c_2\,}
\sum_{i:\,m_1^{1-\epsilon_0}\le N_i\le m_2\,}
\mathop{{\sum}^*}_{N_j\le q\le N_{j+1}\atop 1\le a\le q }
\sum_{TN_i\le|l|\le TN_{i+1}}
\left|S_N(\frac{a}{q}+\frac{l}{TN})\right|^2\ll\\\ll|\Omega_N|^2
\sum_{j:\,c_1^{1-\epsilon_0}\le N_j\le c_2\,}
\sum_{i:\,m_1^{1-\epsilon_0}\le N_i\le m_2\,}
\LL^{3\gamma-1+22\epsilon_0}Q^{6\gamma-1+26\epsilon_0}.
\end{gather*}
For the sums over $i$ and $j$ to be bounded by a constant, it is sufficient to have $\gamma<\frac{1}{6}-6\epsilon_0.$ Lemma is proved.
\end{proof}
\end{Le}
\begin{Le}\label{lemma-17-9}
For $\gamma\le\frac{1}{6}-10\epsilon_0,\,\epsilon_0\in(0,\frac{1}{2500})$ the following inequality holds
\begin{gather}\label{17-15}
\frac{1}{N}\mathop{{\sum}^*}_{1\le a\le q\le\xi_1\atop q>N^{(1-3\gamma)/3-7\epsilon_0}}
\int\limits_{\frac{N^{2/3-3\epsilon_0}}{q^{2}}\le|K|\le\frac{\KN}{q}}
\left|S_N(\frac{a}{q}+\frac{K}{N})\right|^2dK\ll\frac{|\Omega_N|^2}{N}.
\end{gather}
\begin{proof}
We use Lemma \ref{lemma-17-5-2} with
\begin{gather*}
c_1=N^{(1-3\gamma)/3-7\epsilon_0},\,c_2=\xi_1,
f_1=\frac{N^{2/3-3\epsilon_0}}{q^{2}},\,f_2=\frac{\KN}{q},\,
m_1=\frac{N^{2/3-3\epsilon_0}}{Q^{2}},\,m_2=\frac{\KN}{Q_1}.
\end{gather*}
Applying Lemma \ref{lemma-16-7} we obtain
\begin{gather*}
\mathop{{\sum}^*}_{N_j\le q\le N_{j+1}\atop 1\le a\le q }
\sum_{TN_i\le|l|\le TN_{i+1}}
\left|S_N(\frac{a}{q}+\frac{l}{TN})\right|^2\ll|\Omega_N|^2T
\Bigl(
\frac{N^{2\gamma+4\epsilon_0}(\LL Q)^{4\epsilon_0}}{\LL_1Q_1}+
N^{2\gamma-1+4\epsilon_0}(\LL Q)^{1+4\epsilon_0}+\\+
N^{3\gamma-1/2+8\epsilon_0}\frac{(\LL Q)^{1/2-\gamma+6\epsilon_0}Q^{2\epsilon_0}}{Q_1}+
N^{3\gamma-3/2+8\epsilon_0}(\LL Q)^{3/2-\gamma+6\epsilon_0}Q^{1+2\epsilon_0}
\Bigr).
\end{gather*}
By Lemma \ref{lemma-17-5-2} we have
\begin{gather*}
\mathop{{\sum}^*}_{c_1\le q\le c_2\atop 1\le a\le q }\,\int\limits_{f_1\le|K|\le f_2}
\left|S_N(\frac{a}{q}+\frac{K}{N})\right|^2dK\ll\\\ll\frac{1}{T}
\sum_{j:\,c_1^{1-\epsilon_0}\le N_j\le c_2\,}
\sum_{i:\,m_1^{1-\epsilon_0}\le N_i\le m_2\,}
\mathop{{\sum}^*}_{N_j\le q\le N_{j+1}\atop 1\le a\le q }
\sum_{TN_i\le|l|\le TN_{i+1}}
\left|S_N(\frac{a}{q}+\frac{l}{TN})\right|^2\ll\\\ll|\Omega_N|^2
\sum_{j:\,c_1^{1-\epsilon_0}\le N_j\le c_2\,}
\sum_{i:\,m_1^{1-\epsilon_0}\le N_i\le m_2\,}
\Bigl(
\frac{N^{2\gamma+7\epsilon_0}}{\LL Q}+
N^{2\gamma-1+4\epsilon_0}(\LL Q)^{1+4\epsilon_0}+\\+
N^{3\gamma-1/2+8\epsilon_0}\frac{(\LL Q)^{1/2-\gamma+6\epsilon_0}Q^{3\epsilon_0}}{Q}+
N^{3\gamma-3/2+8\epsilon_0}(\LL Q)^{3/2-\gamma+6\epsilon_0}Q^{1+2\epsilon_0}
\Bigr).
\end{gather*}
Using the following estimates (we recall that $\LL=N_{i+1}$)
\begin{gather}\label{17-16}
\sum_{i:\,m_1^{1-\epsilon_0}\le N_i\le m_2\,}\LL^{\alpha}\ll
\left\{
              \begin{array}{ll}
                N^{\epsilon_0}m_2^{\alpha}, & \hbox{if $\alpha>0$;} \\
                N^{\epsilon_0}m_1^{\alpha}, & \hbox{if $\alpha<0$.}
              \end{array}
\right.
\end{gather}
we obtain
\begin{gather*}
\mathop{{\sum}^*}_{c_1\le q\le c_2\atop 1\le a\le q }\,\int\limits_{f_1\le|K|\le f_2}
\left|S_N(\frac{a}{q}+\frac{K}{N})\right|^2dK\ll\\\ll|\Omega_N|^2
\sum_{j:\,c_1^{1-\epsilon_0}\le N_j\le c_2\,}
\Bigl(
N^{2\gamma-2/3+12\epsilon_0}Q+
N^{2\gamma-1/2+14\epsilon_0}+
N^{5\gamma/2-1/4+14\epsilon_0}Q^{-1+10\epsilon_0}+\\+
N^{5\gamma/2-3/4+14\epsilon_0}Q^{1+6\epsilon_0}
\Bigr)
\ll|\Omega_N|^2\Bigl(
N^{2\gamma-2/3+12\epsilon_0}N^{2\gamma+12\epsilon_0}+\\+
N^{2\gamma-1/2+15\epsilon_0}+
N^{5\gamma/2-1/4+14\epsilon_0}N^{-(1-3\gamma)/3+7\epsilon_0}+
N^{5\gamma/2-3/4+14\epsilon_0}N^{2\gamma+14\epsilon_0}\Bigr)\ll|\Omega_N|^2N^{-\epsilon_0}.
\end{gather*}
The last inequality holds if $\gamma<\frac{1}{6}-10\epsilon_0.$ Lemma is proved.
\end{proof}
\end{Le}


\section{The proof of Theorem \ref{uslov}.}
Let $\gamma<\frac{1}{6}.$ We choose $\epsilon_0$ such that $\epsilon_0\in(0,\frac{1}{2500})$ and
$\gamma\le\frac{1}{6}-10\epsilon_0.$ Then it follows from Lemma \ref{lemma-17-3}~--\ref{lemma-17-9} that the first integral in the right side of \eqref{17-4} is less than $\frac{|\Omega_N|^2}{N},$ that is,
\begin{gather}\label{17-30}
\frac{1}{N}\mathop{{\sum}^*}_{0\le a\le q\le\QN\atop q>Q_0}\int\limits_{\frac{Q_0}{q}\le|K|\le\frac{\KN}{q}}
\left|S_N(\frac{a}{q}+\frac{K}{N})\right|^2dK\ll\frac{|\Omega_N|^2}{N}.
\end{gather}.
Substituting \eqref{17-30} and the results of Lemma \ref{lemma-17-3-2} in Lemma \ref{lemma-17-2}, we obtain
\begin{gather}\label{17-31}
\int_0^1\left|S_N(\theta)\right|^2d\theta\ll\frac{|\Omega_N|^2}{N}\quad\mbox{при}\quad
\gamma<\frac{1}{6}.
\end{gather}
Thus the inequality \eqref{8-7} is proved. This, as it was proved in the section §\ref{section idea BK}, is enough for proving Theorem \ref{uslov}. This completes the proof of the theorem.


\textit{I.D. Kan}\\
\textit{Moscow State university}\\
\textit{e-mail address: igor.kan@list.ru}\\
\\
\textit{D.A. Frolenkov}\\
\textit{Steklov Mathematical Institute}\\
\textit{Russia, Moscow, Gubkina str., 8, 119991,}\\
\textit{e-mail address: frolenkov\underline{ }adv@mail.ru}
\end{document}